\documentclass[reqno,10pt,centertags,draft]{amsart}
\usepackage{amsmath,amsthm,amscd,amssymb,latexsym,upref}
\date{\today}

\newcommand{\bbN}{{\mathbb{N}}}
\newcommand{\bbR}{{\mathbb{R}}}

\newcommand{\bbZ}{{\mathbb{Z}}}
\newcommand{\bbC}{{\mathbb{C}}}
\newcommand{\bbQ}{{\mathbb{Q}}}
\newcommand{\bbT}{{\mathbb{T}}}

\newcommand{\cD}{{\mathcal D}}
\newcommand{\cE}{{\mathcal E}}

\newcommand{\cH}{{\mathcal H}}

\newcommand{\cL}{{\mathcal L}}
\newcommand{\cM}{{\mathcal M}}

\newcommand{\cR}{{\mathcal R}}

\newcommand{\gE}{\mathfrak{E}}

\newcommand{\no}{\notag}
\newcommand{\lb}{\label}
\newcommand{\f}{\frac}

\newcommand{\ol}{\overline}
\newcommand{\ti}{\tilde}
\newcommand{\wti}{\widetilde}

\newcommand{\Oh}{O}
\newcommand{\oh}{o}

\newcommand{\loc}{\text{\rm{loc}}}

\newcommand{\Arg}{\text{\rm{Arg}}}

\newcommand{\dom}{\text{\rm{dom}}}

\newcommand{\supp}{\text{\rm{supp}}}
\newcommand{\dist}{\text{\rm{dist}}}
\newcommand{\diam}{\text{\rm{diam}}}

\newcommand{\bi}{\bibitem}
\newcommand{\hatt}{\widehat}
\newcommand{\beq}{\begin{equation}}
\newcommand{\eeq}{\end{equation}}
\newcommand{\ba}{\begin{align}}
\newcommand{\ea}{\end{align}}

\newcommand{\tr}{\text{\rm{tr}}}

\newcommand{\abs}[1]{\lvert#1\rvert}

\renewcommand{\Re}{\text{\rm Re}}
\renewcommand{\Im}{\text{\rm Im}}
\renewcommand{\ln}{\text{\rm ln}}
\renewcommand{\ge}{\geqslant}
\renewcommand{\le}{\leqslant}

\DeclareMathOperator{\SY}{SY(\cE)}
\DeclareMathOperator{\AP}{AP(\bbR)}

\allowdisplaybreaks
\numberwithin{equation}{section}

\newtheorem{theorem}{Theorem}[section]

\newtheorem{lemma}[theorem]{Lemma}

\theoremstyle{definition}
\newtheorem{definition}[theorem]{Definition}
\newtheorem{remark}[theorem]{Remark}

\newtheorem{example}[theorem]{Example}

\begin{document}

\title[Spectral Properties of reflectionless
Schr\"odinger operators]{Spectral Properties of a Class
of Reflectionless Schr\"odinger Operators}
\author[F.\ Gesztesy and P.\ Yuditskii]{Fritz Gesztesy and Peter
Yuditskii}
\address{Department of Mathematics,
University of Missouri, Columbia, MO 65211, USA}
\email{fritz@math.missouri.edu}
\urladdr{http://www.math.missouri.edu/personnel/faculty/gesztesyf.html}
\address{Department of Mathematics, Bar Ilan University, 52900 Ramat Gan, 
Israel}
\email{yuditski@macs.biu.ac.il}
\thanks{Based upon work partially supported by the US National Science
Foundation under Grant No.\ DMS-0405526 and the Austrian Science Fund FWF, 
project number: P16390-N04.}
\date{\today}
\subjclass[2000]{Primary 34B20, 34L05, 34L40; Secondary 34B24, 34B27, 47A10.}

\begin{abstract}
We prove that one-dimensional reflectionless Schr\"odinger
operators with spectrum a homogeneous set in the sense of Carleson,
belonging to the class introduced by Sodin and
Yuditskii, have purely absolutely continuous spectra. This class includes
all earlier examples of reflectionless almost periodic Schr\"odinger
operators.

In addition, we construct examples of reflectionless Schr\"odinger
operators with more general types of spectra, given by the complement of a
Denjoy--Widom-type domain in $\bbC$, which exhibit a singular component.
\end{abstract}

\maketitle

\section{Introduction}\label{s1}

In this paper we consider a certain class of reflectionless
self-adjoint Schr\"odinger operators $H=-d^2/dx^2+V$ in $L^2(\bbR;dx)$ and
study the spectral properties of its members. In particular, assuming that
the spectrum $\sigma(H)$ of $H$ is given by a homogeneous set $\cE$ in the
sense of Carleson, we prove that $H$ has purely absolutely continuous
spectrum. Subsequently, we consider reflectionless Schr\"odinger operators
with more general types of spectra given by the complement of a
Denjoy--Widom-type domain in $\bbC$  and construct a class of examples that
exhibits a singular component in its spectrum.

To put our results in proper perspective we first recall that a {\it
homogeneous set} $\cE\subset\bbR$, bounded from below, is of the type (cf.\
Carleson \cite{Ca83})
\begin{equation}
\cE= [E_0,\infty)\big\backslash\bigcup_{j\in J} (a_j,b_j), \quad
J\subseteq\bbN, \lb{1.1}
\end{equation}
for some $E_0\in\bbR$ and $a_j<b_j$, where $(a_j,b_j)\cap
(a_{j'},b_{j'})=\emptyset$, $j,j'\in J$, $j\neq j'$ such that the following
condition is fulfilled:
\begin{align}
\begin{split}
& \text{There exists an $\varepsilon>0$
such that for all $\lambda\in\cE$}  \lb{1.2} \\
& \text{and all $\delta>0$, \,
$|\cE\cap(\lambda-\delta,\lambda+\delta)|\geq \varepsilon\delta$.}
\end{split}
\end{align}
We will also assume that $\cE$ has finite gap length, that is, $\sum_{j\in J}
(b_j-a_j) < \infty$ is supposed to hold.

{\it Reflectionless} (self-adjoint) Schr\"odinger operators $H$ in
$L^2(\bbR;dx)$ can be characterized by the fact that for all $x\in\bbR$ and
for a.e.\
$\lambda\in \sigma_{\rm ess}(H)$, the diagonal Green's function of $H$ has
purely imaginary normal boundary values,
\begin{equation}
G(\lambda+i0,x,x) \in i\bbR.    \lb{1.3}
\end{equation}
Here $\sigma_{\rm ess}(H)$ denotes the essential spectrum of $H$ (we assume
$\sigma_{\rm ess}(H)\neq\emptyset$) and
\begin{equation}
G(z,x,x')=(H-zI)^{-1}(x,x'), \quad z\in\bbC\backslash\sigma(H), \lb{1.4}
\end{equation}
denotes the integral kernel of the resolvent of $H$.

$H_0=-d^2/dx^2$ and the $N$-soliton potentials $V_N$, $N\in\bbN$,
that is,  exponentially decreasing solutions in $C^\infty(\bbR)$ of some (and
hence infinitely many) equations of the stationary Korteweg--de Vries (KdV)
hierarchy, yield well-known examples of reflectionless Schr\"odinger
operators $H_N=-d^2/dx^2+V_N$. However, since such operators have $N\in\bbN$
strictly negative eigenvalues, the spectrum of $H_N$, $N\in\bbN$, is not a
homogeneous set. The prototype of reflectionless Schr\"odinger operators
with a homogeneous spectrum are in fact the set of periodic Schr\"odinger
operators. Indeed, if $V_a$ is periodic with some period $a>0$,
that is, $V_a(x+a)=V_a(x)$ for a.e.\ $x\in\bbR$, then standard Floquet
theoretic considerations show that the spectrum of $H_a=-d^2/dx^2+V_a$ is a
countable union of compact intervals (which may degenerate into a union of
finitely-many compact intervals and a half-line) and hence a homogeneous set,
and at the same time, the diagonal Green's function of $H_a$ is purely
imaginary for every point in the open interior of $\sigma(H_a)$. More
generally, also certain classes of quasi-periodic and almost periodic
potentials give rise to reflectionless Schr\"odinger operators with
homogeneous spectra. The prime example of such quasi-periodic potentials is
represented by the class of real-valued bounded algebro-geometric KdV
potentials corresponding to an underlying (compact) hyperelliptic Riemann
surface (see, e.g., \cite[Ch.\ 3]{BBEIM94}, \cite{DMN76},
\cite[Ch.\ 1]{GH03}, \cite{Jo88}, \cite[Chs.\ 8, 10]{Le87},
\cite[Ch.\ 4]{Ma86},
\cite[Ch.\ II]{NMPZ84} and the literature cited therein). More general
classes of almost periodic and reflectionless Schr\"odinger operators were
studied by  Avron and Simon \cite{AS81}, Carmona and Lacroix \cite[Ch.\
VII]{CL90}, Chulaevskii \cite{Ch84}, Craig \cite{Cr89},
Deift and Simon \cite{DS83},  Egorova \cite{Eg92}, Johnson \cite{Jo82},
Johnson and Moser \cite{JM82}, Kotani \cite{Ko84}--\cite{Ko87b},
Kotani and Krishna \cite{KK88}, Levitan \cite{Le82}--\cite{Le85},
\cite[Chs.\ 9, 11]{Le87},  Levitan and Savin \cite{LS84}, Moser \cite{Mo81},
Pastur and Figotin \cite[Chs.\ V, VII]{PF92}, Pastur and Tkachenko
\cite{PT89}, and more recently, by Sodin and Yuditskii
\cite{SY95a}--\cite{SY96}, the starting point of our investigation.

In Section \ref{s2} we introduce the class $\SY$ of potentials
associated with one-dimensional Schr\"odinger operators $H=-d^2/dx^2+V$,
whose spectra coincide with the  homogeneous set $\cE$. The class $\SY$ is
possibly a slight generalization of the class of potentials studied by Sodin
and Yuditskii \cite{SY95}, \cite{SY96} in the sense that we do not assume
continuity of the potential $V$ from the outset. For precise details on
$\SY$ we refer to Definition \ref{d2.3}. Here we just mention that $V\in\SY$
if $\sigma(H)=\cE$ and the half-line Weyl--Titchmarsh functions
$m_{\pm}(z,x_0)$ (associated with the restriction of $H$ to the half-lines
$(x_0,\pm\infty)$) satisfy a certain pseudo-analytic continuation property
across $\cE$. We then derive the essential boundedness
of such potentials using the trace formula proved in \cite{GS96},
\begin{equation}
\SY\subset L^\infty(\bbR;dx),  \lb{1.5}
\end{equation}
and set the stage for the remainder of this paper by showing that the
corresponding Schr\"odinger operators are reflectionless.

In Section \lb{s3} we prove by elementary methods that the absolutely
continuous spectrum of $H$ with $V\in\SY$ coincides with the prescribed
homogeneous set $\cE$,
\begin{equation}
\sigma_{\rm ac}(H)=\cE,  \lb{1.6}
\end{equation}
and that $H$ has uniform multiplicity equal to two. To obtain this result we
exploit the notion of essential supports of measures and the essential
closure of sets, a topic summarized in Appendix \ref{A}.

In Section \ref{s4} we prove the absence of a singular component in the
spectrum of $H$ and together with
\eqref{1.6} this then implies that
\begin{equation}
\sigma(H)=\sigma_{\rm ac}(H)=\cE, \quad \sigma_{\rm sc}(H)
=\sigma_{\rm pp}(H)=\emptyset,  \lb{1.7}
\end{equation}
our principal result for potentials $V$ in the class $\SY$. The proof of
$\sigma_{\rm sc}(H)=\sigma_{\rm pp}(H)=\emptyset$ is more involved and
requires some results from the theory of Hardy spaces, most notably,
the space $H^1(\cE)$. The necessary prerequisites for
this topic are summarized in Appendix \ref{B}.

Finally, in Section \ref{s5} we construct a class of reflectionless
Schr\"odinger operators with more general types of spectra given by the
complement of a Denjoy--Widom-type domain in $\bbC$. In particular, we
construct a class of examples which exhibits a particular accumulation point
of spectral bands and gaps. As a result, $H$ acquires an eigenvalue and
hence a singular component in its spectrum.

Throughout the bulk of this paper we repeatedly exploit properties of
Herglotz functions and their exponential representations as well as
certain elements of Weyl--Titchmarsh and spectral multiplicity theory
for self-adjoint Schr\"odinger operators on half-lines and on $\bbR$. A
nutshell-type treatment of both topics is presented in Appendix \ref{B}.

While we focus in this paper on one-dimensional Schr\"odinger operators,
we emphasize that all methods employed apply to Jacobi, Dirac, and CMV-type
operators and yield an analogous set of results.

\section{The Sodin--Yuditskii Class of Reflectionless Potentials}
\label{s2}

In this section we describe the Sodin--Yuditskii
class $\SY$ of reflectionless potentials associated with the homogeneous set
$\cE$ bounded from below and recall some of the basic properties of
one-dimensional Schr\"odinger operators with such potentials.

Before we analyze the class $\SY$ in some detail, we start with some
general considerations of one-dimensional Schr\"odinger operators. Let
\begin{equation}
V\in L^1_{\loc}(\bbR;dx), \quad V \, \text{ real-valued,}  \lb{2.0}
\end{equation}
and assume that the
differential expression
\begin{equation}
L=-d^2/dx^2+V(x), \quad x\in\bbR   \lb{2.0a}
\end{equation}
is in the limit point case at
$+\infty$ and $-\infty$. We denote by $H$ the corresponding self-adjoint
realization of $L$ in $L^2(\bbR;dx)$ given by
\begin{align}
\begin{split}
& Hf=L f, \\
& f\in\dom(H)=\big\{g\in L^2(\bbR;dx)\,|\, g,g' \in AC_{\loc}(\bbR); \,
L g \in L^2(\bbR;dx)\big\}.  \lb{2.1}
\end{split}
\end{align}
Let $g(z,\cdot)$ denote the diagonal Green's function of $H$,
that is,
\begin{equation}
g(z,x)=G(z,x,x), \quad G(z,x,x')=(H-zI)^{-1}(x,x'), \quad
z\in\bbC\backslash\sigma(H), \; x,x'\in\bbR.  \lb{2.2}
\end{equation}
Since for each $x\in\bbR$, $g(\cdot,x)$ is a Herglotz function (i.e., it
maps the open complex upper half-plane analytically to itself),
\begin{equation}
\xi(\lambda,x)=\f{1}{\pi}\lim_{\varepsilon\downarrow 0}
\Im[\ln(g(\lambda+i\varepsilon,x))] \, \text{ for a.e.\
$\lambda\in\bbR$}  \lb{2.3}
\end{equation}
is well-defined for each $x\in\bbR$. In particular, for all $x\in\bbR$,
\begin{equation}
0 \leq \xi(\lambda,x) \leq 1 \,
\text{ for a.e.\ $\lambda\in\bbR$.}   \lb{2.4}
\end{equation}

In the following we will frequently use the convenient abbreviation
\begin{equation}
h(\lambda_0+i0)=\lim_{\varepsilon\downarrow 0}
h(\lambda_0 +i\varepsilon), \quad \lambda_0\in\bbR, \lb{2.4a}
\end{equation}
whenever the limit in \eqref{2.4a} is well-defined and hence \eqref{2.3}
can then be written as $\xi(\lambda,x)=(1/\pi)\Arg(g(\lambda+i0,x))$.
Moreover, we will use the convention that whenever the phrase a.e.\ is used
without further qualification, it always refers to Lebesgue measure on
$\bbR$.

Associated with $H$ in $L^2(\bbR;dx)$ we also introduce the two half-line
Schr\"odinger operators $H_{\pm,x_0}$ in $L^2([x_0,\pm\infty);dx)$ with
Dirchlet boundary conditions at the finite endpoint $x_0\in\bbR$,
\begin{align}
& H_{\pm,x_0} f= Lf, \no \\
& f\in\dom(H_{\pm,x_0})=\big\{g\in L^2([x_0,\pm\infty);dx)\,|\, g, g'\in
AC([x_0, x_0\pm R]) \, \text{for all $R>0$;} \no \\
& \hspace*{4.5cm}
\lim_{\varepsilon\downarrow 0} g(x_0\pm\varepsilon)=0; \, Lg\in
L^2([x_0,\pm\infty);dx)\big\}.  \lb{2.5}
\end{align}
The {\it half-line Weyl--Titchmarsh functions} associated with
$H_{\pm,x_0}$ (to be discussed in \eqref{2.30}--\eqref{2.32}), are
denoted by $m_{\pm}(z, x_0)$, $z\in\bbC\backslash\sigma(H_{\pm,x_0})$.

Next, we recall the trace formula first proved in \cite{GS96} and
slightly refined in \cite{CGHL00} (see also \cite{Ry01a} for an
interesting extension).
\begin{theorem} \lb{t2.1}
Let $H$ be the self-adjoint Schr\"odinger operator defined in \eqref{2.1}
and assume in addition that $H$ is bounded from below,
$E_0=\inf(\sigma(H))>-\infty$. Then,
\begin{equation}
V(x)=E_0 +\lim_{z\to i\infty}\int_{E_0}^\infty
d\lambda\,
z^2(\lambda-z)^{-2}[1-2\xi(\lambda,x)] \, \text{ for a.e.\
$x\in\bbR$.} \lb{2.6}
\end{equation}
\end{theorem}
\begin{proof}
We briefly sketch the main arguments. The exponential Herglotz
representation of $g(z,x)$ yields (cf.\ \eqref{2.14b}, \eqref{B.8m},
\eqref{B.8o})
\begin{align}
g(z,x)&=[m_-(z,x)-m_+(z,x)]^{-1}  \lb{2.7} \\
&=\exp\bigg[c(x) + \int_{\bbR} d\lambda \, \xi(\lambda,x)
\bigg(\f{1}{\lambda-z}-\f{\lambda}{1+\lambda^2}\bigg)\bigg], \lb{2.8} \\
& \hspace*{-.95cm} c(x)=\Re[\ln(g(i,x))], \quad
z\in\bbC\backslash\sigma(H), \; x\in\bbR
\no
\end{align}
and hence by \eqref{2.4},
\begin{align}
\f{d}{dz}\ln[g(z,x)]=\int_{E_0}^\infty d\lambda \, (\lambda-z)^{-2}
\xi(\lambda,x), \quad z\in\bbC\backslash\sigma(H), \; x\in\bbR. \lb{2.9}
\end{align}
Next, we denote by $\cL_V$ the set of Lebesgue points of $V$, that is,
\begin{equation}
\cL_V=\bigg\{x\in\bbR\,\bigg|\, \int_{-\varepsilon}^{\varepsilon} dx'\,
|V(x+x')-V(x)|\underset{\varepsilon\downarrow 0}{=}\oh(\varepsilon)
\text{ as $\varepsilon\downarrow 0$}\bigg\},  \lb{2.10}
\end{equation}
and we suppose that $x\in \cL_V$. Let $C_\varepsilon\subset\bbC_+$ be the
sector  along the positive imaginary axis with vertex at zero and opening
angle $\varepsilon$ with $0<\varepsilon<\pi/2$.  Then, as proved in
\cite[Theorem~4.8]{CG01}, $m_\pm(z,x)$ has an asymptotic expansion of
the form $(\Im(z^{1/2})\ge0$, $z\in\bbC)$
\begin{equation}
m_\pm(z,x)\underset{\substack{\abs{z}\to\infty\\
z\in C_\varepsilon}}{=}\pm i z^{1/2} \mp (i/2) V(x)z^{-1/2}
+o(|z|^{-1/2})    \lb{2.11}
\end{equation}
as $\abs{z}\to\infty$. By \eqref{2.7}, $g(z,x)$ has an asymptotic
expansion in $C_\varepsilon$ of the form
\begin{equation}
g(z,x)\underset{\substack{\abs{z}\to\infty\\ z\in
C_\varepsilon}}{=}
              (i/2) z^{-1/2} + (i/4) V(x) z^{-3/2}+o(|z|^{-3/2}) \lb{2.12}
\end{equation}
as $\abs{z}\to\infty$. Differentiation of \eqref{2.12} with respect to
$z$ then yields
\begin{equation}
-\f{d}{dz}\ln[g(z,x)] \underset{z\to i\infty}{=}
\f{1}{2} z^{-1}+\f12 V(x)z^{-2}+o(|z|^{-2}). \lb{2.13}
\end{equation}
Thus,
\begin{align}
-\f{d}{dz}\ln[g(z,x)] &=\f12 (z-E_0)^{-1}
+\f12 \int_{E_0}^\infty d\lambda\,
(\lambda-z)^{-2}[1-2\xi(\lambda,x)] \no \\
&\underset{z\to i\infty}{=}\f12  z^{-1}
+ \f12 V(x)z^{-2}+o(|z|^{-2})  \lb{2.14}
\end{align}
proves \eqref{2.6} for $x\in \cL_V$. Since $\bbR\backslash \cL_V$ has
Lebesgue measure zero, \eqref{2.6} is proved.
\end{proof}

For subsequent purpose we also note the universal asymptotic $z$-behavior
of $m_{\pm}(z,x)$ and $g(z,x)$ valid for all $x\in\bbR$,
\begin{align}
& m_{\pm}(z,x)\underset{\substack{\abs{z}\to\infty\\ z\in
C_\varepsilon}}{=}
         \pm i z^{1/2}[1 +o(1)], \lb{2.14a}  \\
& g(z,x)\underset{\substack{\abs{z}\to\infty\\ z\in
C_\varepsilon}}{=} (i/2) z^{-1/2}[1+o(1)]. \lb{2.14b}
\end{align}

Since Schr\"odinger operators bounded from below play a special role in
our considerations, we take a closer look at them next. Given $V\in
L^1_{\loc}(\bbR;dx)$, $V$ real-valued, and $L =-d^2/dx^2+V$ as in
\eqref{2.0} and \eqref{2.0a}, we define the associated minimal
Schr\"odinger operator $H_{\rm min}$ in $L^2(\bbR;dx)$ by
\begin{align}
&H_{\rm min}f= L f,  \no \\
&f\in\dom(H_{\rm min})=\big\{g\in L^2(\bbR;dx)\,|\,
g,g' \in AC_{\loc}(\bbR); \, \supp\,(g) \text{ compact;}  \lb{2.19} \\
& \hspace*{8.4cm} L g \in L^2(\bbR;dx)\big\}.   \no
\end{align}
By a well-known result of Hartman \cite{Ha48} (see also Rellich
\cite{Re51}, \cite{Ge93}, and the literature cited in \cite{CG03}),
$H_{\rm min}$ is essentially self-adjoint, or equivalently, the
differential expression $L$ is in the limit point case at $+\infty$
and $-\infty$) if $H_{\rm min}$ is bounded from below. In this case, the
operator $H=\ol{H_{\rm min}}$ (the operator closure of $H_{\rm min}$) is
the unique self-adjoint extension of $H_{\rm min}$ in $L^2(\bbR;dx)$ and
it coincides with the maximally defined operator in \eqref{2.1}.

\begin{definition}  \lb{d2.2}
Let $\cE\subset\bbR$ be a closed set bounded from below which we may write
as
\begin{equation}
\cE= [E_0,\infty)\big\backslash\bigcup_{j\in J} (a_j,b_j), \quad
J\subseteq\bbN, \lb{2.15}
\end{equation}
for some $E_0\in\bbR$ and $a_j<b_j$, where $(a_j,b_j)\cap
(a_{j'},b_{j'})=\emptyset$, $j,j'\in J$, $j\neq j'$. Then $\cE$ is called
{\it homogeneous} if
\begin{align}
\begin{split}
& \text{there exists an $\varepsilon>0$
such that for all $\lambda\in\cE$}  \lb{2.16} \\
& \text{and all $\delta>0$, \,
$|\cE\cap(\lambda-\delta,\lambda+\delta)|\geq \varepsilon\delta$.}
\end{split}
\end{align}
Moreover, we say that $\cE$ is of {\it finite gap length} if
\begin{equation}
\sum_{j\in J} (b_j-a_j) < \infty.  \lb{2.17}
\end{equation}
\end{definition}

Homogeneous sets were originally discussed by Carleson \cite{Ca83}; we
also refer to \cite{JM85}, \cite{PY03}, and \cite{Zi89}.

Next we introduce the Sodin--Yuditskii class $\SY$ of potentials
associated with a homogeneous set $\cE$ of finite gap length.

\begin{definition}  \lb{d2.3}
Let $\cE\subset\bbR$ be a homogeneous set of finite gap length of the form
\eqref{2.15} and pick $x_0\in\bbR$. Then $V\in L^1_{\loc}(\bbR;dx)$ belongs
to the Sodin--Yuditskii class $\SY$ associated with $\cE$ if
\begin{align}
& (i) \; \text{ $V$ is real-valued.}  \lb{2.21} \\
& (ii) \text{ $H_{\rm min}$ is bounded from below and its unique
self-adjoint}  \no \\
& \qquad \text{extension $H$ has spectrum $\sigma(H) = \cE$.} \lb{2.22} \\
& (iii) \; \text{For a.e.\ $\lambda\in\cE$, \,
$\lim_{\varepsilon\downarrow 0} m_+(\lambda+i\varepsilon,x_0)=
\lim_{\varepsilon\downarrow 0} m_-(\lambda-i\varepsilon,x_0)$.}  \lb{2.23}
\end{align}
\end{definition}

Next we will demonstrate that the Sodin--Yuditskii class $\SY$ is independent of the
choice of $x_0\in\bbR$ in \eqref{2.23} and that $\SY\subset L^\infty(\bbR;dx)$. Moreover, potentials $V\in\SY$ are reflectionless in a sense that will be made precise below. These results have short proofs which we will present. Later in this section we will indicate that all potentials in $\SY$ are actually continuous on $\bbR$ and uniformly (i.e., Bohr) almost periodic. Since the latter result requires quite different potential theoretic techniques we will quote them without proofs.   

\begin{theorem} \lb{t2.4}
Let $V\in \SY$. Then, $V$ is reflectionless in the sense that
\begin{equation}
\text{for each $x\in\bbR$, } \, \xi(\lambda,x)=1/2 \, \text{ for
a.e.\ $\lambda\in\cE$.}    \lb{2.24}
\end{equation}
Moreover,
\begin{equation}
\SY \subset L^\infty(\bbR;dx)   \lb{2.25}
\end{equation}
and $\SY$ is independent of the choice of $x_0\in\bbR$ in
Definition \ref{d2.3}\,$(iii)$.
\end{theorem}
\begin{proof}
We introduce the usual fundamendal system of solutions
$\phi(z,\cdot,x_0)$ and $\theta(z,\cdot,x_0)$,
$z\in\bbC$, with respect to a fixed reference point $x_0\in\bbR$, of
solutions of
\begin{equation}
L \psi(z,x) = z \psi(z,x), \quad (z,x)\in \bbC\times\bbR, \lb{2.26}
\end{equation}
satisfying the initial conditions at the point $x=x_0$,
\begin{equation}
\phi(z,x_0,x_0)=\theta'(z,x_0,x_0)=0, \quad
\phi'(z,x_0,x_0)=\theta(z,x_0,x_0)=1. \lb{2.27}
\end{equation}
Then for any fixed $x, x_0\in \bbR$, $\phi(z,x,x_0)$ and
$\theta(z,x,x_0)$ are entire with respect to $z$,
\begin{equation}
W(\theta(z,\cdot,x_0),\phi(z,\cdot,x_0))(x)=1, \quad
z\in\bbC,  \lb{2.28}
\end{equation}
and
\begin{equation}
\ol{\phi(z,x,x_0)}=\phi(\ol z,x,x_0), \quad
\ol{\theta(z,x,x_0)}=\theta(\ol z,x,x_0), \quad
(z,x,x_0)\in\bbC\times\bbR^2. \lb{2.29}
\end{equation}
Particularly important solutions of \eqref{2.26} are the so called
{\it Weyl--Titchmarsh solutions} $\psi_{\pm}(z,\cdot,x_0)$,
$z\in\bbC\backslash\sigma(H)$, uniquely characterized by
\begin{equation}
\psi_{\pm}(z,\cdot,x_0)\in L^2([x_0,\pm\infty);dx), \;
z\in\bbC\backslash\sigma(H), \quad \psi_{\pm}(z,x_0,x_0)=1. \lb{2.30}
\end{equation}
The normalization in \eqref{2.30} shows that
$\psi_{\pm}(z,\cdot,x_0)$ is of the type
\begin{equation}
\psi_{\pm}(z,x,x_0)=\theta (z,x,x_0)
+m_{\pm}(z,x_0)\phi(z,x,x_0),
\quad  z\in\bbC\backslash\sigma(H), \; x\in\bbR, \lb{2.31}
\end{equation}
with $m_{\pm}(z,x_0)$, the Weyl--Titchmarsh $m$-functions associated with the
Dirichlet half-line Schr\"odinger operators $H_{\pm,x_0}$ in \eqref{2.5}.
We recall that $\pm m_{\pm,x_0}$ are Herglotz functions and that
\begin{equation}
\ol{m_{\pm}(z,x_0)}=m_{\pm}(\ol z,x_0), \quad
z\in\bbC\backslash\bbR.  \lb{2.32}
\end{equation}
We also note that
\begin{align}
g(z,x)&=\f{\psi_{+}(z,x,x_0)\psi_{-}(z,x,x_0)}
{W(\psi_{+}(z,x,x_0),\psi_{-}(z,x,x_0))}
=\f{\psi_{+}(z,x,x_0)\psi_{-}(z,x,x_0)}{m_{-}(z,x_0)
-m_{+}(z,x_0)},  \lb{2.33} \\
m_{\pm}(z,x)&=\f{\psi_{\pm}'(z,x,x_0)}{\psi_{\pm}(z,x,x_0)},
\quad z\in\bbC\backslash\cE, \; x\in\bbR.
\lb{2.33a}
\end{align}
By \eqref{2.23} and \eqref{2.32} one concludes that
\begin{equation}
m_+(\lambda+i0,x_0)=m_-(\lambda-i0,x_0)=\ol{m_-(\lambda+i0,x_0)} \,
\text{ for a.e.\ $\lambda\in\cE$}  \lb{2.34}
\end{equation}
and together with \eqref{2.29} (implying that $\phi(\lambda,x,x_0)$ and
$\theta(\lambda,x,x_0)$ are real-valued for all
$(\lambda,x,x_0)\in\bbR^3$) this proves that
\begin{equation}
\text{ for each $x\in\bbR$, \, $g(\lambda+i0,x)$ is purely imaginary for
a.e.\ $\lambda\in\cE$.}  \lb{2.35}
\end{equation}
By \eqref{2.3} this proves \eqref{2.24}.

Since for each $x\in\bbR$, $g(z,x)$ is necessarily real-valued in spectral
gaps $[E_0,\infty)\backslash\cE$, $\xi(\lambda,x)\in\{0,1\}$ for
$\lambda\in (a_j,b_j)$, $j\in J$, whenever $g(\lambda,x)\neq 0$.
Moreover, since
\begin{equation}
\f{d}{dz}(H-z)^{-1}=(H-z)^{-2}, \quad z\in\bbC\backslash\sigma(H),
\lb{2.36}
\end{equation}
one has for fixed $x\in\bbR$,
\begin{equation}
\f{d}{d\lambda}g(\lambda,x)=\int_{\bbR} dx'\,
G(\lambda,x,x')^2 > 0, \quad \lambda\in (-\infty,E_0)\cup\bigcup_{j\in J}
(a_j,b_j).   \lb{2.37}
\end{equation}
Thus, for fixed $x\in\bbR$, $g(\cdot,x)$ is strictly monotonically
increasing on each interval $(a_j,b_j)$, $j\in J$ and hence one obtains
the following behavior of $\xi(\cdot,x)$ on $[E_0,\infty)\backslash\cE$:
\begin{align}
\xi(\lambda,x)&=\begin{cases} 1, & \lambda\in(a_j,\mu_j(x)) \\
0, & \lambda\in(\mu_j(x),b_j) \end{cases} \; \text{ if $\mu_j(x)\in
(a_j,b_j)$,  $j\in J$,}  \lb{2.38} \\
\xi(\lambda,x)&=\begin{cases} 1, & \text{$\lambda\in (a_j,b_j)$ \, if
$g(\lambda,x)<0$, $\lambda\in (a_j,b_j)$, $j\in J$,}\\
0, & \text{$\lambda\in(a_j,b_j)$ \, if
$g(\lambda,x)>0$, $\lambda\in (a_j,b_j)$,  $j\in J$,}  \end{cases}
\lb{2.39}
\end{align}
where
\begin{equation}
g(\mu_j(x),x)=0, \, \text{ whenever } \, \mu_j(x)\in (a_j,b_j), \; j\in J.
\lb{2.40}
\end{equation}
For convenience we will also define
\begin{equation}
\mu_j(x)=\begin{cases} b_j & \text{if $g(\lambda,x)<0$,
$\lambda\in(a_j,b_j)$, $j\in J$,} \\
a_j & \text{if $g(\lambda,x)>0$, $\lambda\in(a_j,b_j)$, $j\in J$.}
\end{cases}  \lb{2.42}
\end{equation}
Insertion of \eqref{2.38} and \eqref{2.39} into the trace formula
\eqref{2.6}, and using the convention \eqref{2.42} then yields
\begin{equation}
V(x)=E_0+\sum_{j\in J} [a_j+b_j-2\mu_j(x)] \, \text{ for a.e.\
$x\in\bbR$.}  \lb{2.43}
\end{equation}
Absolute convergence of the trace formula in \eqref{2.43} is assured by
the finite gap length condition \eqref{2.17} of $\cE$. Since by
definition,
\begin{equation}
\mu_j(x)\in [a_j,b_j], \quad j\in J, \lb{2.44}
\end{equation}
condition \eqref{2.17} proves
\begin{equation}
|V(x)|\leq |E_0|+\sum_{j\in J} |b_j-a_j| < \infty \, \text{ for a.e.\
$x\in\bbR$.}
\end{equation}
In particular, $V\in L^\infty(\bbR;dx)$, proving \eqref{2.25}.

Since by \eqref{2.33a},
\begin{equation}
m_\pm(z,x)=\f{\psi_\pm^{\prime} (z,x,x_0)}{\psi_\pm(z,x,x_0)}=
\f{\theta^{\prime}(z,x,x_0)+m_\pm(z,x_0)\phi^{\prime}(z,x,x_0)}{\theta(z,x,x_0)
+m_\pm(z,x_0)\phi(z,x,x_0)},   \lb{2.45}
\end{equation}
one infers that for each $x\in\bbR$,
\begin{equation}
\lim_{\varepsilon\downarrow 0}
m_+(\lambda+i\varepsilon,x)=\ol{\lim_{\varepsilon\downarrow 0}
m_-(\lambda+i\varepsilon,x)}= \lim_{\varepsilon\downarrow 0}
m_-(\lambda-i\varepsilon,x) \, \text{ for a.e.\ $\lambda\in\cE$,}  \lb{2.46}
\end{equation}
since for all $(x,x_0)\in\bbR^2$, $\theta(z,x,x_0)$ and $\phi(z,x,x_0)$
are entire with respect to $z$ and real-valued for all $z\in\bbR$ (cf.\
\eqref{2.29}). Thus, \eqref{2.23} holds for all $x_0\in\bbR$ and Definition
\ref{d2.3}\,$(iii)$ is independent of the choice of $x_0\in\bbR$.
\end{proof}

We will also call a Schr\"odinger operator $H$ reflectionless if its
potential coefficient $V$ is reflectionless since hardly any confusion can
arise in this manner.

Incidentally, inserting \eqref{2.38} and \eqref{2.39} into \eqref{2.8}
yields the absolutely convergent product expansion for the diagonal
Green's function $g$,
\begin{equation}
g(z,x)=\f{i}{2(z-E_0)^{1/2}} \prod_{j\in J}
\f{[z-\mu_j(x)]}{[(z-a_j)(z-b_j)]^{1/2}}, \quad z\in\bbC\backslash\cE, \;
x\in\bbR.  \lb{2.50}
\end{equation}

At first sight, our class $\SY$ of potentials $V$ appears to be larger than the
class studied by Sodin and Yuditskii in \cite{SY95} (see also \cite{SY96}), since we did not assume continuity of $V$ from the outset. A careful examination of the proofs in 
Sodin and Yuditskii in \cite{SY95} and \cite{SY97}, however, shows 
that the assumption of continuity of $V$ is unnecessary and automatically implied by conditions $(i)$--$(iii)$ in Definition \ref{d2.3}. Moreover, all elements $V\in\SY$ are uniformly almost periodic. 

In the following we denote by  $\AP$ the set of uniformly (i.e., Bohr) almost periodic (and hence continuous) functions on $\bbR$. 

\begin{theorem} \lb{t2.5}
In addition to the boundedness and reflectionless properties of the class $\SY$ described in Theorem \ref{t2.4} one has that
\begin{equation}
\SY \subset \AP.   \lb{2.51}
\end{equation}
In particular,
\begin{equation}
\SY \subset C(\bbR)   \lb{2.51a}
\end{equation}
and hence the class $\SY$ introduced in Definition \ref{d2.3} coincides with the one studied in \cite{SY95}.
\end{theorem}

While a proof of \eqref{2.51}, combining a variety of results of \cite{SY95} and \cite{SY97}, is beyond the scope of this paper, we nevertheless briefly sketch some of the steps involved, simultaneously filling a gap in the proof of the approximation theorem in Section 5 of \cite{SY95}: \\
$(i)$ We recall the notion of Dirichlet data $\cD_j(x)$ and Dirichlet divisors $\cD(x)$
\begin{align}
\begin{split}
&\cD_j(x)=(\mu_j(x),\sigma_j(x)), \quad \mu_j(x)\in [a_j,b_j], \quad \sigma_j(x)\in\{1,-1\}, 
\; j\in J, \lb{2.51b}  \\
&\cD(x)=\{\cD_j(x)\}_{j\in J}, \quad x\in\bbR,
\end{split}
\end{align}
where $\sigma_j(x)=1$ (resp., $\sigma_j(x)=-1$) if $\mu_j(x)\in(a_j,b_j)$ is a right (resp., left) Dirichlet eigenvalue associated with the half-line $[x,\infty)$ (resp., $(-\infty,x]$), that is, $\mu_j(x)\in(a_j,b_j)$ is a pole of $m_+(z,x)$ (resp., $m_-(z,x)$). If $\mu_j(x)$ coincides with one of the endpoints $a_j$ or $b_j$ one identifies $(\mu_j(x),1)$ and 
$(\mu_j(x), -1)$ and simply writes $(\mu_j(x))$. To avoid that $\sigma_j(x)$ formally becomes undefined whenever 
$\mu_j(x)$ hits an endpoint $a_j$ or $b_j$, one can follow \cite{Cr89} and change variables from $(\mu_j(x),\sigma_j(x))$ to $\varphi(x)$ according to
\begin{align}
&\mu_j(x)=[(a_j+b_j)+(b_j-a_j)\cos(\varphi(x))]/2,   \lb{2.51c} \\
&\mu_j(x)\in[a_j,b_j], \quad 0<\varphi_j(x)<\pi \, \text{ for } \, \sigma_j(x)=1, \; 
\pi<\varphi_j(x)<2\pi \, \text{ for } \, \sigma_j(x)=-1,  \no 
\end{align}
without ambiguity as $\mu_j(x)$ equals $a_j$ or $b_j$. As $x$ varies in $\bbR$, 
$\varphi_j(x)$ traces out $[0,2\pi]$ and hence we may identify the corresponding motion of $\cD_j(x)$ with a circle $\bbT_j$ and that of $\cD(x)$ with a torus $\cD(\cE)$ of dimension equal to the cardinality $|J|$ of the index set $J$. The torus $\cD(\cE)$ will be equipped with the product topology which in turn may be generated by the norm
\begin{align}
\|\varphi\|=\sum_{j\in J} 2^{-j} |\varphi_j|_{{\rm mod} (2\pi)}, \quad 
\varphi=\{\varphi_j\}_{j\in J}\in \cD(\cE), \quad \varphi_j\in\bbT_j, \; j\in J.  \lb{2.51d}
\end{align}
$(ii)$ Fix $x_0\in\bbR$ and a divisor $\cD(x_0)\in\cD(\cE)$. Then $\cD(x_0)$ uniquely determines a potential $V\in\SY$ a.e.\ on $\bbR$. The proof of this fact utilizes the Borg--Marchenko uniqueness result for Schr\"odinger operators on the half-lines $(-\infty,x_0]$ and $[x_0,\infty)$ (to the effect that $m_{\pm}(z,x_0))$ uniquely determine the potential  $V$ a.e.\ on $[x_0,\pm\infty)$), the Abel map defined in terms of the harmonic measure of subsets of $\cE$ at $z\in\bbC\backslash\cE$ with respect to $\bbC\backslash\cE$, and the absence of singular inner factors of $\wti m_{\pm}(z,x_0)$ and $\wti m_+(z,x_0)
\wti m_-(z,x_0)$, where $\wti m_{\pm}(z,x_0)=m_{\pm}(z,x_0)-m_{\pm}(-1,x_0)$. 

\smallskip
\noindent
$(iii)$ As a final step one proves the  finite-band approximation theorem. In this context we note, that the proof of the corresponding claim in \cite[Sect.\ 5]{SY95} is valid
only in the sense of convergence on all compact subsets of the real axis, even though
\cite{SY95} contains the erroneous statement that the proof presented yields uniform convergence on $\bbR$. Indeed (following the notation in \cite{SY95} in this paragraph with the exception that we denote $E$ in \cite{SY95} consistently by $\cE$), 
the last displayed formula in \cite[p.\ 652]{SY95} only claims pointwise convergence for all $t\in\bbR$ whereas the first formula in \cite[p.\ 653]{SY95} already claims uniform convergence. (E.g., the sequence $\{e^{i (t/n)}\}_{n\in\bbN}$ pointwise approximates the constant function $f(t)=1$, $t\in\bbR$, as $n\uparrow\infty$, but not uniformly on the whole real axis). However, this error is readily fixed: One only has to change the choice of  approximation of the spectrum $\cE$ by finite band sets $\{\cE^{(N)}\}_{N\in\bbN}$. To have uniform convergence as $N\uparrow\infty$ with respect to $t\in\bbR$ in the last displayed formula on p. 652 in \cite{SY95},  one has to keep fixed the initial frequencies of the approximating sets $\delta(\cE^{(N)})$: That is, the initial coordinates of the vector
$\delta(\cE^{(N)})$  should {\it coincide} with the corresponding coordinates
of $\delta(\cE)$.  Then $A^{(N)}(D^{(N)})+\delta(\cE^{N})t$ indeed converges to $A(D)+\delta(\cE)t$ as $N\uparrow\infty$, uniformly in $t\in\bbR$, with respect to the chosen topology in $\pi^*(\cE)$.

This choice of approximation is well-known in spectral theory. The spectral set $\cE$ is usually described in terms of a conformal map on a comb-like domain. One defines
\begin{equation}
\Pi=\{w=u+iv \in\bbC \,|\, u>0, \, v>0 \}\big\backslash\bigcup_{j\in J}
\{w=u_j+iv\in\bbC \,|\, 0<v\le h_j \}  \lb{5.15}
\end{equation}
with $\sum_{j\in J} h_j<\infty$.
Let $w=\Theta(\cdot)$ be a conformal map of the upper half--plane $\bbC_+$
onto $\Pi$ with the normalizations $E_0\mapsto 0$, $-\infty\mapsto i\infty$,
and in addition,
\begin{equation}
\Theta(z)\underset{z\to-\infty}{=} z^{1/2}[1+\oh(1)].  \lb{5.16}
\end{equation}
Then $\cE=\Theta^{-1}(\bbR_+)$, pre-images of slits form the system of
intervals $(a_j,b_j)$, which represent the spectral gaps, 
and  the set $\{u_j\}_{j\in J}$ forms the collection of frequencies of the corresponding almost periodic potential.

Let $\{J_N\}_{N\in\bbN}$ be an exhaustion of $J$ by finite subsets:
\begin{equation}
\cdots J_N\subset J_{N+1} \subset \dots \subset J.
\end{equation}
Next, define
\begin{equation}
\Pi_N=\{w=u+iv\in\bbC \,|\, u>0, \, v>0 \}\big\backslash\bigcup_{j\in J_N}
\{w=u_j+iv\in\bbC \,|\, 0<v\le h_j \},  \lb{5.15f}
\end{equation}
introduce $\Theta_N$ as the conformal map onto $\Pi_N$ under the same normalizations, and let $\cE^{(N)}=\Theta_N^{-1}(\bbR_+)$, $N\in\bbN$.

Finally we define the  map $\phi_N$ (a conformal map onto its image) by the diagram 
\begin{equation}\label{18j31}
\begin{array}{lll}
\bbC_+
&\stackrel{\Theta}{\longrightarrow}& \Pi \\
\Big{\downarrow}\phi_N & &
\Big\downarrow{\rm id}\\
\bbC_+ &\stackrel{\Theta_N}{\longrightarrow}& \Pi_N
\end{array}
\end{equation}
where ${\rm id}$ denotes the identity map. 

Let $V\in\SY$ with uniquely associated divisor $\cD(x_0)$ for some fixed $x_0\in\bbR$. Associated with $\cD(x_0)=\{\cD_j(x_0)\}_{j\in J}$ are the ``projections'' 
\begin{equation}
\cD^{(N)}(x_0)=\{\phi_{N}(\cD_j(x_0))\}_{j=1}^N, \quad 
\phi_{N}((\mu_j,\sigma_j))=((\phi_{N}(\mu_j),\sigma_j)), \; 1\leq j\leq N, 
\end{equation} 
associated with $\cE^{(N)}$. 
Applying the classical Jacobi inversion formula to the compact hyperelliptic Riemann surface associated with $\cE^{(N)}$, one uniquely determines $\cD^{(N)}(x)$ for all 
$x\in\bbR$ from the initial data $\cD^{(N)}(x_0)$ and then constructs an algebro-geometric finite-band potential $V^{(N)}$ by using the trace formula 
\begin{equation}
V^{(N)}(x)=E_0+\sum_{j=1}^N [\phi_N(a_j)+\phi_N(b_j)-2\phi_N(\mu_j(x))], \quad x\in\bbR.  \lb{2.51f}
\end{equation}
In addition, one proves that the set of algebro-geometric finite-band potentials $V^{(N)}$ is precompact in the topology of uniform convergence on $\bbR$, and that each limit potential belongs to the class $\SY$. In particular, upon embedding $\cD^{(N)}(x_0)$ into $\cD(\cE)$ by introducing 
\begin{equation}
{\wti \cD}^{(N)}(x_0)=\{(\mu_j(x_0),\sigma_j(x_0))\}_{j=1}^N \cup 
\{(b_k)\}_{k=N+1}^{\infty} \in \cD(\cE),      \lb{2.51g}
\end{equation}
there exists a subsequence ${\wti \cD}^{(N_{\ell})}(x_0)$, $\ell\in\bbN$, which converges to $\cD(x_0)$ as $\ell\uparrow\infty$. Using the linearization property of the Abel map, one can then prove that actually ${\wti \cD}^{(N_{\ell})}(x)$ converges to $\cD(x)$ as 
$\ell\uparrow\infty$ uniformly with respect to $x\in\bbR$. By the trace formula 
\eqref{2.51f}, the corresponding potentials $V^{(N_{\ell})}$ converges to $V$ as 
$\ell\uparrow\infty$ uniformly on $\bbR$. Since $V^{(N_{\ell})}$ are continuous and quasi-periodic, one thus concludes that $V$ is continuous and almost periodic.

The actual proofs of the statements in items $(i)$--$(iii)$ in \cite{SY95} and \cite{SY97} rely in part on potential theoretic techniques, most notably, the theory of character-automorphic Hardy spaces. However, due to the above correction in choosing the proper approximating sequence, one additional result is needed to ensure convergence
of  $V^{(N)}(x_0)$ (cf.\ \eqref{2.51f}) to $V(x_0)$ for a fixed divisor $\cD(x_0)$. 
(In the original paper \cite{SY95}, $\phi_N$ played the role of the identity map
at this place and hence the corresponding claim was obvious.)

\begin{lemma}Let  
$\phi_N$ be defined by \eqref{18j31}.
Then
\begin{equation}\label{9feb4}
\phi_N (b_j)-\phi_N (a_j)\ge b_j-a_j \, \text{ for all $j\in J_N$.}
\end{equation}
On the other hand,
\begin{equation}\label{9feb5}
\sum_{j\in J_N}[\phi_N (b_j)-\phi_N (a_j)]\le \sum_{j\in J}(b_j-a_j).
\end{equation}
In particular, this implies the uniform estimate
\begin{equation}
\sum_{\{j\in J_N \,|\, j\ge N\}}[\phi_N(b_j)-\phi_N(a_j)] \le  
\sum_{\{j\in J \,|\, j\ge N\}}(b_j-a_j)
\end{equation}
with respect to $N\in\bbN$. 
\end{lemma}
\begin{proof}
We introduce the intermediate domain 
\begin{align}
 \Pi_N^{\tau} &=\{w=u+iv\in\bbC\,|\, u>0, \, v>0 \} \no \\ 
& \quad \;\; \big\backslash
\bigcup_{j\in J_N\cup\{k\in J \,|\, u_k \ge\tau\}}
\{w=u_j+iv\in\bbC \,|\, 0<v\le h_j \},   
\end{align}
and assume that $\tau>0$ is a suficiently large parameter, $\tau\not=u_j$, $j\in J$, that later on will be sent to $+\infty$, . The corresponding
function 
$\phi_N^{\tau}(z)=(\Theta_N^{\tau})^{-1}(\Theta(z))$ maps the upper half plane in the upper half plane and therefore
possesses the representation
\begin{equation}\label{9feb7}
\phi_N^{\tau}(z)=z+C^\tau_N +\int_0^{x(\tau)}\frac{d\rho^\tau_N (x)}{x-z},
\end{equation}
where $d \rho_N^\tau$ is a nonnegative measure, $C_N^\tau\in\bbR$, and 
$x(\tau)=\Theta^{-1}(\tau)$. We note that for each interval $(a,b)$, which is free of charge of the measure $d\rho$, we have
\begin{equation}
\phi_N^{\tau}(b)-\phi_N^{\tau}(a)=b-a+\int_0^{x(\tau)} 
\frac{(b-a) d\rho^\tau_N(x)}{(x-b)(x-a)}\ge b-a.
\end{equation}
That is, the map is expanding on $(a,b)$. In particular,
\begin{equation}\label{9feb9}
\phi_N^{\tau}(b_j)-\phi_N^{\tau}(a_j)\ge b_j-a_j, \quad j\in J_N.
\end{equation}
On the other hand, the set
\begin{equation}
\cE^\tau=[0,x(\tau)]\setminus \bigcup_{\{j\in J \,|\, u_j\le \tau\}}(a_j,b_j)
\end{equation}
is also free of charge of the measure $d\rho$, that is, the measure of the image of this set is also larger than the measure of the set $\cE^\tau$. Thus (taking into account that 
$\phi_N^{\tau}(0)=0$),
\begin{equation}\label{9feb11}
\phi_N^{\tau}(x(\tau))-
\sum_{j\in J_N}(\phi_N^{\tau}(b_j)-\phi_N^{\tau}(a_j))\ge x(\tau)-
 \sum_{\{j\in J \,|\, u_j\le \tau\}} (b_j-a_j),
\end{equation}
 and hence
 \begin{equation}
\sum_{j\in J_N} [\phi_N^{\tau}(b_j)-\phi_N^{\tau}(a_j)]+ x(\tau)-\phi_N^{\tau}(x(\tau))\le
 \sum_{\{j\in J \,|\, u_j\le \tau\}}(b_j-a_j)\le\sum_{j\in J}(b_j-a_j).
\end{equation}
Using again the integral representation \eqref{9feb7} one gets $x(\tau)\ge\phi_N^{\tau}(x(\tau))$ and therefore,
\begin{equation}\label{9feb13}
\sum_{j\in J_N} [\phi_N^{\tau}(b_j)-\phi_N^{\tau}(a_j)]
 \le\sum_{j\in J}(b_j-a_j).
\end{equation}

Keeping $N$ fixed and sending $\tau$ to infinity we get \eqref{9feb4} from
\eqref{9feb9}, and \eqref{9feb5} from \eqref{9feb13} by Caratheodory's theorem 
\cite{Ca12}, \cite[Sect.\ II.5]{Go69}. 
\end{proof}

Finally, we conclude this section with a converse to Theorem \ref{t2.4} due to
Sodin and Yuditskii \cite{SY95}.

\begin{theorem} [\cite{SY95}] \lb{t2.6}
Let $V\in L^1_{\loc}(\bbR;dx)$ be real-valued and assume that
$H_{\min}$ is bounded from below and its unique self-adjoint operator $H$
has spectrum $\cE$, where $\cE$ is a homogeneous set of finite gap length of
the form \eqref{2.15}. In addition, suppose that $V$ is reflectionless in the
sense that
\begin{equation}
\text{for each $x\in\bbR$, } \, \xi(\lambda,x)=1/2 \, \text{ for
a.e.\ $\lambda\in\cE$.}    \lb{2.52}
\end{equation}
Then,
\begin{equation}
V\in \SY   \lb{2.53}
\end{equation}
and
\begin{equation}
\text{for all $x\in\bbR$ and for a.e.\ $\lambda\in\cE$, } \,
\lim_{\varepsilon\downarrow 0} m_+(\lambda+i\varepsilon,x)=
\lim_{\varepsilon\downarrow 0} m_-(\lambda-i\varepsilon,x).  \lb{2.54}
\end{equation}
\end{theorem}
\begin{proof}
In the Appendix of \cite{SY95} it is proved that \eqref{2.52} implies that
\begin{equation}
\text{for a.e.\ $\lambda\in\cE$, } \,
\lim_{\varepsilon\downarrow 0} m_+(\lambda+i\varepsilon,x_0)=
\lim_{\varepsilon\downarrow 0} m_-(\lambda-i\varepsilon,x_0),  \lb{2.55}
\end{equation}
implying \eqref{2.53}. Using \eqref{2.45}, and noticing that
$\theta(\lambda,x,x_0)$, $\theta'(\lambda,x,x_0)$, $\phi(\lambda,x,x_0)$,
and $\phi'(\lambda,x,x_0)$ are all real-valued for $\lambda\in\cE$, one
concludes \eqref{2.54} (this argument has also been used to arrive at
\eqref{2.46}).
\end{proof}


In this section we shall prove that the set $\cE$, the spectrum of a
reflectionless Schr\"odinger operator $H$ with potential $V$ in the
Sodin--Yuditskii class $\SY$, coincides with the absolutely continuous
spectrum, $\sigma_{\rm ac}(H)$, of $H$. In fact, we will provide two
elementary proofs of this fact.

We start by recalling the following result on essential supports of
the absolutely continuous spectrum of Schr\"odinger operators on the
real line proven in \cite{GS96} (see also  \cite[p.\ 383]{AD56}). For
completeness we will provide its short proof.

Let $V\in L^1_{\loc}(\bbR;dx)$ be real-valued, assume that the
differential expression $L=-d^2/dx^2+V$ is in the limit point case at
$+\infty$ and $-\infty$, and denote by $H$ the corresponding self-adjoint
realization of $L$ in $L^2(\bbR;dx)$ introduced in \eqref{2.1}. In
addition, let $g(z,x)$, $z\in\bbC\backslash\sigma(H)$, $x\in\bbR$, and
$\xi(\lambda,x)$ for a.e.\ $\lambda\in\bbR$, $x\in\bbR$, be defined as in
\eqref{2.2} and \eqref{2.3}. By $\ol{A}^e$ we denote the essential closure
of the Lebesgue measurable set $A\subset\bbR$ (cf.\ Appendix \ref{A} for
details).

\begin{theorem} [\cite{GS96}] \lb{t3.1}
For each $x\in\bbR$, the set
\begin{equation}
\{\lambda\in\bbR\,|\, 0<\xi(\lambda,x)<1\}  \lb{3.1}
\end{equation}
is an essential support of the absolutely continuous spectrum,
$\sigma_{\rm ac}(H)$, of $H$. In particular, for each $x\in\bbR$,
\begin{equation}
\sigma_{\rm ac}(H)=\ol{\{\lambda\in\bbR\,|\, 0<\xi(\lambda,x)<1\}}^e.
\lb{3.2}
\end{equation}
\end{theorem}
\begin{proof}
Pick $x\in\bbR$ and denote by $m_\pm(\cdot,x)$ the two half-line
Weyl--Titchmarsh
functions associated with the $H_{\pm,x}^D$ as defined in \eqref{2.5}. Using
\eqref{2.3} and \eqref{2.7}, the following three sets coincide up to sets of
Lebesgue measure zero:
\begin{align}
& (i) \hspace*{3mm} \text{$\{\lambda\in\bbR\,|\,
0<\xi(\lambda,x)<1\}$,}  \lb{3.3} \\
& (ii) \hspace*{2mm}
\text{$\{\lambda\in\bbR\,|\, 0< \Im[g(\lambda+i0,x)]$  exists
finitely$\}$,}  \lb{3.4} \\
& \text{$(iii)$ $S_{+,x}\cup S_{-,x}$}, \quad
\text{$S_{\pm,x}=\{\lambda\in\bbR\,|\, 0< \Im[m_\pm(\lambda+i0,x)]$ exists
finitely$\}$.}  \lb{3.5}
\end{align}
By \eqref{B.11}, $S_{\pm,x}$ are essential supports of $\sigma_{\rm
ac}(H_{\pm,x})$.
But since $H$ and $H_{-,x}\oplus H_{+,x}$ have the same absolutely continuous
spectrum (their resolvents only differ by a rank-one perturbation), \eqref{3.3}
represents an essential support of $\sigma_{\rm ac}(H)$ proving the first claim
concerning \eqref{3.1}. Equation \eqref{3.2} is then clear from \eqref{A.14}.
\end{proof}

Next, let $\SY$ be the Sodin--Yuditskii class introduced in Definition
\ref{d2.3}. Given $V\in\SY$, we consider the Schr\"odinger operator $H$
in $L^2(\bbR;dx)$ by
\begin{equation}
Hf=-f''+Vf, \quad f\in \dom(H)=H^{2,2}(\bbR)   \lb{3.6}
\end{equation}
with spectrum, $\sigma(H)$, of $H$ given by
\begin{equation}
\sigma(H)=\cE, \quad \cE= [E_0,\infty)\big\backslash\bigcup_{j\in J}
(a_j,b_j), \quad J\subseteq \bbN,  \lb{3.7}
\end{equation}
according to \eqref{2.15} and \eqref{2.22}.

The principal result of this section then reads as follows.

\begin{theorem}  \lb{t3.2}
Let $V\in\SY$ and denote by $H$ the associated Schr\"odinger operator
defined in \eqref{3.6}. Then, the absolutely continuous spectrum of $H$
coincides with $\cE$,
\begin{equation}
\sigma_{\rm ac}(H) = \cE.   \lb{3.8}
\end{equation}
Moreover, $\sigma_{\rm ac}(H)$ has uniform multiplicity equal to two.
\end{theorem}
\begin{proof}
Since $V$ is reflectionless by Theorem \ref{t2.4}, one has
\begin{equation}
\text{for each $x\in\bbR$, } \;
\xi(\lambda,x)=(1/\pi)\Im[\ln(g(\lambda+i0,x)]= 1/2 \,
\text{ for a.e.\ $\lambda\in\cE$.}  \lb{3.9}
\end{equation}
By \eqref{3.2}, this implies
\begin{equation}
\sigma_{\rm ac}(H)=\ol{\{\lambda\in\bbR \,|\, 0<\xi(\lambda,x_0)<1\}}^e
\supseteq \ol{\cE}^e  \lb{3.10}
\end{equation}
for some $x_0\in\bbR$. By the definition of homogeneity of $\cE$ in
\eqref{2.16}, one infers,
\begin{equation}
\text{for all $\lambda\in\cE$ and all $\varepsilon>0$, } \;
|(\lambda-\varepsilon,\lambda+\varepsilon)\cap\cE|>0  \lb{3.11}
\end{equation}
and hence,
\begin{equation}
\cE\subseteq \ol{\cE}^e.  \lb{3.12}
\end{equation}
Thus one obtains by \eqref{3.10} that
\begin{equation}
\sigma_{\rm ac}(H)\supseteq \ol{\cE}^e \supseteq \cE.  \lb{3.13}
\end{equation}
Since $\sigma(H)=\cE$ by \eqref{3.7}, this proves \eqref{3.8}.

Equations \eqref{2.7} and \eqref{2.23} imply
\begin{equation}
-1/g(\lambda+i0,x_0)=\pm 2i \, \Im[m_\pm(\lambda+i0,x_0)] \,
\text{ for a.e.\ $\lambda\in\cE$.}    \lb{3.14}
\end{equation}
Finally, combining \eqref{2.23}, \eqref{3.14}, and \eqref{B.38} then yields
that $\sigma_{\rm ac}(H)=\cE$ has uniform spectral multiplicity two since
\begin{equation}
\text{for a.e.\ $\lambda\in\cE$, } \;
0<\pm\Im[m_\pm(\lambda+i0,x_0)]<\infty.   \lb{3.15}
\end{equation}
\end{proof}

That $\sigma(H)$ has uniform multiplicity equal to two is in accordance
with Theorem\ 9.1 in Deift and Simon \cite{DS83}.

Although it was not necessary to use the following information in the
proof of Theorem \ref{t3.2}, we note that by \eqref{A.10b}, one has
\begin{equation}
\ol{\cE}^e \subseteq \ol\cE = \cE.  \lb{3.16}
\end{equation}
Thus, combining \eqref{3.12} and \eqref{3.16} leads to
\begin{equation}
\cE=\ol{\cE}^e=\ol \cE.   \lb{3.17}
\end{equation}

One can also give an alternative proof of \eqref{3.8} based on the
reflectionless property of $V$ as follows (still under the assumptions
of Theorem \ref{t3.2}):

\begin{proof}[Alternative proof of \eqref{3.8}]
Fix $x_0\in\bbR$. Since
$V\in\SY$, \eqref{2.24} yields that
\begin{equation}
\xi(\lambda,x_0)=\f{1}{\pi}\Im[\ln(g(\lambda+i0,x_0))]=\f{1}{2} \,
\text{ for a.e.\ $\lambda\in\cE$}  \lb{3.18}
\end{equation}
and hence
\begin{equation}
\Re[g(\lambda+i0,x_0)]=0 \, \text{ for a.e.\ $\lambda\in\cE$.} \lb{3.19}
\end{equation}
If there exists a measurable set $A\subseteq \cE$ of positive Lebesgue
measure, $|A|>0$, on which $\Im(g)$ vanishes, that is,
\begin{equation}
\Im[g(\lambda+i0,x_0)]=0 \, \text{ for a.e.\ $\lambda\in A$,} \lb{3.20}
\end{equation}
then \eqref{3.19} and \eqref{3.20} yield the existence of a subset
$B\subseteq \cE$ of positive Lebesgue measure, $|B|>0$, such that
\begin{equation}
g(\lambda+i0,x_0)=0 \, \text{ for a.e.\ $\lambda\in B$.} \lb{3.21}
\end{equation}
Since $g(\cdot,x_0)$ is a Herglotz function, the uniqueness property of
Herglotz functions in Theorem \ref{tB.2}\,$(ii)$ yields the contradiction
$g\equiv 0$. Thus, no such set $B\subseteq\cE$ exists and one concludes
that
\begin{equation}
\Im[g(\lambda+i0,x_0)]>0 \, \text{ for a.e.\ $\lambda\in \cE$.} \lb{3.22}
\end{equation}
Moreover, since $g(\lambda+i0,x_0)$ exists finitely for a.e.\
$\lambda\in\bbR$ by Theorem \ref{tB.2}\,$(i)$, this yields
\begin{equation}
0<\Im[g(\lambda+i0,x_0)]<\infty \, \text{ for a.e.\ $\lambda\in \cE$.}
\lb{3.23}
\end{equation}
By the asymptotic behavior of $g$ in \eqref{2.14b}, this shows that
\begin{equation}
g(z,x_0)=[m_-(z,x_0)-m_+(z,x_0)]^{-1}  \lb{3.24}
\end{equation}
is a Herglotz function of the type (cf.\ \eqref{B.3}, \eqref{B.8a}, and
\eqref{B.22}--\eqref{B.24})
\begin{equation}
g(z,x_0)=\int_{\bbR} \f{d\Omega_{0,0}(\lambda,x_0)}{\lambda -z},
\quad z\in\bbC_+.  \lb{3.25}
\end{equation}
Combining \eqref{3.23} and \eqref{B.5}, the absolutely continuous part
$d\Omega_{0,0, \rm ac}(\cdot,x_0)$ of the measure
$d\Omega_{0,0}(\cdot,x_0)$ is supported on $\cE$,
\begin{equation}
\supp\,[d\Omega_{0,0, \rm ac}(\cdot,x_0)]=\cE.  \lb{3.26}
\end{equation}
(This also follows from combining, \eqref{3.22}, \eqref{A.14}, and
\eqref{B.11}.) Next, we replace $m_\pm$ in \eqref{B.26}, \eqref{B.27} by
$m_\pm(\cdot,x_0)$ (cf.\ the paragraph preceding Theorem \ref{tB.7})
and introduce the trace measure $d\Omega^{\rm tr}(\cdot,x_0)$
associated with the trace $M^{\rm tr}(\cdot,x_0)$ of the $2\times 2$
Weyl--Titchmarsh matrix $M(\cdot,x_0)$ of the Schr\"odinger operator $H$
\begin{align}
&M^{\rm tr}(z,x_0)=M_{0,0}(z,x_0)+M_{1,1}(z,x_0)
=\f{1+m_-(z,x_0)m_+(z,x_0)}{m_-(z,x_0)-m_+(z,x_0)}
\lb{3.27} \\
& \hspace*{1.6cm} =a(x_0) + \int_{\bbR} d\Omega^{\rm tr}(\lambda,x_0)
\bigg(\frac{1}{\lambda -z)}-\frac{\lambda} {1+\lambda^2}\bigg), \quad
z\in\bbC_+, \lb{3.28} \\
& \, a(x_0)=\Re\big[M^{\rm tr}(i,x_0)\big], \quad d\Omega^{\rm
tr}(\cdot,x_0) = d\Omega_{0,0}(\cdot,x_0) + d\Omega_{1,1}(\cdot,x_0). \no
\end{align}
By \eqref{3.26} one infers that
\begin{equation}
\supp\,\big[d\Omega^{\rm tr}_{\rm ac}(\cdot,x_0)\big]=\cE  \lb{3.29}
\end{equation}
and hence
\begin{equation}
\sigma_{\rm ac}(H)=\cE  \lb{3.30}
\end{equation}
holds as a consequence of \eqref{B.28}.  
\end{proof}

In addition to
\begin{equation}
M_{0,0}(z,x)=g(z,x)=[m_-(z,x)-m_+(z,x)]^{-1}, \quad z\in\bbC_+, \;
x\in\bbR  \lb{3.31}
\end{equation}
we will also analyze
\begin{equation}
M_{1,1}(z,x)=h(z,x)=\f{m_-(z,x)m_+(z,x)}{m_-(z,x)-m_+(z,x)}, \quad
z\in\bbC_+, \;  x\in\bbR.  \lb{3.32}
\end{equation}
Recalling $m_{\pm}=\psi'_{\pm}/\psi_{\pm}$ (cf.\ \eqref{2.45}), we note
that $h$ is given by
\begin{equation}
h(z,x)=\f{\psi'_+(z,x,x_))\psi'_-(z,x,x_))}
{W(\psi_+(z,x,x_)),\psi_-(z,x,x_)))}, \quad z\in\bbC_+, \; x\in\bbR.
\lb{3.33}
\end{equation}
In analogy to the product expansion for $g$ in \eqref{2.50} one
then obtains
\begin{equation}
h(z,x)=\f{i}{2}\f{[z-\nu_0(x)]}{(z-E_0)^{1/2}} \prod_{j\in J}
\f{[z-\nu_j(x)]}{[(z-a_j)(z-b_j)]^{1/2}}, \quad z\in\bbC\backslash\cE, \;
x\in\bbR,   \lb{3.34}
\end{equation}
where (in analogy to $\mu_j(x)\in [a_j, b_j]$, $j\in J$, $x\in\bbR$, in
connection with $g(\cdot,x)$)
\begin{equation}
\nu_0(x)\leq E_0, \quad \nu_j(x)\in[a_j, b_j], \quad j\in J, \; x\in\bbR.
\lb{3.35}
\end{equation}
Moreover, \eqref{2.14a} and \eqref{3.32} imply the universal asymptotic
expansion, valid for all $x\in\bbR$,
\begin{equation}
h(z,x)\underset{\substack{\abs{z}\to\infty\\ z\in
C_\varepsilon}}{=} (i/2) z^{1/2}[1+o(1)]. \lb{3.36}
\end{equation}

\section{Absence of Singular Spectrum} \label{s4}

In this section we will prove the absence of the singular spectrum of
a reflectionless Schr\"odinger operator $H$ with potential $V$ in the
Sodin--Yuditskii class $\SY$, that is, we intend to prove that
\begin{equation}
\sigma(H)=\sigma_{\rm ac}(H)=\cE, \quad
\sigma_{\rm sc}(H)=\sigma_{\rm pp}(H)=\emptyset.  \lb{4.1}
\end{equation}
Unlike our two elementary proofs of $\sigma_{\rm ac}(H)=\cE$ in Section
\ref{s3} (cf.\ \eqref{3.8}), the proof of \eqref{4.1} relies on certain
techniques developed in harmonic analysis and potential theory associated
with domains $(\bbC\cup\{\infty\})\backslash\cE$.

We start with an elementary lemma which will permit us to reduce the
discussion of unbounded homogeneous sets $\cE$ (typical for Schr\"odinger
operators) to the case of a compact homogeneous sets $\wti\cE$ (as discussed
by Peherstorfer and Yuditskii \cite{PY03} and Sodin and Yuditskii \cite{SY97}
in connection with Jacobi operators).

\begin{lemma}  \lb{l4.1}
Let $m$ be a Herglotz function with representation
\begin{align}
\begin{split}
& m(z)=c+\int_{\bbR} d\omega(\lambda)
\bigg(\f{1}{\lambda-z}-\f{\lambda}{1+\lambda^2}\bigg), \quad z\in\bbC_+,
\\
& c=\Re[m(i)], \quad \int_{\bbR} \f{d\omega(\lambda)}
{1+\lambda^2} <\infty, \quad \bbR\backslash
\supp\,(d\omega)\neq\emptyset.   \lb{4.2}
\end{split}
\end{align}
Consider the change of variables
\begin{align}
\begin{split}
& z \mapsto \zeta=(\lambda_0-z)^{-1}, \;\;
z=\lambda_0-\zeta^{-1}, \quad z\in\bbC\cup\{\infty\}, \\
& \text{for some fixed $\lambda_0\in \bbR\backslash\supp \, (d\omega)$}.
\lb{4.3}
\end{split}
\end{align}
Then,
\begin{align}
& \wti m(\zeta)=m(z(\zeta))
=\wti c +\int_{\bbR} \f{d\wti\omega(\eta)}{\eta-\zeta}, \quad
\zeta\in\bbC_+, \lb{4.4} \\
& d\wti\omega(\eta)
=x^2d\omega\big(\lambda_0-\eta^{-1}\big)\big|_{\supp\,(d\wti\omega)},
\lb{4.5} \\
& \supp\, (d\wti\omega) \subseteq
\big[-\dist\,(\lambda_0,\supp\,(d\omega))^{-1},
\dist\,(\lambda_0,\supp\,(d\omega))^{-1}\big],  \lb{4.6} \\
& \wti c=c+\int_{\supp\,(d\wti\omega)} d\wti\omega(\eta)\,
\f{\lambda_0-\big(1+\lambda_0^2\big)\eta}
{1-2\lambda_0 \eta+\big(1+\lambda_0^2\big)\eta^2}. \lb{4.7}
\end{align}
In particular, $d\wti\omega$ is purely absolutely continuous if and only
if $d\omega$ is, and
\begin{equation}
\text{if
$d\omega(\lambda)=\omega'(\lambda)d\lambda\big|_{\supp \, (d\omega)}$,
then $d\wti\omega(\eta)
=\omega'\big(\lambda_0-\eta^{-1}\big)d\eta\big|_{\supp \,
(d\wti\omega)}$.}  \lb{4.8}
\end{equation}
\end{lemma}
\begin{proof}
This is a straightforward computation. We note that
\begin{equation}
\int_{\bbR} \f{d\omega(\lambda)}{1+\lambda^2} <\infty \, \text{ is
equivalent to } \int_{\bbR} d\wti\omega(\eta)<\infty.
\lb{4.9}
\end{equation}
\end{proof}

Next, we need the notion of a compact homogeneous set and
hence slightly modify Definition \ref{d2.2} as follows.

\begin{definition}  \lb{d4.2}
Let $\wti\cE\subset\bbR$ be a compact set which we may write as
\begin{equation}
\wti\cE= [E_0,E_1]\big\backslash\bigcup_{j\in \wti J}
\big(\wti a_j,\wti b_j\big),
\quad \wti J\subseteq\bbN, \lb{4.10}
\end{equation}
for some $E_0, E_1\in\bbR$, $E_0<E_1$, and $\wti a_j<\wti b_j$, where
$\big(\wti a_j,\wti b_j\big)\cap \big(\wti a_{j'},\wti
b_{j'}\big)=\emptyset$, $j,j'\in \wti J$, $j\neq j'$. Then $\wti \cE$ is
called {\it homogeneous} if
\begin{align}
\begin{split}
& \text{there exists an $\varepsilon>0$
such that for all $\eta\in\wti\cE$}  \lb{4.11} \\
& \text{and all $0<\delta<\diam \,\big(\wti\cE\,\big)$, \,
$\big|\wti\cE\cap(\eta-\delta,\eta+\delta)\big|\geq
\varepsilon\delta$.}
\end{split}
\end{align}
\end{definition}

The following theorem is a direct consequence of a result of
Zinsmeister \cite[Theorem\ 1]{Zi89} (cf.\ Theorem \ref{tB.4d}) and also a
special case of a result of Peherstorfer and Yuditskii \cite{PY03} 
proved by entirely different methods based on character automorphic Hardy
functions.

\begin{theorem} [\cite{Zi89}, Theorem\ 1, \cite{PY03}, Lemma\ 2.4] \lb{t4.4}
Let $\wti\cE\subset\bbR$ be a compact homogenous set and $r$ a Herglotz
function with representation
\begin{align}
\begin{split}
& r(\zeta)=a+\int_{\wti\cE} \f{d\mu(\eta)}{\eta-\zeta}, \quad \zeta
\in\bbC_+, \\
& a\in\bbR, \quad d\mu \, \text{ a finite measure, \,
$\supp \, (d\mu) \subseteq \wti\cE$.}  \lb{4.12}
\end{split}
\end{align}
Denote by $r(\eta+i0)=\lim_{\varepsilon\downarrow 0}r(\eta+i\varepsilon)$
the a.e.\ normal boundary values of $r$ and assume that
\begin{equation}
\Re[r(\cdot+i0)] \in L^1\big(\wti\cE;d\eta\big).  \lb{4.13}
\end{equation}
Then, $d\mu$ is purely absolutely continuous and hence
\begin{equation}
d\mu(\eta)=d\mu_{\rm ac}(\eta)=\f{1}{\pi} \Im[r(\eta+i0)]\,d\eta.
\lb{4.14}
\end{equation}
\end{theorem}
\begin{proof}
We note that $r$ is real-valued on $\bbR\backslash\wti\cE$ and
$\Im[r(\cdot +i0)]\in L^1_{\loc}(\bbR;d\eta)$ by general principles. Since
$\wti\cE$ is compact, one has $\Im[r(\cdot +i0)]\in L^1(\bbR;d\eta)$.
Thus, the integrability condition \eqref{4.13} is equivalent to
\begin{equation}
r(\cdot +i0) \in L^1\big(\wti\cE;d\eta\big).  \lb{4.15}
\end{equation}
Temporarily denoting by $d\eta$ the Lebesgue measure on $\bbR$, and
introducing the auxiliary Herglotz function
\begin{equation}
w(\zeta)=a+\bigg[\int_{\wti\cE} d\mu(\eta)\bigg/ \int_{\wti\cE} d\eta\bigg]
\int_{\wti\cE} \f{d\eta}{\eta-\zeta}, \quad \zeta\in\bbC_+,   \lb{4.16}
\end{equation}
we consider the function
\begin{equation}
F(\zeta)=r(\zeta)-w(\zeta)=\int_{\wti\cE} \f{d\nu(\eta)}{\eta-\zeta},
\quad \zeta\in\bbC_+  \lb{4.17}
\end{equation}
with the finite signed measure $d\nu$ given by
\begin{equation}
d\nu(\eta)=d\mu(\eta)-
\bigg[\int_{\wti\cE} d\mu(\eta)\bigg/ \int_{\wti\cE} d\eta\bigg]
\chi_{\wti\cE}(\eta) d\eta.
\lb{4.18}
\end{equation}
Then $F$ is analytic on $\bbC\backslash\wti\cE$ and
$\lim_{|z|\to\infty} zF(z)=0$. The symmetry property
$F(\ol z)=\ol{F(z)}$, the integrability condition \eqref{4.15}, and
the known existence of the nontangential limits
$r_{\pm}(\eta)=\lim_{\zeta\in\bbC_{\pm}, \, \zeta\to\eta} r(\zeta)$ for a.e.\
$\eta\in\wti\cE$, imply that
\begin{equation}
d\nu \in \wti{H^1}(\wti\cE)   \lb{4.19}
\end{equation}
(cf.\ the definition of $\wti{H^1}(\wti\cE)$ in \eqref{B.8K}).
Theorem \ref{tB.4d} then proves that $d\nu$, and hence $d\mu$, is purely
absolutely continuous.
\end{proof}

Combining Lemma \ref{l4.1} and Theorem \ref{t4.4}, we obtain the principal
result of this section. (Incidentally, this yields yet another proof of
\eqref{3.8}.)

\begin{theorem} \lb{t4.5}
Let $V\in\SY$ and denote by $H$ the associated Schr\"odinger operator
defined in \eqref{3.6}. Then, the spectrum of $H$ is purely absolutely
continuous and coincides with $\cE$,
\begin{equation}
\sigma(H)=\sigma_{\rm ac}(H) = \cE, \quad \sigma_{\rm sc}(H)=
\sigma_{\rm pp}(H)=\emptyset.   \lb{4.36}
\end{equation}
Moreover, $\sigma(H)$ has uniform multiplicity equal to two.
\end{theorem}
\begin{proof}
By Theorem \ref{t3.2} we need to prove the absence of the singular
spectrum of $H$,
\begin{equation}
\sigma_{\rm sc}(H)=\sigma_{\rm pp}(H)=\emptyset.   \lb{4.37}
\end{equation}
Recalling the absolutely convergent product expansion for the diagonal
Green's function $g$,
\begin{equation}
g(z,x)=\f{i}{2(z-E_0)^{1/2}} \prod_{j\in J}
\f{[z-\mu_j(x)]}{[(z-a_j)(z-b_j)]^{1/2}}, \quad z\in\bbC\backslash\cE, \;
x\in\bbR,  \lb{4.38}
\end{equation}
we pick $\lambda_0<\nu_0(x)$ (implying $\lambda_0 <E_0$) and introduce the
change of variables
\begin{equation}
z \mapsto \zeta=(\lambda_0-z)^{-1}, \;\;
z=\lambda_0-\zeta^{-1}, \quad z\in\bbC\cup\{\infty\}. \lb{4.39}
\end{equation}
This results in
\begin{equation}
\wti g(\zeta,x)=g(z(\zeta),x)=C(x)\f{\zeta^{1/2}}{(\zeta-\wti E_0)^{1/2}}
\prod_{j\in J}
\f{[\zeta-\wti\mu_j(x)]}{[(\zeta-\wti a_j)(\zeta-\wti b_j)]^{1/2}}, \quad
\zeta\in\bbC\backslash\wti\cE, \; x\in\bbR,  \lb{4.40}
\end{equation}
where
\begin{align}
& \wti E_0=(\lambda_0-E_0)^{-1}, \no \\
& \wti a_j=(\lambda_0-a_j)^{-1}, \quad \wti b_j=(\lambda_0-b_j)^{-1},
\quad  j\in J,  \lb{4.41} \\
& \wti \mu_j(x)=[\lambda_0-\mu_j(x)]^{-1}, \quad j\in J, \; x\in\bbR, \no
\end{align}
and
\begin{equation}
\wti\cE=\big[\wti E_0,0\big]\big\backslash\bigcup_{j\in J}
\big(\wti a_j,\wti b_j\big).  \lb{4.42}
\end{equation}
Since
\begin{equation}
g(z,x)\underset{z\downarrow -\infty}{=}\f{1}{2|z|^{1/2}}[1+\oh(1)], \quad
x\in\bbR,  \lb{4.43}
\end{equation}
one concludes
\begin{equation}
\wti g(\zeta,x)\underset{\zeta\downarrow
0}{=}\f{|\zeta|^{1/2}}{2}[1+\oh(1)], \quad x\in\bbR,  \lb{4.44}
\end{equation}
and hence
\begin{equation}
C(x)>0  \lb{4.45}
\end{equation}
in \eqref{4.40}. Next, we denote by $g(\lambda+i0,x)
=\lim_{\varepsilon\downarrow 0} g(\lambda+i\varepsilon,x)$ and $\wti
g(\eta+i0,x) =\lim_{\varepsilon\downarrow 0} \wti g(\eta+i\varepsilon,x)$
the a.e.\ normal boundary values of $g(\cdot,x)$ and $\wti g(\cdot,x)$,
$x\in\bbR$. By \eqref{2.35},
\begin{equation}
\text{ for each $x\in\bbR$, \, $g(\lambda+i0,x)$ is purely imaginary for
a.e.\ $\lambda\in\cE$}  \lb{4.46}
\end{equation}
and hence
\begin{equation}
\text{ for each $x\in\bbR$, \, $\wti g(\eta+i0,x)$ is purely imaginary
for a.e.\ $\eta\in\wti\cE$.}  \lb{4.47}
\end{equation}
In particular,
\begin{equation}
\text{ for each $x\in\bbR$, \, $\Re[\wti g(\eta+i0,x)]=0$ for a.e.\
$\eta\in\wti\cE$}  \lb{4.48}
\end{equation}
and thus the analog of \eqref{4.13} is satisfied in the context of
$\wti g(\cdot,x)$. Applying Theorem \ref{t4.4} thus proves
\begin{align}
\wti g(\zeta,x)=C(x)+\f{1}{\pi} \int_{\wti\cE} \f{[-i\wti
g(\eta+i0,x)]d\eta}{\eta-\zeta}, \quad \zeta\in\bbC\backslash\wti\cE, \;
x\in\bbR.  \lb{4.49}
\end{align}
Applying Lemma \ref{l4.1} to $g(\cdot,x)$ as represented in \eqref{3.25},
\begin{equation}
g(z,x)=\int_{\cE} \f{d\Omega_{0,0}(\lambda,x)}{\lambda -z}, \quad
z\in\bbC\backslash\cE, \; x\in\bbR,  \lb{4.50}
\end{equation}
then yields for the measure $d\Omega_{0,0}(\cdot,x)$,
\begin{align}
\begin{split}
& \text{for each $x\in\bbR$, \,
$d\Omega_{0,0}(\lambda,x) = d\Omega_{0,0, \rm ac}(\lambda,x)
=\f{1}{\pi}[-i g(\lambda+i0,x)]d\lambda$,} \\
& \text{$\supp\,(d\Omega_{0,0}(\cdot,x))=\cE$.}  \lb{4.51}
\end{split}
\end{align}
Next, we repeat this analysis for $h(\cdot,x)$. Introducing,
\begin{equation}
\wti\nu_0(x)=[\lambda_0-\nu_0(x)]^{-1}, \quad
\wti \nu_j(x)=[\lambda-\nu_j(x)]^{-1}, \quad j\in J, \; x\in\bbR,
\lb{4.52}
\end{equation}
$h(\cdot,x)$ transforms into
\begin{equation}
\begin{split}
& \wti h(\zeta,x)=h(z(\zeta),x)= -D(x)\f{[\zeta-\wti\nu_0(x)]}
{\zeta^{1/2}(\zeta-\wti E_0)^{1/2}}  \prod_{j\in J}
\f{[\zeta-\wti\nu_j(x)]}{[(\zeta-\wti a_j)(\zeta-\wti b_j)]^{1/2}}, \\
& \hspace*{8.25cm} \zeta\in\bbC\backslash\wti\cE, \; x\in\bbR.  \lb{4.53}
\end{split}
\end{equation}
The asymptotic behavior \eqref{3.36} then proves $D(x)>0$. At this point
one can again apply Lemma \ref{l4.1} and Theorem \ref{t4.4} to $\wti
h(\cdot,x)$ and obtains
\begin{align}
\begin{split}
& \text{for each $x\in\bbR$, \,
$d\Omega_{1,1}(\lambda,x) = d\Omega_{1,1, \rm ac}(\lambda,x)
=\f{1}{\pi}[-i h(\lambda+i0,x)]d\lambda$,} \\
& \text{$\supp\,(d\Omega_{1,1}(\cdot,x))=\cE$.}  \lb{4.54}
\end{split}
\end{align}
Combining \eqref{4.51} and \eqref{4.54} then yields
\begin{align}
\begin{split}
& d\Omega^{\rm tr}(\cdot,x_0) = d\Omega^{\rm tr}_{\rm ac}(\cdot,x_0)
=\f{1}{\pi}\big[-i M^{\rm tr}(\lambda+i0,x_0)\big]d\lambda, \\
&\supp\,\big[d\Omega^{\rm tr}(\cdot,x_0)\big]=\cE  \lb{4.55}
\end{split}
\end{align}
for the trace measure $d\Omega^{\rm tr}(\cdot,x_0)$ associated with $H$.
By \eqref{B.28} this completes the proof of \eqref{4.37}.
\end{proof}

\section{Examples of Reflectionless Schr\"odinger Operators \\ with a
Singular Component in their Spectrum}  \label{s5}

In our final section we will construct examples of reflectionless
Schr\"odinger operators with a singular component in their spectra by
weakening the hypothesis that $\cE\subset\bbR$ is a homogeneous
set. Instead, we will assume the spectrum to be the complement of a
Denjoy--Widom-type domain in $\bbC$.

We start by taking a closer look at Herglotz functions closely related to
reflectionless Schr\"odinger operators.

Let $\gE$ be a closed subset of the nonnegative real axis which we write
as
\begin{equation}
\gE=[0,\infty)\big\backslash\bigcup_{j\in\bbN}(a_j,b_j).  \lb{5.1}
\end{equation}
To a positive measure $d\sigma$ supported on $\gE$,
such that
\begin{equation}
\int_{\gE}\frac{d\sigma(\lambda)}{1+\lambda^2}<\infty,  \lb{5.2}
\end{equation}
we associate the Herglotz function
\begin{equation}
r(z)= \Re(r(i))+\int_{\gE} d\sigma(\lambda)\bigg(\f{1}{\lambda-z}
-\f{\lambda}{1+\lambda^2}\bigg), \quad z\in \bbC_+,  \lb{5.3}
\end{equation}
and recall the property $\Im (r(z))> 0$, $z\in\bbC_+$. In this note we
consider  a specific subclass of such  functions which is directly related
to the diagonal terms of the $2\times 2$ Weyl--Titchmarsh matrix of
reflectionless one-dimensional Schr\"odinger operators and hence of
particular importance to spectral theory. Indeed, motivated by the property
\eqref{2.35} of the diagonal Green's function
$g(\cdot,x)$, $x\in\bbR$, being purely imaginary a.e. on $\cE$, we
introduce the following definition.

\begin{definition}  \lb{d5.1}
Let $r$ be a Herglotz function and $\gE$ be a closed subset of
$[0,\infty)$ of the type \eqref{5.1}. Then $r$ is said to belong to the
class $\cR(\gE)$ if
\begin{equation}
r(\lambda+i0)\in i\bbR \, \text{ for a.e.\ $\lambda\in \gE$.}  \lb{5.4}
\end{equation}
\end{definition}
We note that $r$ is increasing on every interval $(a_j,b_j)$
and thus there exist unique $\mu_j\in[a_j,b_j]$, $j\in\bbN$, such that
\begin{equation}
\begin{cases} r(\lambda) \le 0, & a_j\le \lambda\le \mu_j, \\
             r(\lambda) \ge 0, & \mu_j\le \lambda\le b_j, \\
\end{cases}  \quad j\in\bbN.   \lb{5.5}
\end{equation}
Together with \eqref{5.4} this implies that the argument $\xi(\lambda)$ of
$r(\lambda+i0)$, that is,
\begin{equation}
\xi(\lambda)=\Im(\log (r(\lambda+i0)))  \lb{5.6}
\end{equation}
is defined for a.e.\ $\lambda\in\bbR$, and one obtains the exponential
Herglotz representation of $r$,
\begin{equation}
\log (r(z))=\log (C) + \int_0^\infty d\lambda\, \xi(\lambda) \bigg(
\frac{1}{\lambda-z}-\frac{\lambda}{1+\lambda^2}\bigg),  \lb{5.7}
\end{equation}
where
\begin{equation}
\xi(\lambda)= \begin{cases} 1, & a_j\le \lambda\le \mu_j,
\\ 0, & \mu_j\le \lambda\le b_j, \\ 1/2, & \lambda\in \gE. \end{cases}
\lb{5.8}
\end{equation}
Next, we turn to the case where the measure $d\sigma$ in the Herglotz
representation \eqref{5.3} of $r$ also satisfies the condition
\begin{equation}
\int_{\gE}\frac{d\sigma(\lambda)}{1+|\lambda|}<\infty,  \lb{5.8A}
\end{equation}
and where $r$ in \eqref{5.3} is of the special form
\begin{equation}
r(z)=\int_{\gE} \f{d\sigma(\lambda)}{\lambda-z}, \quad z\in\bbC_+. \lb{5.8B}
\end{equation}
Henceforth, we will assume that the sum of all lengths of the intervals
$(a_j,b_j)$ is finite, that is,
\begin{equation}
\sum_{j\in\bbN}(b_j-a_j)<\infty.   \lb{5.8a}
\end{equation}
In this case the representation \eqref{5.7} can be rewritten in the form
\begin{equation}
r(z)=\f{i}{2 z^{1/2}}
\prod_{j\in\bbN}\frac{(z-\mu_j)}{[(z-a_j)(z-b_j)]^{1/2}},
\lb{5.9}
\end{equation}
where the constant $C$ in \eqref{5.7} is chosen such that
\begin{equation}
r(z)\underset{z\downarrow -\infty}{=} (i/2) z^{-1/2}[1+\oh(1)].
\lb{5.10}
\end{equation}
Under the assumption of a certain regularity of the set
$\gE$,  there exists a specific choice of the zeros $\{\mu_j\}_{j\in\bbN}$ in
\eqref{5.9} which will lead to a purely absolutely continuous measure in the
Herglotz representation \eqref{5.3}. To identify this particular choice of
zeros, we denote by
\begin{equation}
m(z)=\frac i {2
z^{1/2}}\prod_{j\in\bbN}\frac{(z-c_j)}{[(z-a_j)(z-b_j)]^{1/2}}, \quad
z\in\bbC_+,  \lb{5.11}
\end{equation}
a Herglotz function such that
\begin{equation}
-\int_{a_j+0}^{c_j} dz\, m(z) =\int_{c_j}^{b_j-0}dz\, m(z) <\infty, \quad
j\in\bbN,  \lb{5.12}
\end{equation}
and, in addition,
\begin{equation}
\int_{-1}^{-0} dz\, m(z) <\infty.  \lb{5.13}
\end{equation}
Next, we make an even stronger assumption and suppose that for the above
choice of $\{c_j\}_{\in\bbN}$, the following {\it Parreau--Widom-type}
condition (cf.\ \cite{Wi71}, \cite{Wi71a}) also holds,
\begin{equation}
\sum_{j\in\bbN}\int_{c_j}^{b_j-0} dz\, m(z) < \infty.  \lb{5.14}
\end{equation}

\begin{remark} \lb{r5.2}
In terms of the
conformal map \eqref{5.15}, 
we have
\begin{equation}
\Theta(z)=\frac 1{i}\int_0^z  dz\, m(z), \quad z\in\bbC_+, \lb{5.17}
\end{equation}
and $c_j=\Theta^{-1}(u_j+ih_j)$, $j\in\bbN$, that is,
$\sum_{j\in\bbN} h_j <\infty$ becomes \eqref{5.14}.
\smallskip
\end{remark}

\begin{theorem}\label{t5.3}
Suppose the Herglotz function $m$ introduced in \eqref{5.11} satisfies
conditions \eqref{5.12}--\eqref{5.14}. Then the  measure $d\rho$ in the
Herglotz representation of $m$,
\begin{equation}
m(z)=\int_{\gE} \frac{d\rho(\lambda)}{\lambda-z}, \quad z\in\bbC_+,
\lb{5.18}
\end{equation}
is absolutely continuous with respect to Lebesgue measure $d\lambda$ on
$\bbR$.
\end{theorem}
\begin{proof}
We recall that any (nonvanishing) Herglotz function is a function of the
{\it Smirnov class} $N_+(\bbC_+)$ (i.e., a ratio of two
uniformly bounded functions such that the denominator is an outer function).
In fact, any Herglotz function is an outer function (cf.\ \cite[p.\
111]{RR94}). The following maximum principle holds for functions
$g\in N_+(\bbC_+)$: $g(\cdot+i0)\in L^p(\bbR;d\lambda)$ implies $g\in
H^p(\bbC_+)$. Next, we will show that $m(\cdot+i0)/(\cdot+i)\in
L^1(\bbR;d\lambda)$.

Due to \eqref{5.13} and \eqref{5.14} the function $m(\cdot+i0)/(\cdot+i)$ is
integrable on $(-1,0)$ and on $\bigcup_{j\in\bbN}(a_j,b_j)$. Next, since
$m(\lambda+i0)\in i\bbR$ for a.e.\ $\lambda\in\gE$ and thus,
$m(\lambda+i0)=i\pi\rho'_{\rm a.c.}(\lambda)$ a.e.\ on $\cE$, we get
\begin{equation}
\frac 1{\pi i}\int_{\gE} \frac{m(\lambda+i0) d\lambda}{1+|\lambda|}\le
\int_{\gE} \frac{d\rho(\lambda)}{1+|\lambda|}<\infty.  \lb{5.19}
\end{equation}
Finally, by \eqref{5.10} the function $m(\cdot+i0)/(\cdot+i)$ is integrable
on $(-\infty,-1]$. Thus,
\begin{equation}
\int_{\bbR}\frac{|m(\lambda+i0)|d\lambda}{1+|\lambda|}<\infty.  \lb{5.19a}
\end{equation}
Next, following an argument in \cite[p.\ 96--97]{RR94}, we introduce the
function
\begin{equation}
n(z)=\f{1}{2\pi i} \int_{\bbR} \f{m(\lambda+i0) d\lambda}{\lambda-z},
\quad z\in\bbC\backslash\bbR.  \lb{5.19b}
\end{equation}
The $n$ is analytic on $\bbC_+\cup\bbC_-$ and
\begin{align}
n(z)-n(\ol z) &=\f{1}{2\pi i} \int_{\bbR} m(\lambda+i0) d\lambda
\bigg(\f{1}{\lambda-z} -\f{1}{\lambda - \ol z}\bigg)  \no \\
&=\f{y}{\pi} \int_{\bbR} \f{m(\lambda+i0) d\lambda}{(\lambda-x)^2+y^2}
\no \\
&= m(z), \quad  z=x+iy, \; y>0,  \lb{5.19c}
\end{align}
where we used Fatou's theorem in the last step (cf.\ \cite[p.\ 86]{RR94}).
Since $n$ and $m$ are analytic on $\bbC_+$, so is $n(\ol z)$, $z\in\bbC_+$.
But then,
\begin{equation}
\ol n(\ol z)=-\f{1}{2\pi i} \int_{\bbR} \f{\ol{m(\lambda+i0)}
d\lambda}{\lambda-z},
\quad z\in\bbC\backslash\bbR,  \lb{5.19d}
\end{equation}
is also analytic on $\bbC_+$ and hence $n(\ol z)$ must be constant for
$z\in\bbC_+$. Since
\begin{equation}
\lim_{y\uparrow\infty} n(-iy) =0,  \lb{5.19e}
\end{equation}
(using \eqref{5.19a} and applying the dominated convergence theorem), one
concludes
\begin{equation}
n(\ol z)=0, \quad z\in\bbC_+,   \lb{5.19f}
\end{equation}
and hence,
\begin{equation}
\frac 1{2\pi i}\int_{\bbR}\frac{m(\lambda+i0)d\lambda}{\lambda-z}=m(z),
\quad \frac 1{2\pi i}\int_{\bbR}\frac{
m(\lambda+i0) d\lambda}{\lambda-\ol z}=0, \quad  z\in\bbC_+.  \lb{5.20}
\end{equation}
Thus, since $m$ is real-valued on $\bbR\backslash\cE$,
\begin{align}
m(z)&=\f{1}{2\pi i} \int_{\bbR} \f{[m(\lambda+i0)-\ol{m(\lambda+i0)}]
d\lambda}{\lambda -z}  \no \\
&= \f{1}{\pi} \int_{\cE} \f{\Im(m(\lambda+i0)) d\lambda}{\lambda -z} \no \\
&= \int_{\cE}  \f{d \rho(\lambda)}{\lambda -z}, \quad z\in\bbC_+,  \lb{5.20a}
\end{align}
and hence, $d\rho(\lambda)=\pi^{-1}\Im (m(\lambda+i0)) d\lambda$.
\end{proof}

The goal of this section is to show that the above Parreau--Widom-type
condition is insufficient in guaranteeing that all measures related to the
class $\cR(\gE)$ are absolutely continuous. More precisely, we will show that
Herglotz functions $r$ of the type \eqref{5.9}, for a certain distribution of
its zeros $\{\mu_j\}_{j\in\bbN}$, may have a singular component.

To this end, we proceed to constructing a closed set $\gE$ such that the
intervals
$(a_j,b_j)$, $j\in\bbN$, accumulate just at a single point, for instance,
the origin,
\begin{equation}
0<\dots <a_{j+1}<b_{j+1}<a_{j}<b_j<\dots <a_2<b_2<a_1<<b_{1}<\infty.
\lb{5.30}
\end{equation}
In this case, the origin is the only point that may support  a
singular component, $\sigma(\{0\})>0$. A good strategy for accomplishing
this goal is to  put all zeros $\mu_j$ in \eqref{5.9} as far as possible from
the origin, that is, we will choose $\mu_j=b_j$, $j\in\bbN$, in the next
lemma.

\begin{lemma}  \lb{l5.5}
Suppose $\gE$ is of the form \eqref{5.1} with $\{a_j\}_{j\in\bbN}$ and
$\{b_j\}_{j\in\bbN}$ satisfying \eqref{5.30}. In addition, assume that the
product $\prod_{j\in\bbN} (b_{j+1}/a_j)$ converges absolutely, that is,
suppose
\begin{equation}
\sum_{j\in\bbN} [1-(b_{j+1}/a_j)] < \infty.  \lb{5.32}
\end{equation}
Define the Herglotz
function $r_0$ by
\begin{equation}
r_0(z)= \frac i
{2 z^{1/2}}\prod_{j\in\bbN}\bigg(\frac{z-b_j}{z-a_j}\bigg)^{1/2}=
\int_{\gE}\frac{d\sigma_0(\lambda)}{\lambda-z}, \quad z\in\bbC_+.  \lb{5.31}
\end{equation}
Then the measure $d\sigma_0$ has a point mass at $0$, that is,
\begin{equation}
\sigma_0(\{0\})>0.   \lb{5.32a}
\end{equation}
\end{lemma}
\begin{proof}
First, we recall that
\begin{equation}
\sigma_0(\{0\})=\lim_{z\uparrow 0}(-z)r_0(z)  \lb{5.33}
\end{equation}
since $d\sigma_0$ is supported on $\gE \subseteq [0,\infty)$.
For a fixed $n\in\bbN$, we split the infinite product in \eqref{5.31} into
three factors
\begin{equation}
\prod_{j=n+1}^{\infty}\bigg(\frac{z-b_j}{z-a_j}\bigg)^{1/2}
\bigg(\frac{z-b_1}{z-a_n}\bigg)^{1/2}\,
\prod_{j=1}^{n-1}\bigg(\frac{z-b_{j+1}}{z-a_j}\bigg)^{1/2}.    \lb{5.34}
\end{equation}
We note that for $z<0$,
\begin{equation}
\frac{z-b_j}{z-a_j}>1    \lb{5.35}
\end{equation}
and
\begin{equation}
\frac{z-b_{j+1}}{z-a_j}>\frac{b_{j+1}}{a_j}.    \lb{5.36}
\end{equation}
Therefore,
\begin{equation}
-zr_0(z)>\f{1}{2}
\bigg[\frac{z(b_1-z)}{z-a_n}\bigg]^{1/2}\,
\prod_{j=1}^{n-1}\bigg(\frac{b_{j+1}}{a_j}\bigg)^{1/2} \ge
\f{1}{2} \bigg[\frac{z(b_1-z)}{z-a_n}\bigg]^{1/2}
\prod_{j\in\bbN}\bigg(\frac{b_{j+1}}{a_j}\bigg)^{1/2}.    \lb{5.37}
\end{equation}
Passing to the limit $n\to\infty$, we get
\begin{equation}
-zr_0(z)\ge \f{1}{2} (b_1-z)^{1/2}
\prod_{j\in\bbN} (b_{j+1}/a_j)^{1/2}.     \lb{5.38}
\end{equation}
Then, passing to the limit $z\uparrow 0$, one obtains
\begin{equation}
\sigma_0(\{0\})=\lim_{z\uparrow 0}(-z)r_0(z) \ge \f{1}{2} b_1^{1/2}
\prod_{j\in\bbN} (b_{j+1}/a_j)^{1/2}>0.    \lb{5.39}
\end{equation}
\end{proof}

Next, we turn to a construction of intervals symmetric with respect to the
origin and to a given half-axis. Let $[b_1,\infty)$, $b_1>0$, be the given
half-axis.  We now find an interval $[b_2,a_1]$, which is symmetric with
respect to the origin and to the half-axis $[b_1,\infty)$. Starting with
intervals $[-1,-\alpha]\cup[\alpha,1]$ under the conformal map
$\zeta\mapsto z$, $\zeta=-k/z+1$, $0<k<b_1$, one gets
\begin{equation}
\alpha=-\frac k {b_1}+1, \quad
-1=-\frac k  {b_2}+1,  \quad
-\alpha=-\frac k  {a_1}+1.   \lb{5.40}
\end{equation}
In particular
\begin{equation} 
\frac 1  {b_2}  - \frac 1  {a_1}=\frac 1  {b_1}.  \lb{5.41}
\end{equation}

Proceeding iteratively, one proves the following result.

\begin{lemma}\label{l5.6}
Let  $\gE_n=[b_1,\infty)\cup\bigcup_{j=2}^n [b_{j},a_{j-1}]
\subset [0,\infty)$, $n\geq 2$, be a finite system of closed intervals.
Define an interval $[b_{n+1},a_{n}]$ symmetric with respect to the
origin and to the given half-axis $[b_{n},\infty)$. Put
$\gE_{n+1}=[b_1,\infty)\cup\bigcup_{j=2}^{n+1}[b_{j},a_{j-1}]$. Let $u_n$ be
the solution of the Dirichlet problem $\Delta u=0$ in
$\bbC\backslash \gE_{n+1}$ with the boundary conditions
$u|_{[b_{n+1},a_{n}]}=1$, $u|_{\gE_n}=0$. Then,
\begin{equation}
u_n(0)\ge 1/2.   \lb{5.42}
\end{equation}
\end{lemma}
\begin{proof}
Let $u_0$ correspond to the extremal case $\gE= [b_{n},\infty)$. Due to the
symmetry between $[b_{n+1},a_{n}]$ and $[b_{n},\infty)$,
$u_0(0)=1/2$. Note that the difference
$u-u_0$ is nonnegative on the boundary of the domain
$\bbC\backslash([b_{n+1},a_{n}]\cup [b_{n},\infty))$
and therefore in its interior. Thus, $u(0)\ge u_0(0)=1/2$.
\end{proof}

Next we turn to some properties of the corresponding function $\Theta$
in \eqref{5.17}: Let $\Theta_n$ be the corresponding functions related to the
system of intervals $\gE_n$ introduced in Lemma \ref{l5.6}. We recall that
its imaginary part,
\begin{equation}
\omega_n (z)=\Im (\Theta_n(z)), \quad z\in \bbC\backslash \gE_n, \lb{5.42a}
\end{equation}
is a positive single-valued harmonic function on the domain
$\bbC\backslash \gE_n$ such that $\omega_n|_{\gE_n}=0$ with the only
singularity
\begin{equation}
\omega_n(z)\underset{z\to -\infty}{=}\Im\big(z^{1/2}\big)[1+\oh(1)].
\lb{5.43}
\end{equation}
Thus, $\omega_{n}(z)-\omega_{n+1}(z)$ is a uniformly bounded
harmonic function (the singularity at infinity cancels). Since
\begin{equation}
\omega_{n}(\lambda \pm i0)-\omega_{n+1}(\lambda \pm i0)\ge 0, \quad
\lambda \in \gE_{n+1},  \lb{5.44}
\end{equation}
the same inequality also holds inside the domain
\begin{equation}
\omega_{n}(z)-\omega_{n+1}(z)\ge 0,\quad z\in \bbC\backslash \gE_{n+1}.
\lb{5.45}
\end{equation}

We also note that $\omega_n(\lambda)$ decreases in $(-\infty,b_n]$.
In fact, $\Theta_n(\lambda)=i\omega_n(\lambda)$  for all $\lambda\in
\bbR\backslash \gE_n$.

\begin{theorem} [Construction of a Denjoy--Widom-type domain
$\bbC\backslash\gE$]
\lb{t5.7}
${}$\\ There exists a closed set $\gE$ of the form \eqref{5.1} with
$\{a_j\}_{j\in\bbN}$  and $\{b_j\}_{j\in\bbN}$ satisfying \eqref{5.30}, such
that for the choice $\mu_j=b_j$, $j\in\bbN$, the corresponding Herglotz
function $r_0$ defined in \eqref{5.31} leads to a Herglotz representation,
where the associated measure $d\sigma_0$ has a point mass at $0$, whereas the
choice $\mu_j=c_j$, $j\in\bbN$, leads to a Herglotz function $m$ defined in
\eqref{5.11} with corresponding measure $d\rho$ in \eqref{5.18} purely
absolutely continuous with respect to Lebesgue measure.
\end{theorem}
\begin{proof}
We follow closely the construction given by Hasumi \cite[p.\ 215--218]{Ha83}.
Starting with $[b_1,\infty)$, $b_1>0$, we will construct a system of
intervals accumulating at the origin such that on one hand the series
\eqref{5.32} converges and such that, on the other hand, the measure $d\rho$
(see Theorem \ref{t5.3}) is absolutely continuous.

First we claim that the following three conditions can be
simultaneously satisfied by choosing $\ell_n$, $0<\ell_n<1/2$, $n\in\bbN$,
sufficiently small:
\begin{equation}
\begin{aligned}
{(a)}&\ &&b_{n+1}=\ell_n b_n,\\
{(b)}&\ &&\frac 1{b_{n+1}}-\frac 1{a_{n}}=\frac 1{b_{n}},  \lb{5.46} \\
{(c)}&\ &&\inf_{n\in\bbN}\{\omega_n(z)\,|\, z\in [b_{n+1},a_{n}]\}\ge \frac
1 2 \omega_n(0).
\end{aligned}
\end{equation}
Indeed, making $\ell_n$ smaller, we get an interval
$ [b_{n+1},a_{n}]$ that approaches the origin. Thus, conditions $(a)$ and
$(b)$ hold. In addition, by continuity of $\omega_n$, condition $(c)$ is
also satisfied.

We note that under this construction, condition \eqref{5.32} is of the form
\begin{equation}
\sum_{j\in\bbN} b_{j+1}\left(\frac 1 {b_{j+1}}-\frac 1{ a_{j}}\right)=
\sum_{j\in\bbN} b_{j+1}\frac {1} {b_{j}}=\sum_{j\in\bbN} \ell_j. \lb{5.47}
\end{equation}
Thus, if necessary, by making $\ell_n$ even smaller
in a such way that $\sum_{n\in\bbN} \ell_n<\infty$, condition
\eqref{5.32} is satisfied.

As a next step we seek an inductive estimate on $\omega_n(0)$.
Due to condition $(c)$ in \eqref{5.46} one infers
\begin{equation}
\omega_n(z)-\omega_{n+1}(z)-\frac 1 2 \omega_n(0) u_n(z)\ge 0, \quad
z\in \gE_{n+1},   \lb{5.48}
\end{equation}
where  $u_n$ is the solution
of the Dirichlet problem $\Delta u =0$ in $\bbC\backslash
\gE_{n+1}$ with the boundary conditions $u|_{ [b_{n+1},a_{n}]}=1$,
$u|_{\gE_n}=0$. Since the same inequality holds inside the
domain, using \eqref{5.42} one gets,
\begin{equation}
\omega_n(0)-\omega_{n+1}(0)\ge \frac 1 2 \omega_n(0) u_n(0)\ge
\frac 1 4 \omega_n(0).  \lb{5.49}
\end{equation}
Thus, $\{\omega_n(0)\}_{n\in\bbN}$ form a sequence that is dominated
by a geometric progression,
\begin{equation}
\omega_{n+1}(0)\le\frac 3 4  \omega_{n}(0)\le
\left( \frac 3 4\right)^{n+1}\omega_{0}(0), \quad n\in\bbN_0.  \lb{5.50}
\end{equation}

Now we are in position to show that all assumptions of
Theorem \ref{t5.3} are satisfied.  We note that
\begin{equation}
\omega(z)=\lim_{n\uparrow\infty}\omega_n(z),   \lb{5.51}
\end{equation}
and thus
\begin{equation}
m(z)=\lim_{n\uparrow \infty} m_n(z), \quad
\Theta(z)=\lim_{n\uparrow\infty}\Theta_n(z).   \lb{5.52}
\end{equation}

Due to the monotonicity property with respect to $n$ one gets
\begin{equation}
\omega(z)\le\omega_n(z),  \lb{5.53}
\end{equation}
in particular,
\begin{equation}
\omega(c_n)\le\omega_n(c_n),   \lb{5.54}
\end{equation}
and due to the monotonicity property of $\omega_n(\cdot)$ on $(-\infty,b_n]$
one obtains (cf.\ \eqref{5.50})
\begin{equation}
\omega_n(c_n)\le\omega_n(0) \leq \bigg(\f{3}{4}\bigg)^n\omega_0(0), \quad
n\in\bbN.    \lb{5.55}
\end{equation}
Thus,
\begin{equation}
\sum_{n\in\bbN}\omega(c_n)<\infty,   \lb{5.56}
\end{equation}
that is, \eqref{5.14} is satisfied. Finally,
\begin{align}
\begin{split}
\int_{-1}^{-0} dz\, m(z)& = \int_{-1}^{-0} dz \lim_{n\uparrow\infty} m_n(z)
\le\lim_{n\uparrow\infty}\int_{-1}^0 dz\, m_n(z)     \lb{5.57}  \\
& = \lim_{n\uparrow\infty}[\omega_n(-1)-\omega_n(0)]\le\omega(-1)<\infty.
\end{split}
\end{align}
\end{proof}

\begin{remark}  \lb{r5.7a}
By inspection (cf.\ Theorem \ref{tB.5}), the measure $d\sigma_0$ in the
Herglotz representation of $r_0$ in Theorem
\ref{t5.7} has purely absolutely continuous spectrum away from zero. This is
in agreement with a result of Aronszajn and Donoghue \cite{AD56} recorded in
Theorem \ref{tB.4}\,$(iv)$ since $\xi(\lambda)=1/2$ on
$\bigcup_{j\in\bbN} (b_{j+1},a_j)\cup(b_1,\infty)$. This result by Aronszajn
and Donoghue applies to open intervals and hence does not exclude an
eigenvalue at
$\lambda=0$ (i.e., $\sigma_0(\{0\})>0$) as constructed in Theorem \ref{t5.7}.
\end{remark}

In order to apply this construction to one-dimensional Schr\"odinger
operators we next recall sufficient conditions for $m_+(z)=m_+(z,0)$,
$z\in\bbC_+$, to be the half-line Weyl--Titchmarsh function associated with
a Schr\"odinger operator on $[0,\infty)$ in terms of the corresponding
measure $d\omega_+$ in the Herglotz representation of $m_+$. Based on the
classical inverse spectral theory approach due to Gelfand and Levitan
\cite{GL51}, the following result discussed in Thurlow \cite{Th79} (see
also, \cite[Sects.\ 2.5, 2.9]{Le87}, \cite{LG64}, \cite[Sect.\ 26.5]{Na68},
\cite{Ro60}) describes sufficient conditions for a monotonically
nondecreasing  function $\omega_+$ on $\bbR$ to be the spectral function of a
half-line Schr\"odinger operator $H_+$ in $L^2([0,\infty);dx)$ with a
Dirichlet boundary condition at $x=0$.

\begin{theorem} [\cite{Th79}] \lb{t5.8}
Let $\omega_+$ be a monotonically nondecreasing function on
$\bbR$ satisfying the following two conditions. \\
\noindent $(i)$
Whenever $f\in C([0,\infty))$ with compact support
contained in $[0,\infty)$ and
\begin{equation}
\int_\bbR d\omega_+(\lambda)\,|F(\lambda)|^2 =0,
    \text{ then $f=0$,} \lb{5.59}
\end{equation}
where
\begin{equation}
F(\lambda)=\int_{0}^{\infty}
dx\,\f{\sin(\lambda^{1/2}(x-x_0))}{\lambda^{1/2}}f(x), \quad
\lambda\in\bbR. \lb{5.60}
\end{equation}
\noindent $(ii)$ Define
\begin{equation}
\wti\omega_+(\lambda)=\begin{cases}
\omega_+(\lambda)-\f{2}{3\pi}\lambda^{3/2}, & \lambda\geq 0, \\
\omega_+(\lambda), & \lambda<0 \end{cases} \lb{5.61}
\end{equation}
and assume the limit
\begin{equation}
\lim_{R\uparrow\infty}\int_{-\infty}^R d\wti\omega_+(\lambda) \,
\f{\sin(\lambda^{1/2} x)}{\lambda^{1/2}}= \Phi(x), \quad x\geq 0, \lb{5.62}
\end{equation}
exists with $\Phi\in L^\infty([0,R];dx)$ for all $R>0$.
Moreover, suppose that for some $r\in\bbN_0$, $\Phi^{(r+1)}\in
L^1([0,R];dx)$ for all $R>0$, and $\lim_{x\downarrow 0}\Phi(x)=0$. \\
Then $d\omega_+$ is the spectral measure of a self-adjoint
Schr\"odinger operator $H_+$ in $L^2([0,\infty);dx)$ associated
with the differential expression
$L_+=-d^2/dx^2 +V_+$, $x>0$, with a Dirichlet boundary condition
at $x=0$, a self-adjoint boundary condition at $\infty$ $($if
necessary$)$, and a real-valued potential coefficient $V_+$ satisfying
$V_+^{(r)}\in L^1([0,R];dx)$ for all $R>0$.
\end{theorem}

Since the analogous result applies to Schr\"odinger operators on
the half-line $(-\infty,0]$, we omit the details.

\begin{remark}  \lb{r5.9}
We add two observations to Theorem \ref{t5.8}: First, whenever the
points of increase of $\omega_+$ have a finite limit point, in particular, if
$\omega_+$ is strictly increasing on an interval, then condition $(i)$ in
Theorem \ref{t5.8} is satisfied. This follows from the fact that $F$ is an
entire function of order $1/2$ which cannot vanish on a set with a finite
accumulation point, or even on an interval, without vanishing identically.
The latter yields $f=0$ by the inverse sine transform. Secondly, if the
minimal operator $H_{+,\rm min}$ associated with $L_+$ and a Dirichlet
boundary condition at $x=0$ is bounded from below, then $L_+$ is in the
limit point case at infinity (i.e., $H_{+,\rm min}$ is essentially
self-adjoint) and no boundary condition is needed at $x=\infty$ (this is
again based on the well-known result of Hartman \cite{Ha48} (see also
Rellich \cite{Re51} and \cite{Ge93}). Here $H_{+,\rm min}$ in
$L^2([0,\infty);dx)$ is defined as
\begin{align}
&H_{+,\rm min}f= L_+ f,  \no  \\
&f\in\dom(H_{+,\rm min})=\big\{g\in L^2([0,\infty);dx)\,|\, \supp\,(g)
\text{ compact; }   \lb{5.63} \\
& \hspace*{3.7cm} \lim_{\varepsilon\downarrow 0}g(\varepsilon)=0; \, L_+ g
\in L^2([0,\infty);dx)\big\}.  \no
\end{align}
In the concrete situations we have in mind
below, both observations clearly apply (as all spectra involved are bounded
from below and contain intervals of absolutely continuous spectrum) and
hence we will focus exclusively on condition $(ii)$ in Theorem \ref{t5.8}
and utilize the fact that all (half-line and full-line) Schr\"odinger
operators considered in this paper are in the limit point case at
$\pm\infty$.
\end{remark}

In connection with Theorem \ref{t5.8}\,$(ii)$, the following special case
of Lemma 8.3.2 in \cite{Le87} turns out to be useful.

\begin{lemma} [\cite{Le87}, Lemma 8.3.2] \lb{l5.10}
Let $f\in L^1([1,\infty);dp)$ and suppose for some $R>1$ and $p>R$, $f$
permits the asymptotic expansion
\begin{equation}
f(p)\underset{p\uparrow\infty}{=} \sum_{k=1}^N f_k p^{-2k}
+ \Oh\big(p^{-(2N+2)}\big) \, \text{ for all $N\in\bbN$}   \lb{5.64}
\end{equation}
for some coefficients $\{f_k\}_{k\in\bbN}\subset\bbC$. Then the function
\begin{equation}
\phi(x)=\int_1^{\infty} dp\, \f{\sin(px)}{p} f(p), \quad x\geq 0,   \lb{5.65}
\end{equation}
satisfies
\begin{equation}
\phi \in C^{\infty} ((0,\infty)), \quad \phi^{(m)} \in
L^{\infty}([0,\infty);dx)  \, \text{ for all $m\in\bbN_0$,}
\quad \lim_{x\downarrow 0}\phi(x)=0,  \lb{5.66}
\end{equation}
and
\begin{equation}
\lim_{x\uparrow\infty} \phi^{(m)} (x) =0 \, \text{ for all $m\in\bbN_0$.}
\lb{5.67}
\end{equation}
\end{lemma}

Next, we construct a one-dimensional Schr\"odinger operator $H_{r_0}$
with spectrum the set $\gE$ in Theorem \ref{t5.7}, whose
diagonal Green's function $g_{r_0}(\cdot,0)$ coincides
with the function $r_0$ in Theorem \ref{t5.7}, and whose point spectrum is
nonempty.

\begin{theorem}  \lb{t5.11}
Let $\gE$ be a closed set of the form \eqref{5.1} with $\{a_j\}_{j\in\bbN}$
and $\{b_j\}_{j\in\bbN}$ satisfying \eqref{5.30} and let $r_0$ be the
Herglotz function defined in \eqref{5.31} with associated measure
$d\sigma_0$ having a point mass at $0$ as constructed in Theorem \ref{t5.7}.
Then, there exists a self-adjoint Schr\"odinger operator
$H_{r_0}$ in $L^2(\bbR;dx)$ with real-valued reflectionless potential
coefficient $V_{r_0}$ satisfying
\begin{equation}
V_{r_0} \in C^{\infty}(\bbR\backslash\{0\}),  \quad
V_{r_0}^{(m)} \in L^1_{\loc}(\bbR;dx) \, \text{ for all $m\in\bbN_0$}
\lb{5.68}
\end{equation}
with spectral properties
\begin{equation}
    \sigma(H_{r_0})=\sigma_{\rm ac}(H_{r_0})=\gE, \quad
\sigma_{\rm pp}(H_{r_0})=\{0\},  \quad \sigma_{\rm sc}(H_{r_0})=\emptyset.
\lb{5.69}
\end{equation}
In particular, $H_{r_0}$ has nonempty singular component $($a zero
eigenvalue$)$ in its spectrum.
\end{theorem}
\begin{proof}
First we note that if $H$ is a reflectionless Schr\"odinger operator with
associated diagonal Green's function $g(\cdot,x)\in \cR(\gE)$, $x\in\bbR$,
then also
\begin{equation}
-g(\cdot,x)^{-1} \in \cR(\gE), \quad x\in\bbR.  \lb{5.70}
\end{equation}
Moreover, by \eqref{2.7}, we also have
\begin{equation}
-g(z,x)^{-1} = m_+(z,x) - m_-(z,x), \quad (z,x)\in \bbC_+\times \bbR,
\lb{5.71}
\end{equation}
with $m_\pm (\cdot,x)$ the half-line Weyl--Titchmarsh functions of $H$.
Relation \eqref{5.71} will be helpful in introducing $m_{\pm,r_0}(z,0)$ for
$H_{r_0}$. To this end we define
\begin{equation}
g_{r_0}(z,0)= r_0(z) = \frac{i}{2 z^{1/2}}
\prod_{j\in\bbN}\bigg(\frac{z-b_j}{z-a_j}\bigg)^{1/2}=
\int_{\gE}\frac{d\sigma_0(\lambda)}{\lambda-z}, \quad z\in\bbC_+  \lb{5.72}
\end{equation}
and write
\begin{align}
-g_{r_{0}}(z,0)^{-1}&=-r_0(z)^{-1} =m_{+,r_{0}}(z)-m_{-,r_{0}}(z)  \no \\
&=2iz^{1/2}\prod_{j\in\bbN} \bigg(\f{z-a_j}{z-b_j}\bigg)^{1/2}  \lb{5.72a} \\
&= \Re(-g(i)^{-1})+\int_{\bbR} d\omega_{r_0} (\lambda)
\bigg(\f{1}{\lambda-z}-\f{\lambda}{1+\lambda^2}\bigg), \quad z\in\bbC_+,
\lb{5.73}
\end{align}
and
\begin{equation}
\pm  m_{\pm, r_{0}}(z)=\pm \Re(m_{\pm, r_{0}}(i))+\int_{\bbR}
d\omega_{\pm, r_{0}}(\lambda)
\bigg(\f{1}{\lambda-z}-\f{\lambda}{1+\lambda^2}\bigg), \quad z\in\bbC_+.
\lb{5.74}
\end{equation}
Here $m_{\pm,r_{0}}(z)=m_{\pm,r_{0}}(z,0)$ will be chosen next to represent
the half-line Weyl--Titchmarsh functions of $H_{r_0}$. Because of
$g_{r_0}(\cdot,0)\in \cR(\gE)$, one infers
\begin{equation}
m_{+,r_{0}}(\lambda+i0) = \ol{m_{-,r_{0}}(\lambda+i0)}
\, \text{ for a.e.\ $\lambda\in\gE$}  \lb{5.75}
\end{equation}
and hence,
\begin{equation}
-g_{r_0}(\lambda+i0,0)^{-1} = \pm 2i \, \Im(m_{\pm,r_{0}}(\lambda+i0))
\, \text{ for a.e.\ $\lambda\in\gE$.}  \lb{5.76}
\end{equation}
By \eqref{5.72}--\eqref{5.74}, \eqref{5.76}, \eqref{A.14}, \eqref{B.11},
\eqref{B.13}, and the fact that any singular continuous component of a
measure must be supported on an uncountable set, one concludes
\begin{align}
& \supp(d\omega_{\pm,r_{0},\rm ac})=\supp(d\omega_{r_0,\rm ac}) =\gE,  \quad
d\omega_{\pm,r_{0},\rm ac}=\f{1}{2} d\omega_{r_0, \rm ac},  \lb{5.77} \\
& d\omega_{\pm,r_{0},\rm sc}=d\omega_{r_0, \rm sc} =0.  \lb{5.78}
\end{align}
Moreover, since
\begin{equation}
\lim_{\varepsilon\downarrow 0} (- \varepsilon) g(-\varepsilon,0)=
\lim_{\varepsilon\downarrow 0} \varepsilon r_0(-\varepsilon) =
\sigma_0(\{0\})>0   \lb{5.79}
\end{equation}
by assumption (cf.\ Theorem \ref{t5.7}), and the Herglotz property of
$\pm m_{\pm,r_0}$ implies that
\begin{equation}
\lim_{\varepsilon\downarrow 0} (\pm \varepsilon)
m_{\pm,r_0}(-\varepsilon) \geq 0,   \lb{5.79a}
\end{equation}
one infers
\begin{equation}
\lim_{\varepsilon\downarrow 0} (- \varepsilon) g(-\varepsilon,0)^{-1} =
\lim_{\varepsilon\downarrow 0} (\pm \varepsilon)
m_{\pm, r_{0}}(-\varepsilon) =0  \lb{5.80}
\end{equation}
and hence
\begin{equation}
\omega_{r_0, \rm pp} (\{0\}) = \omega_{\pm,r_{0},\rm pp}(\{0\})=0.  \lb{5.81}
\end{equation}
Since by \eqref{5.72} and \eqref{B.14},
\begin{equation}
\supp(\omega_{r_0, \rm pp})\cap (\bbR\backslash\{0\}) =
\supp(\omega_{\pm,r_{0},\rm pp})\cap (\bbR\backslash\{0\})=\emptyset,
\lb{5.82}
\end{equation}
one finally concludes
\begin{equation}
d\omega_{r_0, \rm pp} = d\omega_{\pm,r_{0},\rm pp} =0.
\lb{5.83}
\end{equation}
Next we introduce
\begin{align}
h_{r_0}(z,0)&=\f{m_{-,r_0}(z)m_{+,r_0}(z)}{m_{-,r_0}(z)-m_{+,r_0}(z)}
=g_{r_0}(z,0)m_{-,r_0}(z)m_{+,r_0}(z)  \no \\
&=\f{i}{2 z^{1/2}} \prod_{j\in\bbN} \bigg(\f{z-b_j}{z-a_j}\bigg)^{1/2}
m_{-,r_0}(z)m_{+,r_0}(z)  \lb{5.84}  \\
&=\Re(h_{r_0}(i)) + \int_{\bbR} d\rho_0(\lambda)\bigg(\f{1}{\lambda-z}
-\f{\lambda}{1+\lambda^2}\bigg), \quad z\in\bbC_+.  \lb{5.85}
\end{align}
A comparison of \eqref{5.72a} and \eqref{5.84} then proves that $d\rho_0$
has no pure points on $(0,\infty)$. Moreover, since $h_{r_0}$ and $g_{r_0}$
are Herglotz, and the measure $d\sigma_0$ associated with $g_{r_0}$ has a
point mass at $0$, considering the expression
\begin{equation}
\f{\varepsilon h_{r_0}(-\varepsilon,0)}{\varepsilon g_{r_0}(-\varepsilon,0)}
= m_{+,r_0}(-\varepsilon)m_{-,r_0}(-\varepsilon)   \lb{5.91}
\end{equation}
implies that
\begin{equation}
\lim_{\varepsilon\downarrow 0} m_{+,r_0}(-\varepsilon)m_{-,r_0}(-\varepsilon)
\geq 0 \, \text{ exists finitely.}  \lb{5.92}
\end{equation}
In particular,
\begin{equation}
\rho_0(\{0\})>0 \, \text{ if and only if } \,
\lim_{\varepsilon\downarrow 0} m_{+,r_0}(-\varepsilon)m_{-,r_0}(-\varepsilon)
>    0.  \lb{5.93}
\end{equation}
Summing up, one concludes
\begin{equation}
\supp(d\rho_0)=\supp(d\rho_{0,\rm ac})=\gE, \quad
\supp(d\rho_{0,\rm pp}) \subseteq \{0\},  \quad
d\rho_{0,\rm sc}=0.  \lb{5.94} \\
\end{equation}
In order to apply Theorem \ref{t5.8} to $m_{+,r_0}$, one considers the
function
\begin{equation}
\Phi(x)=\lim_{R\uparrow\infty} \int_0^R d\ti \omega_+(\lambda)
\f{\sin(\lambda^{1/2}x)}{\lambda^{1/2}},  \quad x\geq 0,  \lb{5.95}
\end{equation}
where
\begin{equation}
d\ti \omega_+(\lambda) = \pi^{-1}
\lambda^{1/2}\left[\prod_{j\in\bbN}
\bigg(\f{\lambda-a_j}{\lambda-b_j}\bigg)^{1/2}\chi_{\gE}(\lambda) -1\right]
d\lambda  \lb{5.96}
\end{equation}
with $\chi_{\gE}$ the characteristic function of the set $\gE$.
Splitting the integral in \eqref{5.95} over the intervals $[0,b_1]$
and $[b_1,\infty)$,
\begin{equation}
\Phi_1(x) = \int_0^{b_1} d\ti \omega_+(\lambda)
\f{\sin(\lambda^{1/2}x)}{\lambda^{1/2}},  \quad
\Phi_2(x)=\lim_{R\uparrow\infty} \int_{b_1}^R d\ti \omega_+(\lambda)
\f{\sin(\lambda^{1/2}x)}{\lambda^{1/2}}, \quad x\geq 0,  \lb{5.97}
\end{equation}
clearly,
\begin{equation}
\Phi_1\in C^{\infty}((0,\infty)), \quad \Phi_1^{(m)} \in
L^{\infty}([0,\infty);dx)) \, \text{ for all $m\in\bbN_0$,} \quad
\lim_{x\downarrow 0} \Phi_1(x)=0.  \lb{5.98}
\end{equation}
Applying Lemma \ref{l5.10} to $\Phi_2$ then yields all the properties in
\eqref{5.98} for $\Phi_2$ as well. Thus, Theorem \ref{t5.8} shows that
$m_{+,r_0}$ is the half-line Weyl--Titchmarsh function of a Schr\"odinger
operator $H_{+,r_0}$ in $L^2((0,\infty);dx)$ with a Dirichlet boundary
condition at $x=0$ and a real-valued potential coefficient $V_{+,r_0}$
satisfying
\begin{equation}
V_{+,r_0}^{(m)} \in L^1 ([0,R];dx)  \, \text{ for all $m\in\bbN_0$ and all
$R>0$.}  \lb{5.99}
\end{equation}
Analogous considerations then prove that $m_{-,r_0}$ is the half-line
Weyl--Titchmarsh function of a Schr\"odinger operator $H_{-,r_0}$ in
$L^2((-\infty,0];dx)$ with a Dirichlet boundary condition at $x=0$ and a
real-valued potential coefficient $V_{-,r_0}$ satisfying
\begin{equation}
V_{-,r_0}^{(m)} \in L^1 ([-R,0];dx)  \, \text{ for all $m\in\bbN_0$ and all
$R>0$.}  \lb{5.100}
\end{equation}
Given $V_{\pm,r_0}$ we next introduce
\begin{equation}
V_{r_0}(x)=\begin{cases} V_{+,r_0}(x), & x>0, \\ V_{-,r_0}(x), & x<0,
\end{cases}  \lb{5.101}
\end{equation}
and note that $V_{r_0}$ satisfies the properties in \eqref{5.68}. Finally,
introducing the differential expression
\begin{equation}
L_{r_0}=-d^2/dx^2+V_{r_0}(x), \quad x\in\bbR,   \lb{5.103}
\end{equation}
and denoting by $H_{r_0}$ the corresponding self-adjoint
realization of $L_{r_0}$ in $L^2(\bbR;dx)$,
\begin{align}
\begin{split}
& H_{r_0}f=L_{r_0} f, \\
& f\in\dom(H_{r_0})=\big\{g\in L^2(\bbR;dx)\,|\, g,g' \in AC_{\loc}(\bbR); \,
L_{r_0} g \in L^2(\bbR;dx)\big\},  \lb{5.104}
\end{split}
\end{align}
a study of the trace measure of $H_{r_0}$,
\begin{equation}
d\Omega^{\tr}_{r_0}(\lambda)=d\sigma_0(\lambda) + d\rho_0(\lambda) \lb{5.105}
\end{equation}
then yields the spectral properties of $H_{r_0}$ in \eqref{5.69}.
\end{proof}

For simplicity we singled out the eigenvalue $0$ in $H_{r_0}$ (which is of
course correlated with the point mass at $0$ of $d\sigma_0$ in
Lemma \ref{l5.5}). A similar construction leads to an eigenvalue of a
reflectionless Schr\"odinger operator at a point $\lambda_0>0$.

For an entirely different construction of reflectionless tridiagonal
matrices (Jacobi operators) on the lattice $\bbZ$ with empty singular
continuous spectra but possibly countably many accumulation points in the
set of eigenvalues as well as in the set of boundary points of intervals of
absolutely continuous spectrum, we refer to
\cite{GKT96}.

\appendix
\section{Essential Closures of Sets \\ and Essential
Supports of Measures}
\lb{A}
\renewcommand{\theequation}{A.\arabic{equation}}
\renewcommand{\thetheorem}{A.\arabic{theorem}}
\setcounter{theorem}{0}
\setcounter{equation}{0}

The following material on essential closures of essential supports of
absolutely continuous measures is probably well-known, but we found no
comprehensive treatment in the literature and hence decided to collect
the relevant facts in this appendix.

For basic facts on measures on $\bbR$ relevant to this appendix we
refer, for instance, to \cite{Ar57}--\cite{AD64}, \cite[p.\ 179]{CFKS87},
\cite{DSS94}, \cite{Gi84}--\cite{GP87}, \cite[Sect.\ V.12]{PF92},
\cite[p.\ 140--141]{RS78}, \cite{Si95}, \cite{SW86}. All measures in this
appendix will be assumed to be nonnegative without explicitly stressing this
fact again.

Since Borel and Borel--Stieltjes measures are incomplete (i.e., not any
subset of a set of measure zero is measurable) we will enlarge the Borel
$\sigma$-algebra to obtain the complete Lebesgue and Lebesgue--Stieltjes
measures. We recall the standard Lebesgue decomposition of a measure
$d\mu$ on $\bbR$ with respect to Lebesgue measure $dx$ on $\bbR$,
\begin{align}
d\mu&=d\mu_{\rm ac}+d\mu_{\rm s}=d\mu_{\rm ac}+d\mu_{\rm sc} +d\mu_{\rm
pp}, \lb{A.1} \\
d\mu_{\rm ac}&=fdx, \quad 0\leq f\in L^1_{\rm loc}(\bbR;dx), \lb{A.2}
\end{align}
where $d\mu_{\rm ac}$, $d\mu_{\rm sc}$, and $d\mu_{\rm pp}$ denote the
absolutely continuous, singularly continuous, and pure point parts of
$d\mu$, respectively. The Lebesgue measure of a Lebesgue measurable set
$\Omega\subset\bbR$ will be denoted by $|\Omega|$.

In the following, the Lebesgue measure of a set $S$ will be denoted by
$|S|$ and all sets whose $\mu$-measure or Lebesgue measure is considered
are always assumed to be Lebesgue--Stieltjes or Lebesgue measurable, etc.

\begin{definition} \lb{dA.1}
Let $d\mu$ be a Lebesgue--Stieltjes measure and suppose $S$ and $S'$
are $\mu$-measurable. \\
$(i)$ $S$ is called a {\it support} of $d\mu$ if
$\mu(\bbR\backslash S)=0$. \\
$(ii)$ The smallest closed support of $d\mu$ is called the {\it
topological support} of $d\mu$ and denoted by $\supp \, (d\mu)$. \\
$(iii)$ $S$ is called an {\it essential} $($or {\it minimal\,$)$ support}
of $d\mu$ (relative to Lebesgue measure $dx$ on $\bbR$) if
$\mu(\bbR\backslash S)=0$, and $S'\subseteq S$ with $S'$
$|\cdot|$-measurable, $\mu(S')=0$ imply $|S'|=0$.
\end{definition}

\begin{remark} \lb{rA.2}
Item $(iii)$ in Definition \ref{dA.1} is equivalent to \\
$(iii')$ $S$ is called an {\it essential} $($or {\it minimal\,$)$ support}
of $d\mu$ (relative to Lebesgue measure $dx$ on $\bbR$) if
$\mu(\bbR\backslash S)=0$, and $S'\subseteq S$, $\mu(\bbR\backslash S')=0$
imply $|S\backslash S'|=0$.
\end{remark}

\begin{lemma} [\cite{Gi84}]\lb{lA.3}
Let $S, S' \subseteq\bbR$ be $\mu$- and $|\cdot|$-measurable. Define the
relation $\sim$ by $S\sim S'$ if
\begin{equation}
\mu(S\Delta S')=|S\Delta S'|=0  \lb{A.3}
\end{equation}
$($where $S\Delta S'=(S\backslash S')\cup (S'\backslash S)$$)$. Then
$\sim$ is an equivalence relation. Moreover, the set of all essential
supports of $d\mu$ is an equivalence class under $\sim$.
\end{lemma}

\begin{example} \lb{eA.4}
Let $d\mu_{\rm pp}$ be a pure point measure and
\begin{equation}
\mu_{\rm pp}(\{x\})=\begin{cases} c(x)>0, & x\in [0,1]\cap \bbQ, \\
0, & x\in [0,1]\backslash\bbQ \text{ or } |x|>1.
\end{cases}
\lb{A.4}
\end{equation}
Then,
\begin{equation}
\supp \, (d\mu_{\rm pp})=[0,1]. \lb{A.5}
\end{equation}
However, since $[0,1]\cap \bbQ$ is an essential support of $d\mu_{\rm pp}$
and since $|[0,1]\cap \bbQ|=0$, also
\begin{equation}
|S_{\mu_{\rm pp}}|=0  \lb{A.6}
\end{equation}
for any other essential support $S_{\mu_{\rm pp}}$ of $d\mu_{\rm pp}$.
\end{example}

\begin{remark} \lb{rA.5}
If $d\mu =d\mu_{\rm ac}$, then $|S\Delta S'|=0$ implies $\mu(S\Delta
S')=0$ and hence any two essential supports of $d\mu$ differ at most by
sets of Lebesgue measure zero. Indeed, one can use the following,
\begin{align}
&S_1=(S_2\cap S_1)\cup (S_1\backslash S_2), \quad
S_2=(S_1\cap S_2)\cup (S_2\backslash S_1),  \lb{A.7} \\
&S_1\cup(S_2\backslash S_1)=S_2\cup(S_1\backslash S_2)  \lb{A.8}
\end{align}
for any subsets $S_j\subset\bbR$, $j=1,2$.
\end{remark}

\begin{definition} \lb{dA.6}
Let $A\subset \bbR$ be Lebesgue measurable. Then the {\it essential
closure} ${\ol A}^e$ of $A$ is defined as
\begin{equation}
{\ol A}^e =\{x\in\bbR\,|\, \text{for all $\varepsilon>0$: }
|(x-\varepsilon,x+\varepsilon)\cap A|>0\}. \lb{A.9}
\end{equation}
\end{definition}

The following is an immediate consequence of Definition \ref{dA.6}.

\begin{lemma} \lb{lA.7}
Let $A, B, C \subset\bbR$ be Lebesgue measurable. Then,
\begin{align}
& \text{$(i)$ If $A\subseteq B$ then ${\ol A}^e\subseteq {\ol B}^e$.}
\lb{A.10A} \\
& \text{$(ii)$ If $|A|=0$ then ${\ol A}^e=\emptyset$.}  \lb{A.10} \\
& \text{$(iii)$ If $A=B\cup C$ with $|C|=0$, then ${\ol A}^e={\ol B}^e$.}
\lb{A.10B}
\end{align}
\end{lemma}

\begin{example} \lb{eA.8} ${}$ \\
$(i)$ Consider $d\mu_{\rm pp}$ in Example \ref{eA.4}. Let
$S_{\mu_{\rm pp}}$ be any essential support of $d\mu_{\rm pp}$. Then
${\ol {S_{\mu_{\rm pp}}}}^e=\emptyset$ by \eqref{A.10}. \\
$(ii)$ Consider $A=[0,1]\cup\{2\}$. Then ${\ol A}^e =[0,1]$.
\end{example}

\begin{lemma} \lb{lA.9}
Let $A\subseteq\bbR$ be Lebesgue measurable. Then,
\begin{align}
& \text{$(i)$ \,\,\, ${\ol A}^e$ is a closed set.}  \lb{A.10a} \\
& \text{$(ii)$ \, ${\ol A}^e \subseteq \ol A$.}  \lb{A.10b}  \\
& \text{$(iii)$ \hspace*{.01mm}
$\ol{\big({\ol A}^e\big)}^e \subseteq {\ol A}^e$.} \lb{A.10c}
\end{align}
\end{lemma}
\begin{proof}
$(i)$ We will show that the set
\begin{equation}
\bbR\backslash{\ol A}^e=\{x\in\bbR\,|\, \text{there is an
$\varepsilon_0>0$ such that } |(x-\varepsilon_0,x+\varepsilon_0)\cap
A|=0\}
\lb{A.11}
\end{equation}
is open. Pick $x_0\in\bbR\backslash{\ol A}^e$, then there is an
$\varepsilon_0>0$ such that $|(x_0-\varepsilon_0,x+\varepsilon_0)\cap
A|=0$. Consider $x_1\in(x_0-(\varepsilon_0/4),x_0+(\varepsilon/4))$ and
the open ball $S(x_1;\varepsilon_0/4)$ centered at $x_1$ with radius
$\varepsilon_0/4$. Then,
\begin{equation}
|S(x_1;\varepsilon_0/4)\cap A|\leq
|(x_0-\varepsilon_0,x_0+\varepsilon_0)\cap A|=0  \lb{A.12}
\end{equation}
and hence $x_1\in\bbR\backslash{\ol A}^e$ and
$S(x_0;\varepsilon/4)\subseteq \bbR\backslash {\ol A}^e$. Thus,
$\bbR\backslash{\ol A}^e$ is open. \\
$(ii)$ Let $x\in {\ol A}^e$. Then for all $\varepsilon>0$,
$|(x-\varepsilon,x+\varepsilon)\cap A|>0$. Choose $\varepsilon_n=1/n$,
$n\in\bbN$, then $(x-\varepsilon_n,x+\varepsilon)\cap A \neq \emptyset$
and we may choose an $x_n\in (x-\varepsilon_n,x+\varepsilon_n)\cap A$.
Since $x_n\to x$ as $n\to \infty$, $x\in \ol A$ and hence
${\ol A}^e\subseteq \ol A$. \\
$(iii)$ By $(ii)$, ${\ol A}^e\subseteq \ol A$. Hence,
$\ol{\big({\ol A}^e\big)}^e \subseteq \ol{{\ol A}^e}={\ol A}^e$ since
${\ol A}^e$ is closed by $(i)$.
\end{proof}

\begin{lemma} \lb{lA.10}
Let $d\mu=d\mu_{\rm ac}$ and $S_1$ and $S_2$ be essential supports of $d\mu$.
Then,
\begin{equation}
{\ol S_1}^e = {\ol{S_2}}^e.  \lb{A.13}
\end{equation}
\end{lemma}
\begin{proof}
Since $|S_1\backslash S_2|=|S_2\backslash S_1|=0$, \eqref{A.13} follows
from \eqref{A.7} and \eqref{A.10}.
\end{proof}

Actually, one can improve on Lemma \ref{lA.10}:

\begin{lemma} \lb{lA.11}
Let $d\mu=d\mu_{\rm ac}=fdx$, $0\leq f\in L^1_{\rm loc}(\bbR)$. If $S$ is
any essential support of $d\mu$, then,
\begin{equation}
{\ol S}^e = \ol{\{x\in\bbR\,|\, f(x)>0\}}^e=\supp \, (d\mu).  \lb{A.14}
\end{equation}
\end{lemma}
\begin{proof}
Since $\{x\in\bbR\,|\, f(x)>0\}$ is an essential support of $d\mu$, it
suffices to prove
\begin{equation}
\ol{\{x\in\bbR\,|\, f(x)>0\}}^e=\supp \, (d\mu).  \lb{A.15}
\end{equation}
We denote $U=\bbR\backslash \supp \, (d\mu)$. Then $U$ is the largest
open set that satisfies $\mu(U)=0$. Next, let
$U'=\bbR\backslash\ol{\{x\in\bbR\,|\, f(x)>0\}}^e$. By Lemma
\ref{lA.9}\,$(i)$, $U'$ is open. \\
``$\supseteq$'': Let $x\in U'$. Then there is an $\varepsilon_0>0$ such
that
\begin{equation}
|(x-\varepsilon_0,x+\varepsilon_0)\cap \{y\in\bbR\,|\, f(y)>0\}|=0
\lb{A.16}
\end{equation}
(cf.\ \eqref{A.11}). Hence,
\begin{equation}
f=0 \; |\cdot|\text{-a.e.\ on $(x-\varepsilon_0,x+\varepsilon_0)$}
\lb{A.16a}
\end{equation}
and thus, $\mu((x-\varepsilon_0,x+\varepsilon_0))=0$. Next one covers $U'$
with open intervals of this form to arrive at $\mu(U')=0$. Since $U$ is
the largest open set satisfying $\mu(U)=0$, one infers $U'\subseteq U$
and hence
\begin{equation}
\ol {\{x\in\bbR\,|\, f(x)>0\}}^e \supseteq \supp \, (d\mu). \lb{A.17}
\end{equation}
``$\subseteq$'': Suppose $x\in U$. Since $\supp \, (d\mu)$ is closed,
there is an $\varepsilon_0>0$ such that
$(x-\varepsilon_0,x+\varepsilon_0)\cap
\supp \, (d\mu)=\emptyset$. Thus, $\mu((x-\varepsilon,x+\varepsilon))=0$
for all $0\leq\varepsilon\leq\varepsilon_0$. Actually, $\mu(B)=0$ for all
$\mu$-measurable $B\subseteq (x-\varepsilon_0,x+\varepsilon_0)$ and hence
$f=0$ $|\cdot|$-a.e.\ on $(x-\varepsilon_0,x+\varepsilon_0)$. Thus,
\begin{equation}
|(x-\varepsilon,x+\varepsilon)\cap \{y\in\bbR\,|\, f(y)>0\}|=0 \,
\text{ for all $0\leq\varepsilon\leq\varepsilon_0$}  \lb{A.18}
\end{equation}
and one obtains $x\in U'$. Thus, $U'\supseteq U$ and hence
\begin{equation}
\ol {\{x\in\bbR\,|\, f(x)>0\}}^e \subseteq \supp \, (d\mu). \lb{A.19}
\end{equation}
\end{proof}

We remark that a result of the type \eqref{A.14} has been noted
in \cite[Corollary 11.11]{Bu97} in the context of general ordinary
differential operators and their associated Weyl--Titchmarsh matrices.
In this connection we also refer to \cite[p.\ 301]{Te00} for a
corresponding result in connection with Herglotz functions and their
associated measures.

\section{Herglotz Functions and Weyl--Titchmarsh
Theory for \\ Schr\"odinger Operators in a Nutshell} \lb{B}
\renewcommand{\theequation}{B.\arabic{equation}}
\renewcommand{\thetheorem}{B.\arabic{theorem}}
\setcounter{theorem}{0}
\setcounter{equation}{0}

The material in this appendix is well-known, but since we use it at
various places in this paper, we thought it worthwhile to collect it
in an appendix.

\begin{definition} \lb{dB.1}
Let $\bbC_\pm=\{z\in\bbC \mid \Im(z)\gtrless 0 \}$.
$m:{\mathbb{C_+}}\to {\mathbb{C}}$ is called a {\it Herglotz}
function (or {\it Nevanlinna} or {\it Pick} function) if $m$ is analytic
on
${\mathbb{C}}_+$ and
$m({\mathbb{C}}_+)\subseteq {\mathbb{C}}_+$.
\end{definition}

\smallskip

\noindent One then extends $m$ to $\bbC_-$ by reflection, that is, one
defines
\begin{equation}
m(z)=\overline{m(\overline z)},
\quad z\in{\mathbb{C}}_-. \lb{B.1}
\end{equation}
Of course, generally, \eqref{B.1} does not represent the analytic
continuation of $m\big|_{\bbC_+}$ into $\bbC_-$.

The fundamental result on Herglotz functions and their representations
on Borel transforms, in part due to Fatou, Herglotz, Luzin,
Nevanlinna, Plessner, Privalov, de la Vall{\'e}e Poussin, Riesz, and
others, then reads as follows.

\begin{theorem} $($\cite{AG81}, Sect.\ 69, \cite{AD56},
\cite{Do74}, Chs.~II, IV,
\cite{KK74}, \cite{Ko98}, Ch.~6, \cite{Pr56}, Chs.~II,
IV, \cite{RR94}, Ch.~5$)$. \label{tB.2}  \\
Let $m$ be a Herglotz function. Then, \\
$(i)$ $m(z)$ has finite normal limits $m(\lambda
\pm i0)=\lim_{\varepsilon
\downarrow 0} m(\lambda \pm i\varepsilon)$ for
a.e.~$\lambda \in {\mathbb{R}}$. \\
$(ii)$ Suppose $m(z)$ has a zero normal limit on a subset of
${\mathbb{R}}$ having positive Lebesgue measure. Then $m\equiv 0$. \\
$(iii)$ There exists a nonnegative measure $d\omega$ on
${\mathbb{R}}$
satisfying
\begin{equation} \lb{B.2}
\int_{{\mathbb{R}}} \frac{d\omega (\lambda )}{1+\lambda^2}<\infty
\end{equation}
such that the Nevanlinna, respectively, Riesz-Herglotz
representation
\begin{align}
\begin{split}
&m(z)=c+dz+\int_{{\mathbb{R}}}
d\omega (\lambda) \bigg(\frac{1}{\lambda -z}-\frac{\lambda}
{1+\lambda^2}\bigg), \quad z\in\bbC_+, \lb{B.3} \\[2mm]
& \, c=\Re[m(i)],\quad d=\lim_{\eta \uparrow
\infty}m(i\eta )/(i\eta ) \geq 0
\end{split}
\end{align}
holds. Conversely, any function $m$ of the type \eqref{B.3} is a
Herglotz function. \\
$(iv)$ Let $(\lambda _1,\lambda_2)\subset
{\mathbb{R}}$, then the
Stieltjes inversion formula for $d\omega$ reads
\begin{equation} \lb{B.4}
\frac{1}{2}\omega\left(\left\{\lambda_1
\right\}\right)+\frac{1}{2}\omega
\left(\left\{\lambda_2\right\}\right) +\omega
((\lambda_1,\lambda_2))=\pi^{-1}\lim_{\varepsilon\downarrow
0}\int^{\lambda_2}_{\lambda_1}d\lambda \, \Im(m(\lambda
+i\varepsilon)).
\end{equation}
$(v)$ The absolutely continuous $({\it ac})$ part $d\omega_{ac}$ of
$d\omega$ with respect to Lebesgue measure $d\lambda$ on
${\mathbb{R}}$ is
given by
\begin{equation}\lb{B.5}
d\omega_{ac}(\lambda)=\pi^{-1}\Im[m(\lambda+i0)]\,d\lambda.
\end{equation}
$(vi)$ Local singularities of $m$ and $m^{-1}$ are necessarily
real and at most of first order in the sense that
\begin{align}
& \lim_{\epsilon\downarrow0}
(-i\epsilon)\, m(\lambda+i\epsilon) \geq 0, \quad
\lambda\in\bbR, \lb{B.6} \\
& \lim_{\epsilon\downarrow0}
(i\epsilon) \, m(\lambda+i\epsilon)^{-1}
\geq 0, \quad \lambda\in\bbR. \lb{B.7}
\end{align}
\end{theorem}

Next, we denote by
\begin{equation}
d\omega =d\omega_{\rm ac}+d\omega_{\rm sc}
+d\omega_{\rm pp} \lb{B.8}
\end{equation}
the decomposition of $d\omega$ into its absolutely continuous
$({\it ac})$, singularly continuous $({\it sc})$, and pure point $({\it
pp})$ parts with respect to Lebesgue measure on ${\mathbb{R}}$.

\begin{theorem} [\cite{AD56}, \cite{KK74}, \cite{Si95a},
\cite{Si95}]  \label{tB.3}
Let $m$ be a Herglotz function with representation
\eqref{B.3}. Then, \\
$(i)$
\begin{align}
\begin{split}
& d=0 \text{ and } \int_{{\mathbb{R}}} \f{d\omega(\lambda)}
{1+|\lambda|^s} <\infty \text{ for some } s\in (0,2) \\
& \text{if and only if }
\int^\infty_1d\eta \, \eta^{-s} \, \Im [m(i\eta )]<\infty.
\lb{B.8a}
\end{split}
\end{align}
$(ii)$ Let $(\lambda_1,\lambda_2)\subset {\mathbb{R}}$,
$\eta_1>0$.
Then there is a constant
$C(\lambda_1,\lambda_2,\eta_1)>0$ such that
\begin{equation} \lb{B.8b}
\eta|m(\lambda+i\eta)|\leq
C(\lambda_1,\lambda_2,\eta_1),\quad (\lambda,\eta)\in
[\lambda_1,\lambda_2]\times (0,\eta_1).
\end{equation}
$(iii)$
\begin{equation} \lb{B.8c}
\sup_{\eta
>0}\eta |m(i\eta)|<\infty \text{ if and only if }
m(z)=\int_{{\mathbb{R}}} \f{d\omega (\lambda)}{\lambda-z}
\text{ and }
\int_{{\mathbb{R}}}d\omega(\lambda)<\infty.
\end{equation}
In this case,
\begin{equation} \lb{B.8d}
\int_{{\mathbb{R}}}d\omega
(\lambda)=\sup_{\eta >0}\eta \, |m(i\eta )|=-i\lim_{\eta \uparrow
\infty}\eta \, m(i\eta).
\end{equation}
$(iv)$ For all $\lambda \in {\mathbb{R}}$,
\begin{align} \lb{B.8e}
&\lim_{\varepsilon\downarrow 0}\varepsilon
\Re[m(\lambda+i\varepsilon )]=0,\\
&\omega(\{\lambda\})=\lim_{\varepsilon\downarrow 0}\varepsilon \,
\Im[m(\lambda+i\varepsilon)]=-i\lim_{\varepsilon \downarrow 0}
\varepsilon \,
m(\lambda+i\varepsilon ). \lb{B.8f}
\end{align}
$(v)$ Let $L>0$ and suppose $0\leq \Im [m(z)]\leq L$
for all $z\in
{\mathbb{C}}_+$. Then $d=0$, $d\omega$ is purely absolutely
continuous, $d\omega = d\omega_{ac}$, and
\begin{equation} \lb{B.8g}
0\leq
\frac{d\omega(\lambda)}{d\lambda}=\pi^{-1}\lim_{\varepsilon
\downarrow 0}\Im
[m(\lambda+i\varepsilon)]\leq \pi^{-1}L \text{ for a.e. }
\lambda\in
{\mathbb{R}}.
\end{equation}
$(vi)$ Let $p\in (1,\infty)$, $[\lambda_3,\lambda_4]
\subset
(\lambda_1,\lambda_2)$, $[\lambda_1,\lambda_2]\subset
(\lambda_5,\lambda_6)$. If
\begin{equation} \lb{B.8h}
\sup_{0<\varepsilon <1}\int^{\lambda_2}_{\lambda_1}d\lambda
\, |\Im [m(\lambda+i\varepsilon)]|^p < \infty,
\end{equation}
then $d\omega=d\omega _{ac}$ is purely
absolutely continuous on $(\lambda_1,\lambda_2)$,
$\frac{d\omega_{ac}}{d\lambda}\in L^p((\lambda_1,
\lambda_2);d\lambda)$, and
\begin{equation} \lb{B.8i}
\lim_{\varepsilon\downarrow 0}\bigg\|\pi^{-1}\Im
[m(\cdot +i\varepsilon )]-\frac{d\omega
_{ac}}{d\lambda}\bigg\|_{L^p((\lambda_3,\lambda_4);d\lambda)}=0.
\end{equation}
Conversely, if $d\omega$ is purely absolutely continuous on
$(\lambda_5,\lambda_6)$, and if
$\frac{d\omega_{ac}}{d\lambda}\in$ \linebreak
$L^p((\lambda_5,\lambda_6);d\lambda)$,
then \eqref{B.8h} holds. \\
$(vii)$ Let $(\lambda_1,\lambda_2)
\subset{\mathbb{R}}$.
Then a local version
of Wiener's theorem reads for $p\in (1,\infty)$,
\begin{align}
\begin{split}
&\lim_{\varepsilon\downarrow 0}\varepsilon^{p-1}
\int^{\lambda_2}_{\lambda_1}d\lambda \, |\Im
[m(\lambda+i\varepsilon )]|^p  \\
& \quad =\frac{\Gamma (\frac{1}{2})\Gamma (p-\frac{1}{2})}
{\Gamma (p)} \bigg[ \frac{1}{2}\omega (\{\lambda_1\})^p
+\frac{1}{2}\omega (\{\lambda_2\})^p+
\sum_{\lambda\in (\lambda_1,\lambda_2)}
\omega (\{\lambda \})^p
\bigg]. \lb{B.8j}
\end{split}
\end{align}
Moreover, for $0<p<1$,
\begin{equation} \lb{B.8k}
\lim_{\varepsilon \downarrow
0}\int^{\lambda_2}_{\lambda_1}d\lambda \, |\pi^{-1}\Im
[m(\lambda + i\varepsilon)]|^p
=\int^{\lambda_2}_{\lambda_1}d\lambda
\left|\frac{d\omega_{ac}(\lambda)}{d\lambda}\right|^p.
\end{equation}
\end{theorem}

Together with $m$, $\ln (m)$ is a Herglotz function. Moreover, since
\begin{equation} \lb{B.8l}
0\leq \Im [\ln (m(z)]=\arg [m(z)]\leq
\pi,\quad z\in {\mathbb{C}}_+,
\end{equation}
the measure $d\hatt{\omega}$ in the representation \eqref{B.3} of
$\ln(m)$, that is, in the exponential representation of $m$, is purely
absolutely continuous by Theorem \ref{tB.3}\,$(v)$,
$d\hatt{\omega}(\lambda)=\xi (\lambda )d\lambda$ for some $0\leq\xi\leq
1$. These exponential representations have been studied in great detail
by Aronszajn and Donoghue \cite{AD56}, \cite{AD64} and we record a few of
their properties below.

\begin{theorem}
[\cite{AD56}, \cite{AD64}] \label{tB.4}
Suppose $m(z)$ is a Herglotz function with representation
\eqref{B.3}. Then \\
$(i)$ There exists a $\xi \in L^\infty
({\mathbb{R}})$, $0\leq \xi \leq 1$ a.e., such that
\begin{align}
\begin{split}
& \ln(m(z))=k+\int_{{\mathbb{R}}}d\lambda \,
\xi(\lambda) \bigg(\f{1}{\lambda-z}-\f{\lambda}
{1+\lambda^2}\bigg), \quad z\in\bbC_+, \lb{B.8m} \\
& k=\Re [\ln(m(i))],
\end{split}
\end{align}
where
\begin{equation}\lb{B.8n}
\xi (\lambda)=\pi^{-1}\lim_{\varepsilon
\downarrow 0}\Im [\ln (m(\lambda+i\varepsilon))]
\text{ a.e.}
\end{equation}
$(ii)$ Let $\ell_1,\ell_2\in{\mathbb{N}}$ and
$d=0$ in \eqref{B.3}. Then
\begin{align}
&\int^0_{-\infty} d\lambda \, \xi (\lambda )\f{|\lambda
|^{\ell_1}}{1+\lambda^2}+\int^\infty_0d\lambda \,
\xi (\lambda) \f{|\lambda
|^{\ell_2}}{1+\lambda^2} <\infty  \no \\
& \text{if and only if }  \int^0_{-\infty}d\omega
(\lambda)\f{|\lambda|^{\ell_1}}{1+\lambda^2}
+\int^\infty_0d\omega
(\lambda)\f{|\lambda|^{\ell_2}}{1+\lambda^2} <\infty  \lb{B.8o} \\
& \text{and } \lim_{z\to
i\infty}m(z)=c-\int_{{\mathbb{R}}}d\omega (\lambda) \, \f{\lambda}
{1+\lambda^2}>0. \no
\end{align}
$(iii)$
\begin{align}
& \xi (\lambda)=0 \text{ for } \lambda <0 \text{ if and
only if }
\no \\
& d=0, \quad [0,\infty) \text{ is a support for } \omega \,\,
( \text{i.e., }
\omega ((-\infty,0))=0), \lb{B.8p}  \\
& \int^\infty_0 \f{d\omega(\lambda)}{1+\lambda}<\infty, \text{ and }
c \geq\int^\infty_0 d\omega (\lambda)\f{\lambda}{1+\lambda^2}. \no
\end{align}
In this case
\begin{equation}  \lb{B.8q}
\lim_{\lambda\downarrow
-\infty}m(\lambda)=c-\int^\infty_0 d\omega (\lambda')
\f{\lambda'}{1+{\lambda'}^{2}}
\end{equation}
and
\begin{equation}  \lb{B.8r}
c>\int^\infty_0 d\omega(\lambda ) \f{\lambda}
{1+\lambda ^2} \text{ if and only if }
\int^\infty_0 \f{d\lambda \,\xi(\lambda)}{1+\lambda}<\infty.
\end{equation}
$(iv)$ Let $(\lambda_1,\lambda_2)\subset
{\mathbb{R}}$
and suppose
$0\leq A\leq \xi (\lambda)\leq B\leq 1$ for a.e.
$\lambda \in (\lambda_1,\lambda_2)$ with $(B-A)<1$. Then
$\omega$ is
purely absolutely continuous in
$(\lambda_1,\lambda_2)$ and $\frac{d\omega}{d\lambda}\in
L^p((\lambda_3,\lambda_4);d\lambda)$ for $[\lambda _3,
\lambda_4]\subset
(\lambda_1,\lambda_2)$ and all $p<(B-A)^{-1}$. \\
$(v)$ The measure $\omega$ is purely singular,
$\omega=\omega_s$,
$\omega_{ac}=0$ if and only if $\xi$ equals the characteristic function
of a measurable subset $A\subseteq{\mathbb{R}}$, that is,
$\xi=\chi_A$.
\end{theorem}

While Theorems \ref{tB.3}\,$(v)$,\,$(vi)$ and \ref{tB.4}\,$(iv)$ decsribe
necessary conditions for $d\omega$ to be purely absolutely continuous on
an interval $(\lambda_1,\lambda_2)$, we need to appeal to a stronger
result in Section \ref{s4}. Next, we recall some basic facts from
\cite[Chs.\ 10, 11]{Du70}, \cite[Ch.\ II]{Ga81}, \cite[Ch.\ VI]{Ko98},
\cite[Ch.\ 5]{RR94} to prepare the ground for this result.

\begin{definition}  \lb{dB4a}
Let $F\colon\bbC_+\to\bbC$ be analytic and $C\in\bbC$ a constant satisfying
$|C|=1$. \\
$(i)$ $F$ is called an {\it outer function} on $\bbC_+$ if
\begin{equation}
F(z)=C\exp\bigg(\f{1}{\pi i} \int_{\bbR} dx\, \log (K(x))
\bigg(\f{1}{x-z}-\f{x}{1+x^2}\bigg)\bigg), \quad z\in\bbC_+, \lb{B.8z}
\end{equation}
where $K>0$ a.e.\ on $\bbR$ and $\int_{\bbR} dx\,
|\log(K(x))|(1+x^2)^{-1}<\infty$. In this case, $K(x)=|F(x+i0)|$ for a.e.
$x\in\bbR$. \\
$(ii)$ $F$ is called a {\it Blaschke product} if
\begin{equation}
F(z)=C\bigg(\f{z-i}{z+i}\bigg)^n \prod_{j\in J} \f{|z_j^2+1|}{z_j^2+1}
\f{z-z_j}{z-\ol{z_j}}\, , \quad z\in\bbC_+, \;
\{z_j\}_{j\in J}\in\bbC_+\backslash\{i\},   \lb{B.8A}
\end{equation}
where $n\in\bbN_0$, and $J\subseteq\bbN$ is a (possibly empty) index set.
One then has
\begin{equation}
\sum_{j\in J} \f{y_j}{x_j^2+(y_j+1)^2}<\infty, \quad z_j=x_j+y_j, \,
j\in J. \lb{B.8B}
\end{equation}
$(iii)$ $F$ is called a {\it singular inner function} on $\bbC_+$ if
\begin{equation}
F(z)=Ce^{iaz}\exp\bigg(i\int_{\bbR} d\mu_{\rm s}(x)\,
\bigg(\f{1}{x-z}-\f{x}{1+x^2}\bigg)\bigg), \quad z\in\bbC_+, \lb{B.8C}
\end{equation}
where $a\in(0,\infty)$ and $d\mu_{\rm s}$ is a nonnegative singular measure
on $\bbR$ satisfying $\int_{\bbR} d\mu_{\rm s}(x) (1+x^2)^{-1}<\infty$. \\
$(iv)$ $F$ is called an {\it inner function} if $F(z)=B_F(z)S_F(z)$,
$z\in\bbC_+$, where $B_F$ is a Blaschke product and $S_F$ is a singular inner
function. In this case, the factors $B_F$ and $S_f$ are unique up to
multiplicative unimodular constants.
\end{definition}

We also briefly recall the Nevanlinna, Smirnov, and Hardy classes associated
with $\bbC_+$:

\begin{definition}  \lb{dB4c}
$(i)$ The {\it Nevanlinna class} $N(\bbC_+)$ is defined as the union of
the identically vanishing function on $\bbC_+$ and the set of
analytic functions $F\colon\bbC_+\to\bbC$, $F\not\equiv 0$, such that
\begin{equation}
F(z)=I_F(z)O_F(z)/S_F(z), \quad z\in\bbC_+, \lb{B.8D}
\end{equation}
where $I_F$ is an inner, $O_F$ is an outer, and $S_F$ is a singular inner
function on $\bbC_+$. \\
$(ii)$ The {\it Smirnov class} $N_+(\bbC_+)$ is defined as the union of
the identically vanishing function on $\bbC_+$ and the set of analytic
functions $F\colon\bbC_+\to\bbC$, $F\not\equiv 0$, such that
\begin{equation}
F(z)=I_F(z)O_F(z), \quad z\in\bbC_+, \lb{B.8Da}
\end{equation}
where $I_F$ is an inner and $O_F$ is an outer function on $\bbC_+$. \\
$(iii)$ The Hardy spaces $H^p(\bbC_+)$, $p\in(0,\infty)$ are defined by
\begin{equation}
H^p(\bbC_+)=\bigg\{F\colon\bbC_+\to\bbC \, \text{analytic} \,|\,
\|F\|^p_{H^P(\bbC_+)}=
\sup_{y>0}\int_{\bbR} dx \, |F(x+iy)|^p <\infty\bigg\}.  \lb{B.8u}
\end{equation}
$(iv)$ The Hardy spaces $H^p(\bbR)$, $p\in(0,\infty)$ are defined by
\begin{align}
&H^p(\bbR)=\Big\{f\in L^p(\bbR;dx)\,\Big|\, \text{there exists an
$F\in H^p(\bbC_+)$ such that} \no \\
& \hspace*{2.6cm}
\text{for a.e.\ $x\in\bbR$, $\lim_{z\to x}F(z)=f(x)$ nontangentially}\Big\}
    \lb{B.8F}  \\
& \text{with $\|f\|^p_{H^p(\bbR)}=\|f\|^p_{L^p(\bbR;dx)}$.}  \no
\end{align}
\end{definition}

Of course, $\|\cdot\|_{H^P(\bbC_+)}$ and $\|\cdot\|_{H^p(\bbR)}$ are
norms only for $p\geq 1$.

Moreover, $F\in N(\bbC_+)$ if and only if $F$ is of the form $F=G/H$, where
$G$ and $H$ are analytic and bounded on $\bbC_+$, and $H$ is nonvanishing on
$\bbC_+$. In addition, $F\in N_+(\bbC_+)$ if and only if $F$ is of the form
$F=H_F/O_F$, where $H_F$ and $O_F$ are analytic and bounded by $1$ on
$\bbC_+$ and $O_F$ is outer on $\bbC_+$. In particular,
$N_+(\bbC_+)$ is the smallest algebra that contains all inner and outer
functions on $\bbC_+$. One has the (strict) inclusions
\begin{equation}
H^p(\bbC_+)\subset N_+(\bbC_+)\subset N(\bbC_+), \quad p>0.   \lb{B.8E}
\end{equation}
Functions in $N(\bbC_+)$ which are not identically vanishing have a.e.\
nontangential boundary values on $\bbR$ which cannot vanish on a set of
positive Lebesgue measure.

Analogous facts hold with $\bbC_+$ replaced by $\bbC_-$.

One verifies that every Herglotz function $F$ is an outer function on
$\bbC_+$ and hence sums and differences of Herglotz functions lie in
$N_+(\bbC_+)$.

For subsequent purpose in Theorem \ref{tB.4d} below we also need to make the
connection with the real Hardy space $\wti {H^1}(\bbR)$ in the sense of
Stein and Weiss defined as follows (see, e.g., \cite[Sect.\ III.5.12]{St93}):
\begin{align}
\begin{split}
\wti{H^1}(\bbR)&=\{f\in L^1(\bbR) \,|\, \text{there exist $F_j\in
H^p(\bbC_+)$ such that $f=f_1+\ol{f_2}$,} \\
& \qquad
\text{where $f_j(x)=\lim_{\varepsilon\downarrow 0} F_j(x+i\varepsilon)$ for
a.e.\  $x\in\bbR$, $j=1,2$}\}.
\end{split}
\end{align}
Moreover, denoting by $\cH$ the Hilbert transform on $L^1(\bbR;dx)$,
\begin{equation}
(\cH f)(x)= \f{1}{\pi} \lim_{\varepsilon\downarrow 0}
\int_{|t-x|\geq \varepsilon} dt \, \f{f(t)}{t-x} \, \text{ for a.e.\
$x\in\bbR$}, \; f\in L^1(\bbR;dx),
\end{equation}
then
\begin{equation}
f\in \wti{H^1}(\bbR) \, \text{ if and only if $f\in L^1(\bbR;dx)$ and
$\cH f\in L^1(\bbR;dx)$}
\end{equation}
and
\begin{equation}
f\in \wti{H^1}(\bbR) \, \text{ implies }  \, \int_{\bbR} dx\, f(x) =0.
\end{equation}

Next, we mention a result due to Zinsmeister \cite{Zi89} which goes beyond
Theorems \ref{tB.3}\,$(v)$, $(vi)$ and \ref{tB.4}\,$(iv)$. For this purpose
we now introduce the space $\wti{H^1}(K)$, where $K\subset\bbR$ is compact
and of positive Lebesgue measure, $|K|>0$, as follows: Let $\cM_{\bbC}(\bbR)$
denote the space of complex-valued (and hence finite) Borel measures on
$\bbR$. Denote by $d\nu\in\cM_{\bbC}$ a complex-valued measure supported on
$K$ normalized by
\begin{equation}
\int_K d\nu(x)=0  \lb{B.8G}
\end{equation}
and introduce the function
\begin{equation}
F(z)=\int_K \f{d\nu(x)}{x-z}, \quad z\in \bbC\backslash K. \lb{B.8H}
\end{equation}
One then verifies that
\begin{equation}
F|_{\bbC_{\pm}}\in H^p(\bbC_\pm) \, \text{ for all $p\in(1/2,1)$}  \lb{B.8I}
\end{equation}
and, of course, $F$ has nontangential limits
\begin{equation}
F_{\pm}(x)=\lim_{z\in\bbC_{\pm}, \, z\to x} F(z) \,
\text{ for a.e.\ $x\in K$.}  \lb{B.8J}
\end{equation}
The space $\wti{H^1}(K)$ is then defined by
\begin{align}
\begin{split}
& \wti{H^1}(K)=\bigg\{d\nu\in\cM_{\bbC}(\bbR) \,\bigg|\,
\supp(d\chi)\subseteq K; \, \int_K d\chi(x)=0;    \\
& \hspace*{1.8cm} \|d\nu\|_{\wti{H^1}(K)} =
\|F_+\|_{L^1(K;dx)}+\|F_-\|_{L^1(K;dx)} <\infty \bigg\}.  \lb{B.8K}
\end{split}
\end{align}

\begin{theorem} [\cite{Zi89}, Theorem\ 1]  \label{tB.4d}
Let $K\subset\bbR$ be compact and of positive Lebesgue measure, $|K|>0$.
Then the following two items are equivalent: \\
$(i)$ $K$ is homogeneous. \\
$(ii)$ Every element $d\nu\in \wti{H^1}(K)$ is of the form
\begin{align}
& d\nu(x)=f(x)dx, \, \text{ where } \, f\in \wti{H^1}(\bbR), \;
\supp(f)\subseteq K;  \lb{B.8L}  \\
& \text{and hence } \,
f\in L^1(\bbR;dx) \, \text{ and } \, \int_K dx \, f(x) =0. \lb{B.8M}
\end{align}
In particular, $d\nu$ is purely absolutely continuous.
\end{theorem}

It is mentioned in \cite{Zi89} that the implication $(i)$ implies $(ii)$ can
be inferred from previous (apparently, unpublished) work by P.\ W.\ Jones.

We return to properties of Herglotz functions.

\begin{theorem} [\cite{Gi84}, \cite{GP87}]  \label{tB.5}
Let $m$ be a Herglotz function with representation \eqref{B.3} and
denote by $\Lambda$ the set
\begin{equation}
\Lambda=\{\lambda\in\bbR\,|\, \Im[m(\lambda+i0)] \, \text{exists
$($finitely or infinitely$)$}\}. \lb{B.9}
\end{equation}
Then, $S$, $S_{\rm ac}$, $S_{\rm s}$, $S_{\rm sc}$, $S_{\rm pp}$ are
essential supports of $d\omega$, $d\omega_{\rm ac}$, $d\omega_{\rm s}$,
$d\omega_{\rm sc}$, $d\omega_{\rm pp}$, respectively, where
\begin{align}
S&=\{\lambda\in\Lambda\,|\, 0<\Im[m(\lambda+i0)]\leq\infty\}, \lb{B.10}
\\
S_{\rm ac}&=\{\lambda\in\Lambda\,|\, 0<\Im[m(\lambda+i0)]<\infty\},
\lb{B.11} \\
S_{\rm s}&=\{\lambda\in\Lambda\,|\, \Im[m(\lambda+i0)]=\infty\},
\lb{B.12} \\
S_{\rm sc}&=\Big\{\lambda\in\Lambda\,|\, \Im[m(\lambda+i0)]=\infty, \,
\lim_{\varepsilon\downarrow 0}
(-i\varepsilon)m(\lambda+i\varepsilon)=0\Big\}, \lb{B.13} \\
S_{\rm pp}&=\Big\{\lambda\in\Lambda\,|\, \Im[m(\lambda+i0)]=\infty, \,
\lim_{\varepsilon\downarrow 0}
(-i\varepsilon)m(\lambda+i\varepsilon)=\omega(\{\lambda\})>0\Big\}.
\lb{B.14}
\end{align}
Moreover, since
\begin{align}
& |\{\lambda\in\bbR\,|\, |m(\lambda+i0)|=\infty, \,
\Im[m(\lambda+i0)]<\infty\}|=0, \lb{B.15}  \\
& \omega(\{\lambda\in\bbR\,|\, |m(\lambda+i0)|=\infty, \,
\Im[m(\lambda+i0)]<\infty\})=0,  \lb{B.16}
\end{align}
also $S_{\rm s}^{\prime}$, $S_{\rm sc}^{\prime}$, $S_{\rm pp}^{\prime}$
are essential supports of $d\omega_{\rm s}$, $d\omega_{\rm sc}$,
$d\omega_{\rm pp}$, respectively, where
\begin{align}
S_{\rm s}^{\prime}&=\{\lambda\in\Lambda\,|\, |m(\lambda+i0)|=\infty\},
\lb{B.17} \\
S_{\rm sc}^{\prime}&=\Big\{\lambda\in\Lambda\,|\,
|m(\lambda+i0)|=\infty, \, \lim_{\varepsilon\downarrow 0}
(-i\varepsilon)m(\lambda+i\varepsilon)=0\Big\}, \lb{B.18} \\
S_{\rm pp}^{\prime}&=\Big\{\lambda\in\Lambda\,|\,
|m(\lambda+i0)|=\infty, \, \lim_{\varepsilon\downarrow 0}
(-i\varepsilon)m(\lambda+i\varepsilon)=\omega(\{\lambda\})>0\Big\}.
\lb{B.19}
\end{align}
In particular,
\begin{equation}
S_{\rm s}\sim S_{\rm s}^{\prime}, \quad S_{\rm sc}\sim S_{\rm
sc}^{\prime}, \quad S_{\rm pp}\sim S_{\rm pp}^{\prime}. \lb{B.20}
\end{equation}
\end{theorem}

Next, consider Herglotz functions $\pm m_\pm$ of the type \eqref{B.3},
\begin{align}
\begin{split}
& \pm m_\pm(z)=c_\pm+d_\pm z+\int_{{\mathbb{R}}}
d\omega_\pm (\lambda) \bigg(\frac{1}{\lambda -z}-\frac{\lambda}
{1+\lambda^2}\bigg), \quad z\in\bbC_+, \lb{B.21} \\
& \,\, c_\pm\in\bbR, \quad d_\pm\geq 0,
\end{split}
\end{align}
and introduce the $2\times 2$ matrix-valued Herglotz function $M$
\begin{align}
&M(z)=\big(M_{j,k}(z)\big)_{j,k=0,1}, \quad z\in\bbC_+, \lb{B.22} \\
&M(z)=\f{1}{m_-(z)-m_+(z)}\begin{pmatrix} 1 & \f{1}{2}[m_-(z)+m_+(z)]
\\
\f{1}{2}[m_-(z)+m_+(z)] & m_-(z)m_+(z) \end{pmatrix}, \quad z\in\bbC_+
\lb{B.23} \\
& \hspace*{.85cm} =C+Dz + \int_{{\mathbb{R}}}
d\Omega (\lambda) \bigg(\frac{1}{\lambda -z}-\frac{\lambda}
{1+\lambda^2}\bigg), \quad z\in\bbC_+,  \lb{B.24} \\
& \, C=C^*, \quad D \geq 0 \no
\end{align}
with $C=(C_{j,k})_{j,k=0,1}$ and $D=(D_{j,k})_{j,k=0,1}$ $2\times 2$
matrices and $d\Omega=(d\Omega_{j,k})_{j,k=0,1}$ a
$2\times 2$ matrix-valued nonnegative measure satisfying
\begin{equation}
\int_\bbR \f{d|\Omega_{j,k}(\lambda)|}{1+\lambda^2} < \infty,
\quad j,k=0,1.  \lb{B.25}
\end{equation}
Moreover, we introduce the trace Herglotz function $M^{\rm tr}$
\begin{align}
&M^{\rm tr}(z)=M_{0,0}(z)+M_{1,1}(z)=\f{1+m_-(z)m_+(z)}{m_-(z)-m_+(z)}
\lb{B.26} \\
& \hspace*{1.05cm} =a+b z + \int_{\bbR} d\Omega^{\rm tr}(\lambda)
\bigg(\frac{1}{\lambda -z}-\frac{\lambda} {1+\lambda^2}\bigg), \quad
z\in\bbC_+, \lb{B.27} \\ & \, a\in\bbR, \;  b\geq 0, \quad d\Omega^{\rm
tr}= d\Omega_{0,0} + d\Omega_{1,1}. \no
\end{align}
Then,
\begin{equation}
d\Omega \ll d\Omega^{\rm tr} \ll d\Omega   \lb{B.28}
\end{equation}
(where $d\mu \ll d\nu$ denotes that $d\mu$ is absolutely
continuous with respect to $d\nu$).

\begin{theorem} [\cite{Gi89}, \cite{Gi98}]  \label{tB.6}
Let $m_\pm$, $M$, and $M^{\rm tr}$ be a Herglotz functions with
representations \eqref{B.21}, \eqref{B.24}, and \eqref{B.27},
respectively, and denote by $\Lambda_\pm$ the sets
\begin{equation}
\Lambda_\pm=\{\lambda\in\bbR\,|\, \Im[m_\pm(\lambda+i0)] \,
\text{exists $($finitely or infinitely$)$}\}. \lb{B.29}
\end{equation}
Then, $S_{\Omega^{\rm tr}}$, $S_{\Omega^{\rm tr},\rm ac}$,
$S_{\Omega^{\rm tr},\rm s}$, $S_{\Omega^{\rm tr},\rm sc}$,
$S_{\Omega^{\rm tr},\rm pp}$ are essential supports of
$d\Omega^{\rm tr}$, $d\Omega_{\rm ac}^{\rm tr}$,
$d\Omega_{\rm s}^{\rm tr}$, $d\Omega_{\rm sc}^{\rm tr}$,
$d\Omega_{\rm pp}^{\rm tr}$, respectively, where
\begin{align}
S_{\Omega^{\rm tr}}&=\bbR\backslash
\{\lambda\in\Lambda_+\cap\Lambda_-\,|\,
m_\pm(\lambda+i0)\in\bbR, \, m_+(\lambda+i0)\neq m_-(\lambda+i0)\},
\lb{B.30} \\
&\quad \sim \{\lambda\in\Lambda_+\cap\Lambda_-\,|\,
0<\Im[m_+(\lambda+i0)]<\infty\} \no \\
& \qquad \cup \{\lambda\in\Lambda_+\cap\Lambda_-\,|\,
0<-\Im[m_-(\lambda+i0)]<\infty\} \no \\
& \qquad \cup \{\lambda\in\Lambda_+\cap\Lambda_-\,|\,
m_+(\lambda+i0)=m_-(\lambda+i0)\in\bbR\} \no \\
& \qquad \cup \{\lambda\in\Lambda_+\cap\Lambda_-\,|\,
|m_+(\lambda+i0)|=|m_-(\lambda+i0)|=\infty\}, \lb{B.31} \\
S_{\Omega^{\rm tr},\rm ac}&=\{\lambda\in\Lambda_+\cap\Lambda_-\,|\,
0<\Im[m_+(\lambda+i0)]<\infty\} \no \\
& \quad \cup \{\lambda\in\Lambda_+\cap\Lambda_-\,|\,
0<-\Im[m_-(\lambda+i0)]<\infty\}, \lb{B.32} \\
S_{\Omega^{\rm tr},\rm s}&=\{\lambda\in\Lambda_+\cap\Lambda_-\,|\,
m_+(\lambda+i0)=m_-(\lambda+i0)\in\bbR\} \no \\
& \quad \cup \{\lambda\in\Lambda_+\cap\Lambda_-\,|\,
|m_+(\lambda+i0)|=|m_-(\lambda+i0)|=\infty\}, \lb{B.33} \\
S_{\Omega^{\rm tr},\rm sc}&=\Big\{\lambda\in S_{\Omega^{\rm tr},\rm s}
\,|\, \lim_{\varepsilon\downarrow 0}
(-i\varepsilon)M^{\rm tr}(\lambda+i\varepsilon)=0\Big\}, \lb{B.34} \\
S_{\Omega^{\rm tr},\rm pp}&=\Big\{\lambda\in S_{\Omega^{\rm tr},\rm s}
\,|\, \lim_{\varepsilon\downarrow 0}
(-i\varepsilon)M^{\rm tr}(\lambda+i\varepsilon)<0\Big\}. \lb{B.35}
\end{align}
\end{theorem}

Theorem \ref{tB.6} can be applied to the Schr\"odinger operator $H$
defined in \eqref{2.1} assuming (for simplicity) that the associated
differential expression $L=-d^2/dx^2+V(x)$ is in the limit point case
at $+\infty$ and $-\infty$ (to avoid defining separated boundary
conditions at $=\infty$ and/or $-\infty$). To this end we now identify
\begin{equation}
m_\pm(z) \, \text{ and } \, m_\pm(z,x_0), \quad z\in\bbC_+,  \lb{B.36}
\end{equation}
where $m_\pm(z,x_0)$ are the half-line Weyl--Titchmarsh $m$-functions
associated with the Schr\"odinger operators $H_{\pm,x_0}$ in
$L^2([x_0,\pm\infty);dx)$ with a Dirichlet boundary condition at
$x_0$ (cf.\ \eqref{2.5}).

One then obtains the following basic result.

\begin{theorem} [\cite{Gi98}, \cite{Ka62}, \cite{Ka63}, \cite{Si05}]
\label{tB.7} ${}$ \\
$(i)$ The spectral multiplicity of $H$ is two if and only if
\begin{equation}
|\cM_2|>0,  \lb{B.37}
\end{equation}
where
\begin{equation}
\cM_2=\{\lambda\in\Lambda_+\,|\,
m_+(\lambda+i0,x_0)\in\bbC\backslash\bbR\}\cap\{\lambda\in\Lambda_-\,|\,
m_-(\lambda+i0,x_0)\in\bbC\backslash\bbR\}.  \lb{B.38}
\end{equation}
If $|\cM_2|=0$, the spectrum of $H$ is simple. Moreover, $\cM_2$ is a
maximal set on which $H$ has uniform multiplicity two. \\
$(ii)$ A maximal set $\cM_1$ on which $H$ has uniform multiplicity one
is given by
\begin{align}
\cM_1&=\{\lambda\in\Lambda_+\cap\Lambda_-\,|\, m_+(\lambda+i0,x_0)=
m_-(\lambda+i0,x_0)\in\bbR\}  \no \\
&\quad \cup \{\lambda\in\Lambda_+\cap\Lambda_-\,|\,
|m_+(\lambda+i0,x_0)|= |m_-(\lambda+i0,x_0)|=\infty\}  \no \\
&\quad \cup \{\lambda\in\Lambda_+\cap\Lambda_-\,|\,
m_+(\lambda+i0,x_0)\in\bbR,
m_-(\lambda+i0,x_0)|\in\bbC\backslash\bbR\}  \no \\
&\quad \cup \{\lambda\in\Lambda_+\cap\Lambda_-\,|\,
m_-(\lambda+i0,x_0)\in\bbR,
m_+(\lambda+i0,x_0)|\in\bbC\backslash\bbR\}.  \lb{B.39}
\end{align}
In particular, $\sigma_{\rm s}(H)=\sigma_{\rm sc}(H)\cup
\sigma_{\rm pp}(H)$ is always simple.
\end{theorem}

\medskip

\noindent {\bf Acknowledgments.}
We are indebted to Konstantin Makarov, Franz Peherstorfer, Misha Sodin, and
Michel Zinsmeister for helpful discussions, and especially thank Konstantin
Makarov for his help in connection with the material in Appendices \ref{A} and
\ref{B}.

Fritz Gesztesy gratefully acknowledges the extraordinary
hospitality of Helge Holden and the Department of Mathematical Sciences of
the Norwegian University of Science and Technology, Trondheim, during two
two-month stays in the summers of 2004 and 2005, where parts of this paper
were written. He also gratefully acknowledges a research leave for the
academic year 2005/06 granted by the Office of Research of the University of
Missouri--Columbia.

Peter Yuditskii wishes to thank Franz Peherstorfer and the Section for Dynamical 
Systems and Approximation Theory, University of Linz, Austria, for the great 
hospitality extended to him from September 2003 to September 2005 and he is particularly grateful to Franz Peherstorfer for creating a stimulating atmosphere and numerous discussions on the topic at hand. He also thanks Stas Kupin and the Department of Mathematics of CMI, Universit\'e de Provence, Marseille, France, for their  kind hospitality during a one-month stay in the summer of 2004, where a part of this
paper was written.



\begin{thebibliography}{99}
%
\bi{AG81} N.~I.~Akhiezer and I.\ M.\ Glazman, {\it Theory of
Operators in Hilbert Space}, Vol. I, Pitman, Boston, 1981.
%
\bi{Ar57} N.~Aronszajn,
{\it On a problem of Weyl in the theory of singular Sturm--Liouville
equations}, Amer. J. Math. {\bf 79}, 597--610 (1957).
%
\bi{AD56} N.~Aronszajn and W.~F.~Donoghue, {\it On exponential
representations of analytic functions in the upper half-plane with
positive imaginary part}, J. Analyse Math. {\bf 5}, 321--388
(1956--57).
%
\bibitem{AD64} N.~Aronszajn and W.~F.~Donoghue, {\it A
supplement to the
paper on exponential representations of analytic functions
in the upper half-plane with positive imaginary parts}, J.
Analyse Math. {\bf 12}, 113--127 (1964).
%
\bi{AS81} J.~Avron and B.~Simon, {\it Almost periodic Schr\"odinger
operators I. Limit periodic potentials}, Commun. Math. Phys.
{\bf 82}, 101--120 (1981).
%
\bi{BBEIM94} E.~D.~Belokolos, A.~I.~Bobenko, V.~Z.~Enol'skii, A.~R.~Its,
and V.~B.~Matveev, {\it Algebro-Geometric Approach to Nonlinear
Integrable {E}quations}, Springer, Berlin, 1994.
%
\bi{Bu97} D.~Buschmann, {\it Spektraltheorie verallgemeinerter
Differentialausdr\"ucke -- Ein neuer Zugang}, Ph.D. Thesis,
University of Frankfurt, Germany, 1997.
%
\bi{Ca12} C.\ Caratheodory, {\it Untersuchungen \"uber die konformen 
Abbildungen von festen und ver\"anderlichen Gebieten}, Math. Ann. {\bf 72}, 
107--144 (1912).
%
\bi{Ca83} L.\ Carleson, {\it On $H^\infty$ in multiply connected
domains}, in {\it Conference on Harmonic Analysis in Honor of Antoni
Zygmund}, Vol.\ II, W.\ Beckner, A.\ P.\ Calder\'on, R.\ Fefferman,
and P.\ W.\ Jones (eds.), Wadsworth, CA, 1983, pp.\ 349--372.
%
\bi{CL90} R.~Carmona and J.~Lacroix, {\it Spectral Theory of
Random Schr\"odinger Operators}, Birkh\"auser, Basel, 1990.
%
\bi{Ch84} V.\ A.\ Chulaevskii, {\it Inverse spectral problem for
limit-periodic Schr\"odinger operators}, Funct. Anal. Appl. {\bf 18},
230--233 (1984).
%
\bi{CG01} S.~Clark and F.~Gesztesy, {\it Weyl--Titchmarsh
$M$-function asymptotics for matrix-valued Schr\"odinger
operators}, Proc. London Math. Soc. {\bf 82}, 701--724 (2001).
%
\bi{CG03} S.~Clark and F.~Gesztesy. On Povzner--Wienholtz-type
self-adjointness results for matrix-valued Sturm--Liouville
operators, Proc.\ Roy.\ Soc.\ Edinburg {\bf 133A}, 747--758 (2003).
%
\bi{CGHL00} S.~Clark, F.~Gesztesy, H.~Holden, and
B.~M.~Levitan, {\it Borg-type theorems for matrix-valued
Schr\"odinger operators}, J. Diff. Eqs. {\bf 167}, 181--210 (2000).
%
\bi{Cr89} W.~Craig, {\it The trace formula for Schr\"odinger
operators on the line}, Commun.~Math.~Phys. {\bf 126},
379--407 (1989).
%
\bi{CFKS87} H.~L.~Cycon, R.~G.~Froese, W.~Kirsch, and B.~Simon, {\it
Schr\"odinger Operators}, Springer, Berlin, 1987.
%
\bi{DS83} P.~Deift and B.~Simon, {\it Almost periodic Schr\"odinger
operators III. The absolutely continuous spectrum in one dimension},
Commun.~Math.~Phys. {\bf 90}, 389--411 (1983).
%
\bi{DSS94} R.~del Rio, B.~Simon, and G.~Stolz, {\it Stability of
spectral types for Sturm-Liouville operators}, Math. Res. Lett.
{\bf 1}, 437--450 (1994).
%
\bibitem{Do74} W.~F.~Donoghue, {\it Monotone Matrix Functions and
Analytic Continuation}, Springer, Berlin, 1974.
%
\bi{DMN76} B.~A.~Dubrovin, V.~B.~Matveev, and S.~P.~Novikov,
{\it Non-linear equations of Korteweg-de Vries type, finite-zone
linear operators, and Abelian varieties}, Russian Math. Surv.
{\bf 31:1}, 59--146 (1976).
%
\bi{Du70} P.\ L.\ Duren, {\it Theory of $H^p$ Spaces}, Academic Press,
New York, 1970.
%
\bi{Eg92} I.\ E.\ Egorova, {\it On a class of almost periodic solutions of
the KdV equation with a nowhere dense spectrum}, Russian Acad. Sci. Dokl.
Math. {\bf 45}, 290--293 (1990).
%
\bi{Ga81} J.\ B.\ Garnett, {\it Bounded Analytic Functions}, Academic
Press, New York, 1981.
%
\bi{GL51} I.~M.~Gelfand and B.~M.~Levitan, {\it On the determination of a
differential equation from its spectral function},  Izv.~Akad.~Nauk
SSR.~Ser.~Mat. {\bf 15}, 309--360 (1951) (Russian); English transl.~in
Amer.~Math.~Soc.~Transl.~Ser. 2 {\bf 1}, 253--304 (1955).
%
\bi{Ge93} F.~Gesztesy, {\em A complete spectral characterization of
the double commutation method}, J. Funct. Anal. {\bf 117}, 401--446
(1993).
%
\bi{GH03} F.\ Gesztesy and H.\ Holden, {\it Soliton Equations and
Their Algebro-Geometric Solutions. Vol. I: $(1+1)$-Dimensional
Continuous Models},  Cambridge Studies in Advanced Mathematics,
Vol.\ 79, Cambridge Univ. Press, 2003.
%
\bi{GKT96} F.~Gesztesy, M.~Krishna, and G.~Teschl, {\it On
isospectral sets of Jacobi operators,} Commun. Math. Phys.
{\bf 181}, 631--645 (1996).
%
\bi{GS96} F.~Gesztesy and B.~Simon, {\it The $\xi$
function}, Acta Math. {\bf 176}, 49--71 (1996).
%
\bi{Gi84} D.~J.~Gilbert, {\it Subordinacy and Spectral Analysis
of Schr\"odinger Operators}, Ph.D. Thesis, University of Hull,
1984.
%
\bi{Gi89} D.~J.~Gilbert, {\it On subordinacy and
analysis of the
spectrum of Schr\"odinger operators with two singular
endpoints}, Proc.
Roy. Soc. Edinburgh {\bf 112A}, 213-229 (1989).
%
\bi{Gi98} D.~J.~Gilbert, {\it On subordinacy and spectral
multiplicity for a class of singular differential operators},
Proc. Roy. Soc. Edinburgh {\bf A 128}, 549--584 (1998).
%
\bi{GP87} D.~J.~Gilbert and D.~B.~Pearson, {\it On
subordinacy and
analysis of the spectrum of one-dimensional Schr\"odinger
operators}, J.
Math. Anal. Appl. {\bf 128}, 30-56 (1987).
%
\bi{Go69} G.\ M.\ Goluzin, {\it Geometric Theory of Functions of a 
Complex Variable}, Amer. Math. Soc., Providence, RI, 1969. 
%
\bi{Ha48} P.~Hartman {\it Differential equations with
non-oscillatory eigenfunctions}, Duke Math. J. {\bf 15}, 697--709
(1948).
%
\bi{Ha83} M.\ Hasumi, {\it Hardy Classes of Infinitely Connected Riemann
Surfaces}, Lecture Notes in Math. {\bf 1027}, Springer, Berlin, 1983.
%
\bi{Jo82} R.~A.~Johnson, {\it The recurrent Hill's
equation}, J. Diff. Eqs. {\bf 46}, 165--193 (1982).
%
\bi{Jo88} R.~A.~Johnson, {\it On the Sato--Segal--Wilson solutions of
the K--dV equation}, Pac. J. Math. {\bf 132}, 343--3355 (1988).
%
\bi{JM82} R.\ Johnson and J.\ Moser, {\it The rotation number for
almost periodic potentials}, Commun. Math. Phys. {\bf 84}, 403--438
(1982).
%
\bi{JM85} P.\ W.\ Jones and D.\ E.\ Marshall, {\it Critical points of
Green's function, harmonic measure, and the corona problem}, Ark. Mat. {\bf
23}, 281--314 (1985).
%
\bi{Ka62} I.~S.~Kac, {\it On the multiplicity of the spectrum of a
second-order differential operator}, Sov. Math. Dokl. {\bf 3},
1035--1039 (1962).
%
\bi{Ka63} I.~S.~Kac, {\it Spectral multiplicity of a second order
differential operator and expansion in eigenfunctions}, Izv. Akad.
Nauk SSSR {\bf 27}, 1081--11112 (1963). Erratum, Izv. Akad.
Nauk SSSR {\bf 28}, 951--952 (1964). (Russian.)
%
\bibitem{KK74} I.~S.~Kac and M.~G.~Krein, {\it $R$-functions--analytic
functions mapping the upper halfplane into itself}, Amer. Math. Soc.
Transl. (2) {\bf 103}, 1-18 (1974).
%
\bibitem{Ko98} P.~Koosis, {\it Introduction to $H_p$ Spaces}, 2nd ed.,
Cambridge Tracts in Mathematics, Cambridge University Press, Cambridge,
1998.
%
\bi{Ko84} S.~Kotani, {\it Ljapunov indices determine absolutely
continuous spectra of stationary random one-dimensional
Schr\"odinger operators}, in ``{\it Stochastic Analysis\/}'',
K.~It{\v o} (ed.), North-Holland, Amsterdam, 1984, pp.~225--247.
%
\bi{Ko87a} S.~Kotani, {\it One-dimensional random Schr\"odinger
operators and Herglotz functions}, in ``{\it Probabilistic
Methods in Mathematical Physics\/}'', K.~It{\v o} and N.~Ikeda
(eds.), Academic  Press, New York, 1987, pp.~219--250.
%
\bi{Ko87b} S.~Kotani, {\it Link between periodic potentials and
random potentials in one-dimensional Schr\"odinger operators}, in
{\it Differential Equations and Mathematical Physics},
I.~W.~Knowles (ed.), Springer, Berlin, 1987, pp.\ 256--269.
%
\bi{KK88} S.~Kotani and M.~Krishna, {\it Almost periodicity of some
random potentials}, J.~Funct.~Anal. {\bf 78}, 390--405 (1988).
%
\bibitem{Le82}  B.~M.~Levitan, {\it Almost periodicity of
infinite-zone potentials}, Math. USSR
Izvestija  {\bf 18}, 249--273 (1982).
%
\bibitem{Le83}  B.~M.~Levitan, {\it Approximation of
infinite-zone potentials by finite-zone potentials},
Math USSR Izvestija {\bf 20}, 55--87 (1983).
%
\bibitem{Le84}  B.~M.~Levitan, {\it An inverse problem for the
Sturm-Liouville operator in the case of finite-zone and
infinite-zone potentials}, Trans. Moscow Math. Soc. {\bf 45:1}, 1--34
(1984).
%
\bi{Le85} B.~M.~Levitan, {\it On the closure of the set of finite-zone
potentials\/}, Math.\ USSR Sbornik {\bf 51}, 67--89 (1985).
%
\bi{Le87} B.~M.~Levitan, {\it Inverse Sturm--Liouville
Problems}, VNU Science  Press, Utrecht, 1987.
%
\bi{LG64} B.~M.~Levitan and M.~G.~Gasymov, {\it Determination of a
differential equation by two of its spectra}, Russ. Math. Surv.
{\bf 19:2}, 1--63 (1964).
%
\bi{LS84} B.~M.~Levitan and A.~V.~Savin, {\it Examples of
Schr\"odinger operators with almost periodic potentials and nowhere
dense absolutely continuous spectrum}, Sov. Math. Dokl. {\bf 29},
541--544 (1984).
%
\bibitem{Ma86} V.~A.~Marchenko, {\it Sturm--Liouville
Operators and Applications}, Birkh\"auser, Basel, 1986.
%
\bi{Mo81} J.~Moser, {\it An example of a Schr\"odinger equation with
almost periodic potential and nowhere dense spectrum}, Comment. Math.
Helvetici {\bf 56}, 198--224 (1981).
%
\bi{Na68} M.~A.~Naimark, {\it Linear Differential Operators, Part II},
F.~Ungar, New York, 1968.
%
\bi{NMPZ84} S.\ Novikov, S.\ V.\ Manakov, L.\ P.\ Pitaevskii, V.\ E.\
Zakharov, {\it Theory of Solitons}, Consultants Bureau, New York, 1984.
%
\bi{PF92} L.\ Pastur and A.\ Figotin, {\it Spectra of Random and
Almost-Periodic Operators}, Springer, Berlin, 1992.
%
\bi{PT89} L.\ A.\ Pastur and V.\ A.\ Tkachenko, {\it Spectral theory of
a class of one-dimensional Schr\"odinger operators with limit-periodic
potentials}, Trans. Moscow Math. Soc. 1989, 115--166.
%
\bi{PY03} F.~Peherstorfer and P.~Yuditskii, {\it Asymptotic behavior
of polynomials orthonormal on a homogeneous set}, J. Analyse Math.
{\bf 89}, 113--154 (2003).
%
\bibitem{Pr56} I.~I.~Priwalow, {\it Randeigenschaften analytischer
Funktionen}, 2nd ed., VEB Verlag, Berlin, 1956.
%
\bibitem{RS78} M.\ Reed and B.\ Simon, {\it Methods of Modern
Mathematical Physics. IV: Analysis of Operators,} Academic Press, New
York, 1978.
%
\bi{Re51} F.~Rellich, {\it Halbbeschr{\"a}nkte gew{\"o}hnliche
Differentialoperatoren zweiter Ordnung}, Math. Ann. {\bf 122}, 343--368
(1951).
%
\bi{Ro60} F.~S.~Rofe-Beketov, {\it Expansions in eigenfunctions of
infinite systems of differential equations in the non-self-adjoint
and self-adjoint cases}, Mat. Sb. {\bf 51}, 293--342 (1960).
(Russian.)
%
\bibitem{RR94} M.~Rosenblum and J.~Rovnyak, {\it Topics in
Hardy Classes and Univalent Functions}, Birkh\"auser, Basel, 1994.
%
\bi{Ry01a} A.~Rybkin, {\it On the trace approach to
the inverse scattering problem in dimension one}, SIAM J. Math.
Anal. {\bf 32}, 1248--1264 (2001).
%
\bibitem{Si95a} B.~Simon, {\it Spectral analysis of rank one
perturbations and applications}, CRM Proceedings and Lecture
Notes {\bf 8}, 109-149 (1995).
%
\bi{Si95} B.~Simon, {\it $L^p$ norms of the Borel transform and the
decomposition of measures}, Proc. Amer. Math. Soc. {\bf 123},
3749--3755 (1995).
%
\bi{Si05} B.\ Simon, {\it On a theorem of Kac and Gilbert}, J. Funct.
Anal. {\bf 223}, 109--115 (2005).
%
\bi{SW86} B.~Simon and T.~Wolff, {\it Singular continuous spectrum
under rank one perturbations and localization for random
Hamiltonians}, Commun. Pure Appl. Math. {\bf 39}, 75--90 (1986).
%
\bi{SY95a} M.~Sodin and P.~Yuditskii, {\it Almost periodic
Sturm-Liouville operators with Cantor homogeneous spectrum and
pseudoextendible Weyl functions}, Russ. Acad. Sci. Dokl. Math. {\bf
50}, 512--515 (1995).
%
\bi{SY95} M.~Sodin and P.~Yuditskii, {\it Almost periodic Sturm-Liouville
operators with Cantor  homogeneous spectrum}, Comment. Math. Helvetici
{\bf 70}, 639--658 (1995).
%
\bi{SY96} M.~Sodin and P.~Yuditskii, {\it Almost-periodic
Sturm-Liouville operators with homogeneous spectrum}, in {\it
Algebraic and Geometric Methods in Mathematical Physics}, A.~Boutel
de Monvel and A.~Marchenko (eds.), Kluwer, 1996, pp.\ 455--462.
%
\bi{SY97} M.~Sodin and P.~Yuditskii, Almost periodic Jacobi matrices
with homogeneous spectrum, infinite dimensional Jacobi inversion, and
Hardy spaces of character-automorphic functions, J. Geom. Anal. {\bf
7}, 387--435 (1997).
%
\bi{St93} E.\ M.\ Stein, {\it  Harmonic Analysis: Real-Varuiable Methods,
Orthogonality, and Oscillatory Integrals}, Princeton University Press,
Princeton, NJ, 1993.
%
\bi{Te00} G.~Teschl, {\it Jacobi Operators and Completely Integrable
Nonlinear Lattices}, Math. Surv. Monographs, Vol.\ 72,
Amer. Math. Soc., Providence, R.I., 2000.
%
\bi{Th79} C.~Thurlow, {\it A generalisation of the inverse spectral
theorem of Levitan and Gasymov}, Proc. Roy. Soc. Edinburgh {\bf 84A},
185--196 (1979).
%
\bi{Wi71} H.\ Widom, {\it The maximum principle for multiple-valued analytic
functions}, Acta Math. {\bf 126}, 63--82 (1971).
%
\bi{Wi71a} H.\ Widom, {\it $H^p$ sections of vector bundles over Riemann
surfaces}, Ann. Math. {\bf 94}, 304--324 (1971).
%
\bi{Zi89} M.\ Zinsmeister, {\it Espaces de Hardy et domaines de Denjoy},
Ark. Mat. {\bf 27}, 363--378 (1989).
%
\end{thebibliography}
\end{document}